\documentclass{article}
\usepackage[a4paper, total={5.5in, 8in}]{geometry}
\usepackage[english]{babel}
\usepackage[T1]{fontenc}
\usepackage[utf8]{inputenc}
\usepackage{amsmath}
\usepackage{amsthm}
\usepackage{amssymb}

\usepackage{algpseudocode}
\usepackage{accents}
\usepackage{mathrsfs}
\usepackage{graphicx,mathtools,tikz,hyperref}
\usepackage{caption}
\usepackage{subcaption}
\usepackage{algorithm}
\usepackage{stackengine}

\usetikzlibrary{positioning}
\PassOptionsToPackage{hyphens}{url}\usepackage{hyperref}

\newtheorem{theorem}{Theorem}
\newtheorem{corollary}{Corollary}[theorem]

\newtheorem{remark}{Remark}

\newcommand*{\PLap}{\ensuremath{\textrm{Lap}_{h_x}^{\textrm{P}}}}
\newcommand*{\NLap}{\ensuremath{\textrm{Lap}_{h_x}^{\textrm{N}}}}
\newcommand*{\Dx}{\ensuremath{D_x}}

\newcommand*{\Dy}{\ensuremath{D_y}}

\newcommand*{\barrhoax}{\ensuremath{\overline{\rho_1^{n+1}}^x}}
\newcommand*{\barrhoay}{\ensuremath{\overline{\rho_{1}^{n+1}}^y}}
\newcommand*{\barrhobx}{\ensuremath{\overline{\rho_2^{n+1}}^x}}
\newcommand*{\barrhoby}{\ensuremath{\overline{\rho_2^{n+1}}^y}}

\providecommand{\keywords}[1]
{
  \small	
  \textbf{\textit{Keywords---}} #1
}

\title{A first-order computational algorithm for reaction-diffusion type equations via primal-dual hybrid gradient method}
\author{Shu Liu$^{\star}$, Siting Liu$^{\star}$,  Stanley Osher$^{\star}$,  Wuchen Li$^{\dagger}$\ \ \\
\textsuperscript{$\star$} Department of Mathematics, University of California, Los Angles\\
{\tt \small  \{shuliu, siting6, sjo\}@math.ucla.edu}\\
\textsuperscript{$\dagger$} Department of Mathematics, University of South Carolina  \\
 \texttt{wuchen@mailbox.sc.edu}\\
}
\date{}

\begin{document}

\maketitle
\begin{abstract}
We propose an easy-to-implement iterative method for resolving the implicit (or semi-implicit) schemes arising in solving reaction-diffusion (RD) type equations. We formulate the nonlinear time implicit scheme as a min-max saddle point problem and then apply the primal-dual hybrid gradient (PDHG) method. Suitable precondition matrices are applied to the PDHG method to accelerate the convergence of algorithms under different circumstances. Furthermore, our method is applicable to various discrete numerical schemes with high flexibility. From various numerical examples tested in this paper, the proposed method converges properly and can efficiently produce numerical solutions with sufficient accuracy. 
\end{abstract}

\section{Introduction}
Reaction-diffusion (RD) equations (systems) have broad applications in many scientific and engineering areas. In material science, the phase-field model is described by typical RD-type equations known as Allen-Cahn \cite{allen1979microscopic} or Cahn-Hilliard equations \cite{cahn1961spinodal}. They are used to model the development of microstructures in a mixture of two or more materials or phases over time; In chemistry, RD systems are used to depict the reaction and diffusion phenomena of chemical species in which a variety of patterns are produced \cite{pearson1993complex, schnakenberg1979simple}; RD systems are also ubiquitous tools in biology: They are widely used for modeling morphogenesis \cite{chou2007numerical}, as well as the evolution of species distribution in ecology system \cite{murray2001mathematical}. In recent years, researchers also found that RD equations are useful in modeling and predicting crimes \cite{crime_RD}. 

Reaction-diffusion (RD) equations are nonlinear parabolic partial differential equations possessing the following general form
\begin{align}
  &  \frac{\partial u(x,t)}{\partial t} = \mathcal{L} u(x,t) + \mathcal{R} u(x,t) \quad \textrm{on } \Omega\subset \mathbb{R}^d,  \label{rd  equ}  \\
  &  \textrm{with initial condition } u(\cdot, 0) = u_0,  \nonumber
\end{align}
where $\mathcal{L}$ is a certain non-positive definite differential operator associated with the diffusion process. For example, $\mathcal{L}$ can be taken as the Laplace operator $\Delta$ or negative biharmonic operator $-\Delta^2$, or more general operators with variable coefficients; $\mathcal{R}$ is a nonlinear operator depicting the reaction process. The RD equation is usually equipped with either the Neumann boundary condition if $\Omega$ is an ordinary region in $\mathbb{R}^d$, or the periodic boundary condition if $\Omega$ is a periodic region $\mathbb{T}^d$. We can also extend the RD equation \eqref{rd  equ} from the 1D function $u$ to multiple dimensional vector function $\boldsymbol{U}$:
\begin{align}
  &  \frac{\partial \boldsymbol{U}(x,t)}{\partial t} = \mathcal{L} \boldsymbol{U}(x,t) + \mathcal{R} \boldsymbol{U}(x,t) \quad \textrm{on } \Omega\subset \mathbb{R}^d,  \label{rd  system}  \\
  &  \textrm{with initial condition } \boldsymbol{U}(\cdot, 0) = \boldsymbol{U}_0.  \nonumber
\end{align}
Equation \eqref{rd  system} can also be equipped with either Neumann or periodic boundary conditions. We also call the equation \eqref{rd  system} reaction-diffusion system.

In recent decades, numerical methods, including finite difference methods \cite{merriman1994motion, eyre1998unconditionally, hundsdorfer2003numerical, li2010unconditionally, ceniceros2013new, christlieb2014high, shen2010numerical, shen2018scalar, shen2019new, liu2021structure, liu2022second, liu2022convergence} and finite element methods \cite{hundsdorfer2003numerical, zhu2009application, fu2023high}, have been developed for computing the reaction-diffusion type equations (systems). Several benchmark problems \cite{jokisaari2017benchmark, church2019high} have also been introduced to verify the proposed methods' effectiveness. 

In order to get rid of the restriction of the Courant–Friedrichs–Lewy (CFL) condition \cite{courant1928partiellen} on small time steps, most of the popular numerical schemes designed for solving reaction-diffusion equations (systems) in the aforementioned works of literature are implicit or semi-implicit. As one uses implicit or semi-implicit schemes for solving RD equations (systems) with nonlinear terms, Newton's method \cite{atkinson1991introduction} is usually needed for solving the series of nonlinear equations arising from time discretization. However, Newton's method encounters several drawbacks that may affect the performance of the proposed numerical scheme, namely,
\begin{itemize}
    \item Newton's method requires the initial guess position to be close enough to the exact solution of the nonlinear equation. Otherwise, Newton's method may diverge.
    \item When solving RD equations on mesh grids by Newton's method, in each iteration, one has to solve a large-scale linear equation involving the Jacobian matrix of a certain nonlinear function. Solving this large-scale linear equation for multiple Newton iterations could be challenging and time-costing.
\end{itemize}

In this paper, we introduce a method based on the Primal-Dual Hybrid Gradient method (PDHG) \cite{zhu2008efficient, chambolle2011first, valkonen2014primal, clason2017primal, jacobs2019solving} for solving the nonlinear updates arising in the time discretization schemes of RD equations (systems) with satisfying speed and accuracy.

We sketch the proposed method as follows. We briefly illustrate the main idea by considering the following fully implicit, semi-discrete scheme of the RD equation at the $t$-th time step:
\begin{equation}
  \frac{u^{t+1}-u^t}{h_t} = \mathcal{L}u^{t+1} + \mathcal{R}u^{t+1}.  \label{semi discrete RD equ}
\end{equation}
In this case, $u^t$ is given, and $h_t>0$ is a stepsize. We need to solve for $u^{t+1}$. Consider the following function $\mathcal{F}$
\begin{equation}
  \mathcal{F}(u) = u - h_t(\mathcal{L}u + \mathcal{R}u) - u^t. 
\end{equation}
The goal is to solve $\mathcal{F}(u)=0$. If we consider the indicator function $\iota$ defined as
\[ \iota(u)=\begin{cases}
  0 \quad u=0\\
  +\infty \quad u\neq 0.
\end{cases}
\]
Then the nonlinear functional equation $\mathcal{F}(u)=0$ is equivalent to the minimization problem
\[\min_{u\in X}~\iota(\mathcal{F}(u)),\]
where $X$ is a certain linear functional space for $u$. Now since $\iota$ can be treated as the Legendre transform of the constant function $0$, i.e., $\iota(u) = \sup_{p\in X^*} ~\{(p, u)\}$ (here $X^*$ denotes the dual space of $X$), we can recast the above minimization problem as a min-max saddle point problem as follows
\begin{equation}
  \min_{u\in X}\max_{p\in X^*} ~ \{(p, ~ \mathcal{F}(u))\}.  \label{abstract min-max}
\end{equation}
We denote $L(u,p)=(p, \mathcal{F}(u))$. To deal with the saddle problem \eqref{abstract min-max}, by leveraging the ideas proposed in the PDHG method, we evolve $p, u$ via the following proximal algorithms with extrapolation on the dual variable $p$.
\begin{align}
   p_{n+1} = & \underset{p\in X^*}{\textrm{argmin}} ~ \left\{\frac{\|p-p_n\|_2^2}{2\tau_p} - L(u_n, p) \right\} = p_n + \tau_p \mathcal{F}(u_n),  \label{abstract PDHG 1}\\
   \widetilde{p}_{n+1} = & p_{n+1} + \omega (p_{n+1} - p_n), \label{abstract PDHG 2}\\
   u_{n+1} = & \underset{u\in X}{\textrm{argmin}} ~ \left\{\frac{\|u - u_n\|_2^2}{2\tau_u} + L(u, \widetilde{p}_{n+1})\right\} = (\textrm{Id} + \tau_u(\widetilde{p}_{n+1}, \partial_u\mathcal{F}))^{-1}(u_n). \label{abstract PDHG 3}
\end{align}
Here the extrapolation coefficient $\omega >0$, $\tau_p, \tau_u$ are time steps used to evolve $u, p$. We should remind the reader to distinguish the PDHG time steps $\tau_p, \tau_u$ from the time step $h_t$ of the reaction-diffusion equation. Recall $u^{t+1}$ as the solution to $\mathcal{F}(u)=0$, if we further assume that $\partial_u\mathcal{F}(u)$, as a linear 
map from $X$ to $X^*$, is injective. Then it is not hard to verify that $u=u^{t+1}, p=0$ is the equilibrium of the above dynamic \eqref{abstract PDHG 1} - \eqref{abstract PDHG 3}. Thus, we may anticipate that, by evolving \eqref{abstract PDHG 1}, \eqref{abstract PDHG 2}, \eqref{abstract PDHG 3}, $u_k, p_k$ could converge to the desired equilibrium point $u^{t+1}, 0$. 

Furthermore, since $\mathcal{F}$ is usually nonlinear, the inversion in \eqref{abstract PDHG 3} cannot be directly evaluated. To mitigate this, we replace $L(u, \widetilde{p}_{k+1})$ by its linearization $\widehat{L}(u, \widetilde{p}_{k+1}) = L(u_k, \widetilde{p}_{k+1})+(\partial_u L(u_k, \widetilde{p}_{k+1}), u-u_k)$ at $u=u_k$. Thus the update of $u_{k+1}$ will have an explicit form
\begin{equation}
    u_{k+1} = u_k - \tau_u (\widetilde{p}_{k+1}, \partial_u\mathcal{F}(u_k)).  \label{abstract PDHG 3 new}
\end{equation}
As a result, we can evolve the discrete-time dynamic \eqref{abstract PDHG 1}, \eqref{abstract PDHG 2}, \eqref{abstract PDHG 3 new} for approximating the solution $u^{t+1}$ of the nonlinear equation $\mathcal{F}(u)=0$, and the explicit updating rules will enable us to deal with large-scale computational problems conveniently and efficiently. This work will mainly focus on applying such a PDHG algorithm to solve various types of reaction-diffusion equations (systems) up to satisfying accuracy and efficiency. Based on the discussion and presentation in this paper, our method may serve as a potential alternative to the widely used Newton-type algorithms for time-implicit updates of reaction-diffusion equations for time-implicit schemes.

It is worth mentioning that instead of designing and analyzing new discretization schemes for RD equations, our paper is mainly devoted to a strategy that can efficiently resolve the ready-made scheme. Thus, in our paper, we will omit most of the discussions on the properties of the numerical scheme but focus more on the implementing details and the effectiveness of the proposed PDHG method.

We clarify that the method is inspired by \cite{liu2022primal}, in which the authors design a similar algorithm for solving multiple types of PDEs accompanied by 1-D examples. This paper will be more specific and focus on computing various 2-D RD equations (systems) with different boundary conditions. 

It is also worth mentioning that introducing damping terms into wave equations to achieve faster stabilization \cite{fahroo1996optimum} shares great similarity with the limiting stepsize version of applying the PDHG method to PDE-solving algorithms. On the other hand, PDHG methods are also utilized in \cite{valkonen2014primal} to solve nonlinear equations with theoretical convergence guarantee under different circumstances. 

Furthermore, people have applied PDHG or first-order methods to compute time-implicit updates of Wasserstein gradient flows \cite{CarrilloCraigWangWei2022_primala} and reaction-diffusions \cite{CarrilloWangWei2023_structure, fu2023high}. Compared to them, the proposed approach work for general non-gradient flow reaction-diffusion equations. 
In recent research \cite{chen2023transformed}, the authors deal with the nonlinear saddle point problems via the transformed PDHG method, with the follow-up research \cite{chen2023accelerated} aiming at solving the nonlinear equations associated with a class of monotone operators. In recent works \cite{zang2020weak, bao2020numerical}, the authors utilize the weak forms of PDEs and deep learning techniques to compute high-dimensional PDEs. The algorithm in \cite{zang2020weak} can be formulated as a min-max saddle point problem and is directly solved by alternative stochastic gradient descent and ascent method. Although the proposed method shares similarities with the aforementioned research. It differs from them in the saddle point problem formulation and the computational scheme.

This paper is organized as follows. In section \ref{PDHG derivation}, we provide a brief introduction to the Primal-Dual Hybrid Gradients method; then we present the details of how we implement the PDHG method to update the finite difference schemes for RD equations (systems). We then provide some existing theoretical results on the convergence of the PDHG method. We demonstrate our numerical examples in section \ref{numerical example}. Our numerical experiments cover well-known Allen-Cahn and Cahn-Hilliard equations; higher-order gradient flow that emerges from functionalized polymer research; RD system known as the Schnakenberg model, which originates from the study of steady chemical patterns; and RD systems involving nonlocal terms depicting the species evolution of wolves and deer. We conclude the work in section \ref{conclude }. Some of the future research directions will also be discussed in section \ref{conclude }. 

\section{PDHG method for reaction-diffusion equation}\label{PDHG derivation}


In recent years, The Primal-Dual Hybrid Gradient (PDHG) method \cite{zhu2008efficient, zhu2009application, chambolle2011first} proves to be an efficient algorithm for solving saddle point problems emerging from imaging. This method is iterative and each of its iterations consists of alternative proximal steps together with a suitable extrapolation. We refer the readers to \cite{chambolle2011first} for further details (both theoretical and experimental) of the method.

\subsection{PDHG method for updating implicit finite difference schemes}
To clearly convey our proposed idea, let us first consider the following reaction-diffusion equation as an illustrative example on 2D periodic region $\Omega = \mathbb{T}^2$. We assume $\Omega$ is square shaped and denote its side length as $L$.
\begin{align}
   \frac{\partial u(x,t)}{\partial t} = \lambda \Delta u(x,t) + f(u(x,t)), \quad u(x,0) = u_0(x). \label{illustrate rd equ}
\end{align}
Here we assume $\lambda$ is a positive constant coefficient; $f:\mathbb{R}\rightarrow\mathbb{R}$ is the nonlinear function depicting the reaction term. Since we assume $\Omega$ to be the periodic region, we use periodic boundary conditions (BC) for equation \eqref{illustrate rd equ}.

Although there are numerous pieces of research on designing numerical schemes for RD equations, to demonstrate how our method works, let us narrow down and focus on the implicit one-step finite difference (FD) scheme. Once we have demonstrated how to implement the method to this implicit scheme, such a method can be easily extended to general numerical schemes.

We discretize the time interval $[0, T]$ into $N_t$ equal subintervals with length $h_t=T/N_t$. Suppose we discretize each side of the region $\Omega$ into $N_x$ subintervals with space stepsize $h_x = L/N_x$, we choose the central difference scheme to discretize the Laplace operator $\Delta$. We denote the discrete Laplace operator with the periodic boundary condition as $\PLap$ which is an $N_x^2\times N_x^2$ block-circulant matrix possessing the following form
\begin{equation*}
  \PLap = \frac{1}{h_x^2}\left[\begin{array}{ccccc}
      L &  I  &   &   &  I \\
      I &  L  & I &   &    \\
        & \ddots & \ddots & \ddots & \\
        &        &   I    &   L   & I \\
      I &        &        &   I   & L
  \end{array}\right]_{N_x \times N_x \textrm{blocks}} \quad L = \left[\begin{array}{ccccc}
      -4 &  1  &   &   &  1 \\
      1 &  -4  & 1 &   &    \\
        & \ddots & \ddots & \ddots & \\
        &        &   1    &   -4   & 1 \\
      1 &        &        &   1    &-4
  \end{array}\right]_{N_x\times N_x} .
\end{equation*}
Here $I$ is the $N_x$ by $N_x$ identity matrix. We denote $U^k\in \mathbb{R}^{N_x^2}$ as the numerical solution of \eqref{illustrate rd equ} at the $k$th time step. We vectorize the $N_x\times N_x$ square array along the column to form the 1D vector $U^k$. That is, for $l = i\cdot N_x + j$ with $1\leq i\leq N_x$ and $1\leq j\leq N_x$, $U^k_l$ is the numerical approximation of $u((j-1)h_x, (i-1)h_x, kh_t)$. 

Now the implicit one-step FD scheme for \eqref{illustrate rd equ} is cast as
\begin{equation}
  U^{k+1} - U^k = h_t(\lambda ~ \PLap U^{k+1} + f(U^{k+1})), \quad U^0 = U_0.\label{illustrate  implicit schm}
\end{equation}
Here $U_0$ denotes the initial condition on mesh grid points. When solving for the numerical solution of \eqref{illustrate rd equ}, one has to sequentially solve a series of nonlinear equations as shown in \eqref{illustrate  implicit schm}. This is the place in which we should apply the PDHG method. Let us denote the function $F:\mathbb{R}^{N_x^2}\rightarrow \mathbb{R}^{N_x^2}$ as
\begin{equation}
F(U) = U - U^k - h_t(\lambda \PLap U + f(U)). \label{F(U) rd equation}
\end{equation}
We want to solve $F(U)=0$. As discussed in the introduction, this is equivalent to minimizing $\iota(F(U))$, which can further be cast as the following min-max saddle point problem
\begin{equation}
  \min_{U\in\mathbb{R}^{N_x^2}}\max_{P\in\mathbb{R}^{N_x^2}} ~\{L(U, P)\}, \label{illustrate saddle }
\end{equation}
where $L$ is defined as $L(U, P) = P^\top F(U)$. As a result, solving the equation $F(U)=0$ finally boils down to the min-max problem \eqref{illustrate saddle }.

Now the PDHG method suggests the following gradient ascent-descent dynamic for solving \eqref{illustrate saddle }.
\begin{align}
  P_{n+1} & = \underset{P\in\mathbb{R}^{N_x^2}}{\textrm{argmin}} \left\{\frac{\|P - P_n\|_2^2}{2\tau_p}+L(U_n, P)\right\} =  P_n + \tau_p F(U_n); \label{illustrate pdhg dynamic 1}\\
  \widetilde{P}_{n+1} & = {P}_{n+1} + \omega (P_{n+1} - P_n); \label{illustrate pdhg dynamic 2}\\
  U_{n+1} & = \underset{U\in\mathbb{R}^{N_x^2}}{\textrm{argmin}} \left\{\frac{\|U - U_n\|_2^2}{2\tau_u}+L(U, \widetilde{P}_{n+1})\right\} = (\textrm{Id} + \tau_u  \nabla_U F(\cdot)^\top \widetilde{P}_{n+1})^{-1} U_n.\label{illustrate pdhg dynamic 3}
\end{align}
Similar to our discussion in the introduction, the third line above involves a nonlinear equation that cannot be directly solved. We thus can replace the term $L(U, \widetilde{P}_{n+1})$ in \eqref{illustrate pdhg dynamic 3} by the linearization $\widehat{L}(U,P)=L(U_n, P)+\nabla_U L(U_n, P)(U-U_n)$, then \eqref{illustrate pdhg dynamic 3} can be explicitly computed as
\begin{equation}
  U_{n+1} = \underset{U\in\mathbb{R}^{N_x^2}}{\textrm{argmin}} \left\{\frac{\|U - U_n\|_2^2}{2\tau_u}+\widehat{L}(U, \widetilde{P}_{n+1})\right\} = U_n - \tau_u  \nabla_U F(U_{n})^\top \widetilde{P}_{n+1}.  \label{illustate pdhg dynamic 3 simplified}
\end{equation}
Let us denote $U_*$ as the solution to $F(U)=0$. It is not hard to tell that $(U_*, 0)$ is a critical point of the functional $L(U, P)$. Furthermore, $(U_*, 0)$ is the equilibrium point of the time-discrete dynamic \eqref{illustrate pdhg dynamic 1}, \eqref{illustrate pdhg dynamic 2}, \eqref{illustate pdhg dynamic 3 simplified}.

To analyze the convergence speed to the equilibrium state $(U_*, 0)$, we first consider the affine case in which $F(U) = AU-b$ with $A$ as an $N_x^2\times N_x^2$ symmetric matrix. We have the following result, similar to the analysis carried out in \cite{liu2022primal}.
\begin{theorem}[Convergence speed in Linear, symmetric case]\label{thm linear convergence}
  We fix the extrapolation coefficient $\omega=2$. Suppose we obtain the sequence $\{(U_n, P_n)\}_{n\geq 0}$ by evolving the PDHG dynamic \eqref{illustrate pdhg dynamic 1}, \eqref{illustrate pdhg dynamic 2}, \eqref{illustate pdhg dynamic 3 simplified} with initial condition $(U_0, P_0)$. Suppose $F(U)=AU$ with $A$ symmetric and non-singular. Denote $\lambda_{\max}$ as the maximum eigenvalue (in absolute value) of $A$, and denote $\kappa$ as the condition number of $A$. Then $\{(U_n, P_n)\}$ will converge to $(U_*, 0)$ linearly if $\tau_u\tau_p\leq \frac{4}{3\lambda^2_{\max}}$. Then the maximum convergence speed is achieved when $\tau_u\tau_p = \frac{\eta_*}{\lambda_{\max}^2}$, with the optimal convergence rate $\gamma_* = \sqrt{1-\frac{\eta_*}{\kappa^2}}$, i.e., we have for any $n\geq 1$, $\|(U_n,P_n)-(U_*,0)\|_2 \leq \gamma_*^n  \|(U_0,P_0)-(U_*,0)\|_2 .$ Here $\eta_*=\eta_*(\kappa)$ is a function of $\kappa$. The range of $\eta_*$ belongs to $[1, \frac{3}{4})$.
\end{theorem}
The proof of the theorem is provided in Appendix \ref{thm pf}. The explicit form of $\eta_*(\kappa)$ and $\gamma_*$ are given in remark \ref{formula of gamma} of Appendix \ref{thm pf}.

The optimal convergence rate $\gamma_*$ will be very close to $1$ if the condition number $\kappa$ is large. Furthermore, one will require $O(\kappa^2)$ iterations for $U_n$ to converge. This can be very expensive when $A$ is a large-scale matrix with a large condition number. For example, we consider the heat equation
\[ \partial_t u(x,t)  =  \Delta u(x,t) \]
with periodic boundary conditions. We apply the one-step implicit scheme to this equation, i.e., we consider solving $(I-\lambda h_t \PLap) U^{n+1} = U^n$ at each time step $n$. Thus $A = I - h_t  \lambda  \PLap $. One can tell that the eigenvalues of $A$ equal $1+4h_t\lambda N^2\sin^2\left(\frac{\pi k}{N}\right)$ for $1\leq k\leq N.$ When $N$ is even, the condition number $\kappa$ of $A$ equals $1+4\lambda N^2h_t$. Since we can get rid of the CFL condition by using the implicit scheme, we can pick $h_t \gg \frac{1}{N^2}$, which leads to $\kappa(A)\gg 1$. The convergence of the primitive PDHG method could be very slow, even for the heat equation. 

As discussed in remark \ref{formula of gamma}, $\gamma_*$ approaches $0$ if the condition number $\kappa $ drops to $1$. Hence, we need to control the condition number of the matrix $A$ to achieve faster convergence speed. This motivates us to introduce the preconditioning technique to the PDHG method. As suggested in both \cite{jacobs2019solving} and \cite{liu2022primal}, we replace the $l^2$ norm used in either $\|U-U_n\|_2$ or $\|P - P_n\|_2$ by the $G-$norm $\|\cdot\|_G$ which is defined as
\[ \|v\|_G = \sqrt{v^\top G v}, \]
with $G$ as a symmetric, positive definite matrix. In this work, we mainly focus on substituting the norm $\|P - P_n\|_2$ with $\|P-P_n\|_G$. The PDHG method involving $G-$norm in its proximal step is sometimes named $G-$prox PDHG \cite{jacobs2019solving}. The dynamic obtained via such $G-$prox PDHG is shown below.
\begin{align}
  P_{n+1} & =  P_n + \tau_p G^{-1}F(U_n); \label{G-prox linearized pdhg dynamic 1}\\
  \widetilde{P}_{n+1} & = P_{n+1} + \omega (P_{n+1} - P_n); \label{G-prox linearized pdhg dynamic 2}\\
  U_{n+1} & = U_n - \tau_u\nabla_UF(U_n)^\top \widetilde{P}_{n+1}.\label{G-prox linearized pdhg dynamic 3}
\end{align}
We should pause here to emphasize to the reader that the above three-line dynamic \eqref{G-prox linearized pdhg dynamic 1}, \eqref{G-prox linearized pdhg dynamic 2}, \eqref{G-prox linearized pdhg dynamic 3} will be the core gadget for our RD equation solver throughout the remaining part of the paper. Before we move on to further details on solving the RD equation \eqref{illustrate rd equ} via $G-$prox PDHG dynamic, let us provide a little more explanation on how \eqref{G-prox linearized pdhg dynamic 1} - \eqref{G-prox linearized pdhg dynamic 3} improve the convergence speed $\gamma_*$. Analogous to Theorem \ref{thm linear convergence}, we have the following corollary for affine function $F$.
\begin{corollary}[Convergence for $G-$prox dynamic]\label{corollary linear convergence}
  Suppose we keep all the assumptions in Theorem \ref{thm linear convergence}, if we further assume that $G$ commutes with $A$, i.e., $GA=AG$, then all the conclusions in Theorem \ref{thm linear convergence} still hold except $\lambda_{\max}$ now denotes the largest (in absolute value) eigenvalue of $A^\top G^{-1}A$, and $\kappa^2$ now is the condition number of $A^\top G^{-1}A$. 
\end{corollary}
It is now clear that if we can find a matrix $G$ that approximates $AA^\top$ well, then $A^\top G^{-1}A$ will be reasonably close to the identity matrix $I$. Thus the condition number of $A^\top G^{-1}A$ will hopefully remain close to $1$. In such cases, by properly choosing the step size $\tau_u, \tau_p$ such that $\tau_u\tau_p$ is close to $1$, we can obtain a rather fast convergence rate $\gamma_*$. 

Up to this stage, although most of our intuition and analysis on $G-$prox PDHG dynamic comes from the case when $F$ is affine, it is natural to extend our treatment to the nonlinear $F(U)$ defined in \eqref{F(U) rd equation}. We may still anticipate the effectiveness of our method since \eqref{F(U) rd equation} can be recast as
\begin{equation*}
  F(U) = (I-\lambda h_t\PLap)U - U^k - h_t f(U),
\end{equation*}
which can be treated as an affine function with a nonlinear perturbation $h_tf(U)$ carrying the small $h_t$ coefficient. We now discuss several details in applying the $G-$prox PDHG dynamic \eqref{G-prox linearized pdhg dynamic 1}-\eqref{G-prox linearized pdhg dynamic 3} to the above $F(U)$. We aim at evolving the following dynamic in order to update the implicit one-step scheme \eqref{illustrate  implicit schm}.
\begin{align}
  P_{n+1} & =  P_n + \tau_p G^{-1}(U_n - \lambda h_t \PLap U_n -h_tf(U_n) - U^k); \label{G-prox RD pdhg dynamic 1}\\
  \widetilde{P}_{n+1} & = P_{n+1} + \omega (P_{n+1} - P_n); \label{G-prox RD pdhg dynamic 2}\\
  U_{n+1} & = U_n - \tau_u(\widetilde{P}_{n+1} - \lambda h_t \PLap\widetilde{P}_{n+1} - h_t f'(U_n)\odot \widetilde{P}_{n+1}).\label{G-prox RD pdhg dynamic 3}
\end{align}

\noindent
\textbf{Initial guess} \quad It is natural to choose the initial value $U_0$ of the dynamic \eqref{G-prox RD pdhg dynamic 1} - \eqref{G-prox RD pdhg dynamic 3} as the computed result from the last time step $k$, i.e., we set $U_0 = U^k$; And we will simply set $P_0 = 0$. A more sophisticated choice for $U_0$ could be the numerical solution at time $k+1$ obtained by a forward Euler scheme or an IMEX scheme \cite{pareschi2000implicit, hundsdorfer2003numerical}.

\vspace{0.3cm}
\noindent
\textbf{Choosing the matrix $G$} \quad   The function $F(U)$ defined in \eqref{F(U) rd equation} is dominated by the affine term $(I-\lambda h_t\PLap)U - U^k$. It is then natural to choose $G = (I-\lambda h_t\PLap)^2$ as the preconditioner matrix. If the nonlinear term $f(U)$ is a highly stiff term, we may also consider absorbing its Jacobian $\nabla_U f(U) = \textrm{diag}(f(U))$ into $G$. So, $G$ can also be chosen as $G=(I-\lambda h_t\PLap -  \textrm{diag}(f(U)))^2$ in such case. However, there might be a trade-off in doing so: if $\textrm{diag}(f(U))$ does not have equal diagonal entries, $I-\lambda h_t\PLap -  \textrm{diag}(f(U))$ cannot be efficiently inverted by Fast Fourier Transform (FFT) or Discrete Cosine Transform (DCT) method. Nevertheless, in most of our experiments, we discover that choosing $G =(I-\lambda h_t\PLap)^2$ is adequate for achieving satisfying convergence speed. For more general reaction-diffusion equations, we discover that the nonlinear function $F(U)$ is usually decomposed as the sum of the linear term $AU$ and the nonlinear term $h_t f(U)$. The linear term $AU$ can be treated as the dominating term of $F(U)$. It is then reasonable to choose $G=AA^\top$ or at least close to $AA^\top$ as a decent preconditioner of our method. This strategy works properly on general RD equations such as the Cahn-Hilliard equation or higher-order equations arising in polymer science. We refer the reader to examples in section \ref{numerical example} for details.

\vspace{0.3cm}

\noindent
\textbf{Application of FFT for fast computation} \quad  Making use of the Fast Fourier Transform (FFT) method, or more precisely, 2-dimensional FFT \cite{golub2013matrix} to accelerate our computation is crucial in our method. There are two places where we should apply FFT. The first is where we compute $G^{-1}u$; the second is where we compute $\PLap u$. We refer the readers to chapter 4.8 of \cite{golub2013matrix} and the references therein for details on implementing FFT. We also refer the reader to the examples in section \ref{numerical example} for applying FFT to general RD equations.

\vspace{0.3cm}

\noindent
\textbf{Choose suitable stepsize} \quad For the affine case discussed in Corollary \ref{corollary linear convergence}, if $G$ is close to $AA^\top$, then $\lambda_{\max}(A^\top G^{-1}A)$ should be a close to $1$, then $\tau_u\tau_p\leq \frac{4}{3\lambda^2_{\max}}$ indicates that we could choose stepsizes $\tau_u,\tau_p$ rather large. In our practice, starting at $\tau_u=\tau_p=0.8$ should be a reasonable choice. One can increase or decrease the stepsize based on the actual performance of the method.

\vspace{0.3cm}

\noindent
\textbf{Stopping criterion} \quad During our computation, we will set up a threshold $\delta$ for our method. After each iteration of the PDHG dynamic, we evaluate the $l^2$ norm of the residual $\textrm{Res}(U_n) = \frac{U_n - U^k}{h_t} - (\lambda \PLap U_n + f(U_n))$, we terminate the PDHG iteration iff $\|\textrm{Res}(U_n)\|_2\leq \delta$.

\begin{remark}[{Neumann boundary condition and DCT}]\label{remark: NLap and DCT}
It is worth providing some further discussions on our treatment of the Neumann boundary condition (BC), i.e., for a particular rectangular region $\Omega\subset\mathbb{R}^2$, $\frac{\partial u}{\partial \vec{n}} = 0$ on $\partial \Omega$. We discretize both sides of $\Omega$ into $N_x-1$ subintervals. Thus the space stepsize $h_x = \frac{L}{N_x-1}$. Such discretization will lead to $N_x^2$ mesh grid points. For point $(ih_x,0)$ on the vertical boundary of $\Omega$, we apply the central difference scheme to the Neumann boundary condition at the midpoint $(ih_x, -\frac{1}{2}h_x)$, which leads to $\frac{U_{i, -1} - U_{i, 0}}{h_x}=0$, thus $U_{i, -1}=U_{i, 0}$. Similar treatments are applied to the other boundaries of $\Omega$. Suppose we consider using the central difference scheme to discretize the Laplacian $\Delta$, let us denote the discretized Laplacian w.r.t. Neumann boundary condition as $\NLap$, then $\NLap$ takes the following form.
\begin{equation}
  \NLap = \frac{1}{h_x^2}\left[\begin{array}{ccccc}
      L_1 &  I  &   &   &    \\
      I &  L_2  & I &   &    \\
        & \ddots & \ddots & \ddots & \\
        &        &   I    &   L_2   & I \\
        &        &        &   I   & L_1
  \end{array}\right]_{N_x \times N_x \textrm{blocks}}.\label{def Lap N}
\end{equation}
Here the block matrices $L_1, L_2$ are
\begin{equation*}
 L_1 = \left[\begin{array}{ccccc}
      -2 &  1  &   &   &    \\
      1 &  -3  & 1 &   &    \\
        & \ddots & \ddots & \ddots & \\
        &        &   1    &   -3   & 1 \\
        &        &        &   1    & -2
  \end{array}\right]_{N_x\times N_x}, \quad L_2 = \left[\begin{array}{ccccc}
      -3 &  1  &   &   &    \\
      1 &  -4  & 1 &   &    \\
        & \ddots & \ddots & \ddots & \\
        &        &   1    &   -4   & 1 \\
        &        &        &   1    & -3
  \end{array}\right]_{N_x\times N_x}
\end{equation*}
Similar to using FFT for the computation involving $\PLap$, we can use the Discrete Cosine Transform (DCT) \cite{DCTstrang} \cite{golub2013matrix} to efficiently evaluate matrix-vector multiplication or solve linear equations involving the matrix $\NLap$. To be more specific, we use the DCT-2 transform introduced in \cite{DCTstrang} in the 2-dimensional scenario which enjoys the $O(N_x^2\log N_x)$ computational complexity.

\end{remark}

We summarize our method in the following algorithm. It is not hard to tell that the total complexity of each inner PDHG method is $O(\sharp\{\textrm{PDHG iter}\} \cdot N_x^2\log N_x\})$. 
\begin{algorithm}[!htb]
\caption{PDHG method for updating implicit one-step FD scheme of RD equation}\label{Algorithm_main}
\begin{algorithmic}[1]
\State \textbf{Input}: Initial condition $u_0$, terminal time $T$, number of time subintervals $N_t$; region size $L$, number of space subintervals $N_x$;  
\State Initialize $h_t = T / N_t$, $h_x = L / N_x$, $\{U^0_{ij}\} = \{u_0(ih_x, jh_x)\}$.
\For{$0\leq k\leq N_t-1$}
\State Set initial condition: $U_0 = U^k$ (or $\widehat{U}^{k+1}$ obtained via explicit Euler or IMEX scheme).
\State $n = 0$.
\While{$\|\textrm{Res}(U_n)\|_2\geq \delta$}
\State \textit{Evolve the G-prox PDHG dynamic \eqref{G-prox RD pdhg dynamic 1} - \eqref{G-prox RD pdhg dynamic 3} with help of 2D FFT(DCT):}
\State \textit{We use FFT for periodic BC and DCT for Neumann BC.}
\State Compute $V_n = \PLap U_n$ via 2D FFT(DCT);
\State Compute $W_n = U_n - \lambda h_t V_n - h_t f(U_n) - U^k$;
\State Solve $GY_n = W_n$ via 2D FFT(DCT);
\State Update $P_{n+1} = P_n + \tau_p Y_n$; $\widetilde{P}_{n+1} = P_{n+1} + \omega(P_{n+1}-P_n)$;
\State Compute $Q_{n+1} = \PLap \widetilde{P}_{n+1}$ via 2D FFT(DCT); 
\State Update $U_{n+1} = U_n - \tau_u(\widetilde{P}_{n+1} - \lambda h_t {Q}_{n+1} - h_t f'(U_n)\odot \widetilde{P}_{n+1})$;
\State $n = n+1$;
\EndWhile
\State Set $U^{k+1} = U_n$;
\EndFor
\State \textbf{Output}: The numerical solution $U^0, U^1,...,U^{N_t}$.
\end{algorithmic}
\label{alg: illust}
\end{algorithm}

In this section, we mainly focus on the one-step implicit finite difference (FD) scheme to illustrate how we apply the PDHG iterations to update the given FD scheme. But we should emphasize that our method is not restricted to such a scheme. One can extend the PDHG method to various types of numerical schemes by formulating the scheme at a certain time step $k$ as a nonlinear equation $F^k(U) = 0$, and construct the functional $L^k(U, P) = P^\top F^k(U)$. Then one can apply the dynamic \eqref{G-prox linearized pdhg dynamic 1} - \eqref{G-prox linearized pdhg dynamic 3} to update the numerical solution from $U^k$ to $U^{k+1}$. In addition, our method is applicable to more general reaction-diffusion equations (systems). Further discussions and details are supplied in section \ref{numerical example}.

\subsection{Discussion on convergence criteria and adaptive $h_t$ method}

Theorem \ref{thm linear convergence} suggests that under the linear case, the convergence of PDHG method relies on condition number $\kappa$, which is directly related to $h_x, h_t$ of our discrete scheme. In practice, we fix $h_x$ and $\tau_u,\tau_p$ in the algorithm. At every time step $n$ we discover that when time stepsize $h_t$ gets larger than a certain threshold value $h_t^*$ which depends on $h_x,\tau_u,\tau_p$ and $n$, PDHG method will hardly converge. The method works well when $h_t$ is slightly smaller than the threshold. The theoretical study on how $h_t^*$ guarantees the convergence of our method will be an important future research direction.

\textbf{Adaptive time stepsize} As discussed above, we cannot guarantee the convergence of the PDHG iteration \eqref{G-prox RD pdhg dynamic 1} - \eqref{G-prox RD pdhg dynamic 3} for any time stepsize $h_t$. Since the aforementioned threshold $h_t^*$ may vary at different time stages, and how $h_t^*$ varies depends on the nature of the equation as well as the discretization scheme. 
Given the potential difficulty of choosing the suitable $h_t$ that guarantees both the numerical accuracy as well as the convergence of the PDHG iterations, we come up with the strategy of using adaptive stepsize $h_t$ throughout the computation. We choose a time stepsize $h_t^0$ that guarantees the numerical accuracy and will serve as the upper bound of all $h_t$ throughout our method, we also set up two integers $N^* > N_* >0 $ as the thresholding integers for enlarging or shrinking the time stepsize $h_t$. We also pick a rescaling coefficient $ \eta \in(0,1)$. During each time step $n$ of the algorithm, we record the total number of PDHG iterations $M_{\mathrm{PDHG}}$, and reset the stepsize $h_t$ for the next time step based on the following rules:
\begin{itemize}
    \item  If $M_{\mathrm{PDHG}}>N^*$, we shrink $h_t$ by rate $\eta$, $h_t = \eta  h_t  $;
    \item If $N^*\geq M_{\mathrm{PDHG}}\geq N_*$, we remain $h_t$ unchanged;
    \item If $M_{\mathrm{PDHG}}<N_*$, we enlarge $h_t$ by rate $\frac{1}{\eta}$, $h_t = \frac{h_t}{\eta }$, if $\frac{h_t}{\eta }\leq h_t^0$; we remain $h_t$ unchanged if $\frac{h_t}{\eta}>h_t^0$.
\end{itemize}

It is also reasonable to fix $h_x, h_t$ but to only shrink the PDHG stepsizes $\tau_u, \tau_p$ when we encounter difficulties in converging. However, according to our experience, shrinking $\tau_u,\tau_p$ usually requires much more PDHG iterations for convergence, which may cause the algorithm less efficient compared with the aforementioned adaptive $h_t$ strategy.

\section{Numerical examples}\label{numerical example}
In this section, we demonstrate some numerical results computed by the proposed method. Throughout our experiments, we always use the extrapolation coefficient $\omega = 1$, and set the initial condition $U_0$ as the numerical solution $U^k$ computed from the last time step $t_k$. One can also try other values of  $ \omega$ or a more sophisticated initial guess of $U_0$. Our experiences show that confining $\omega$ around $1$ will probably provide the best performance of the PDHG method. Choosing $U_0$ as $\widehat{U}^{k+1}$ obtained by a specific explicit or IMEX scheme may slightly shorten the convergence time of the PDHG iteration. However, it is worth mentioning that when we are dealing with stiff equations, such treatment may introduce instability to the PDHG dynamic which may lead to the blow-up of the method.

Furthermore, as we have emphasized before, the mission of this paper is to verify the correctness and effectiveness of the proposed PDHG method in resolving the equation $F(U) = 0$ at each time step. Thus, in this work, we will mainly focus on the straightforward one-step implicit scheme \eqref{illustrate  implicit schm} in all numerical examples by omitting further discussions and experiments on more sophisticated numerical schemes.

We use fixed time stepsize $h_t$ in our numerical examples unless we emphasize that the adaptive $h_t$ method is applied in the experiments.

Our numerical examples are computed in MATLAB on a laptop with 11th Gen Intel(R) Core(TM) i5-1135G7 @ 2.40GHz 2.42 GHz CPU.

\subsection{Allen-Cahn equations}
The Allen-Cahn equation \cite{allen1979microscopic} is a typical reaction-diffusion equation taking the following form
\begin{equation}
  \frac{\partial u(x,t)}{\partial t} = a\Delta u(x,t) - bW'(u(x,t)), ~\textrm{on } \Omega\subset\mathbb{R}^2, \quad u(\cdot, 0) = u_0.  \label{AC equ}
\end{equation}
Here $a, b>0$ are positive coefficients, 
\begin{equation}
  W(u) = \frac{(u^2-1)^2}{4}  \label{AC CH potential}
\end{equation}
is a double-well potential function with its derivative $W'(u) = u^3-u$. We will always assume periodic boundary conditions in our discussion.

The Allen-Cahn equation can be treated as the $L^2-$gradient flow of the following functional $\mathcal{E}(u)$.
\begin{equation}
    \mathcal{E}(u) = \int_\Omega \frac{1}{2}a|\nabla u|^2 + bW(u)~dx.  \label{AC, CH energy}
\end{equation}
\subsubsection{Examples with shrinking level set curve}
In the first example, we consider $\Omega = [-L, L]^2$ with $L=0.25$. We consider taking $a=\epsilon$, $b = \frac{1}{\epsilon}$ with $\epsilon = 0.01$. Let use consider $u_0=2\chi_{B}-1$. Here $\chi_E$ denotes the indicator function of measurable set $E$, i.e., $\chi_E(x)=1$ if $x\in E$, and $\chi_E(x)=0$ otherwise. We denote $B$ as the disk centered at $O$ with a radius equal to $0.2$. Suppose we use periodic boundary conditions for \eqref{AC equ}. It is well-known that the zero-level-set curve of the solution $u(x,t)$ to Allen-Cahn equation behaves similarly to the mean curvature flow as time $t$ increases \cite{fastreactionslowdiffusion, merriman1992diffusion, merriman1994motion}. In this case, we can treat the initial level set curve as the circle centered at the origin with radius $r(0)=0.2$. As $t$ increases, the circle radius will shrink at the rate of $\epsilon$ times circle curvature, i.e., $\dot r(t) = -\epsilon\kappa(t) = -\frac{\epsilon}{r(t)}$. Solving this equation leads to $r(t) = \sqrt{r(0)^2 - 2\epsilon t}$. Thus the level set circle will vanish at finite time $t = \frac{r(0)^2}{2\epsilon} = 2$. 

\begin{figure}[htb!]
\begin{subfigure}{.161\textwidth}
  \centering
  \includegraphics[width=\linewidth]{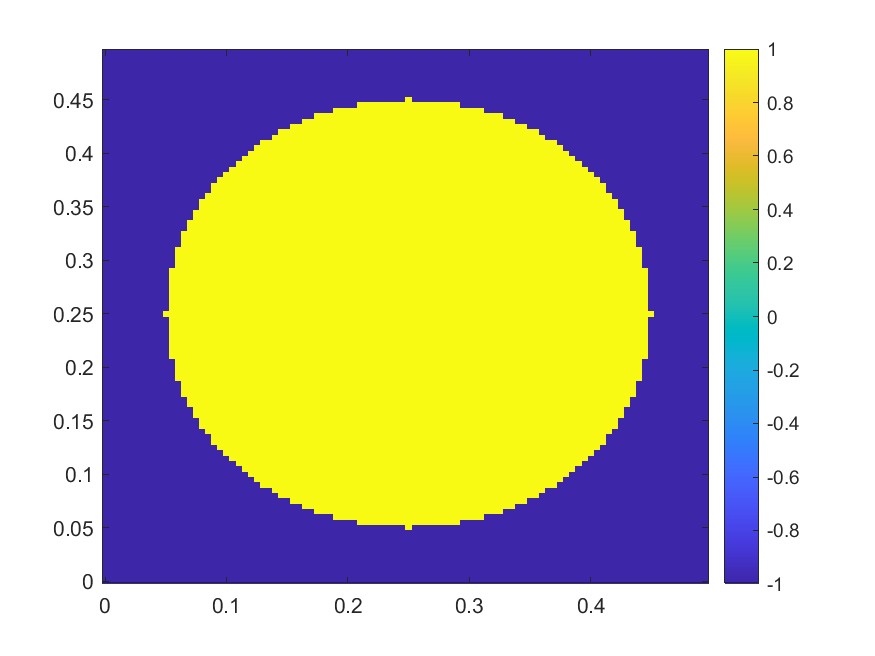}
  \caption{$t=0.0$}
\end{subfigure}
\begin{subfigure}{.161\textwidth}
  \centering
  \includegraphics[width=\linewidth]{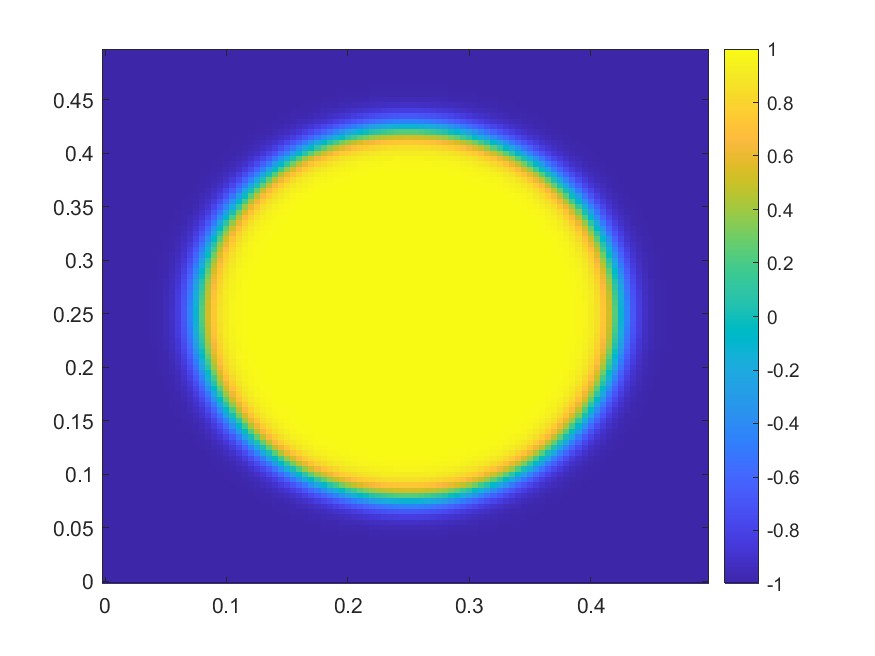}
  \caption{$t=0.5$}
\end{subfigure}
\begin{subfigure}{.161\textwidth}
  \centering
  \includegraphics[width=\linewidth]{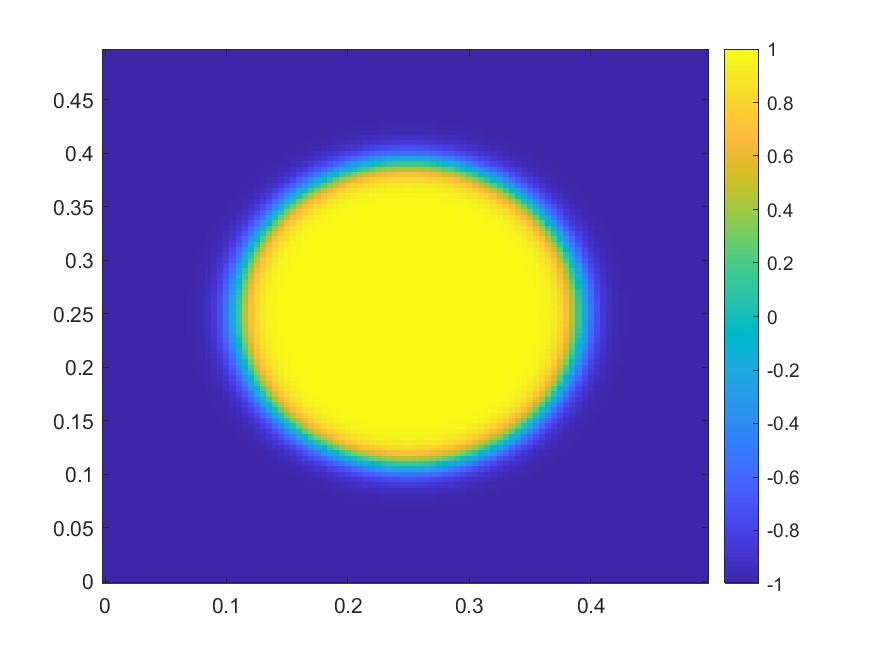}
  \caption{$t=1.0$}
\end{subfigure}
\begin{subfigure}{.161\textwidth}
  \centering
  \includegraphics[width=\linewidth]{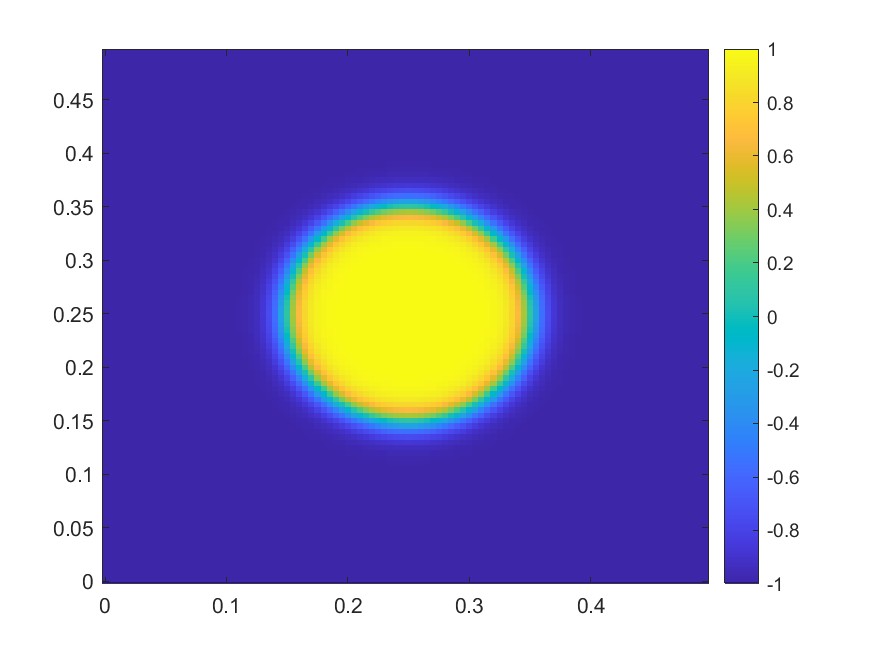}
  \caption{$t=1.5$}
\end{subfigure}
\begin{subfigure}{.161\textwidth}
  \centering
  \includegraphics[width=\linewidth]{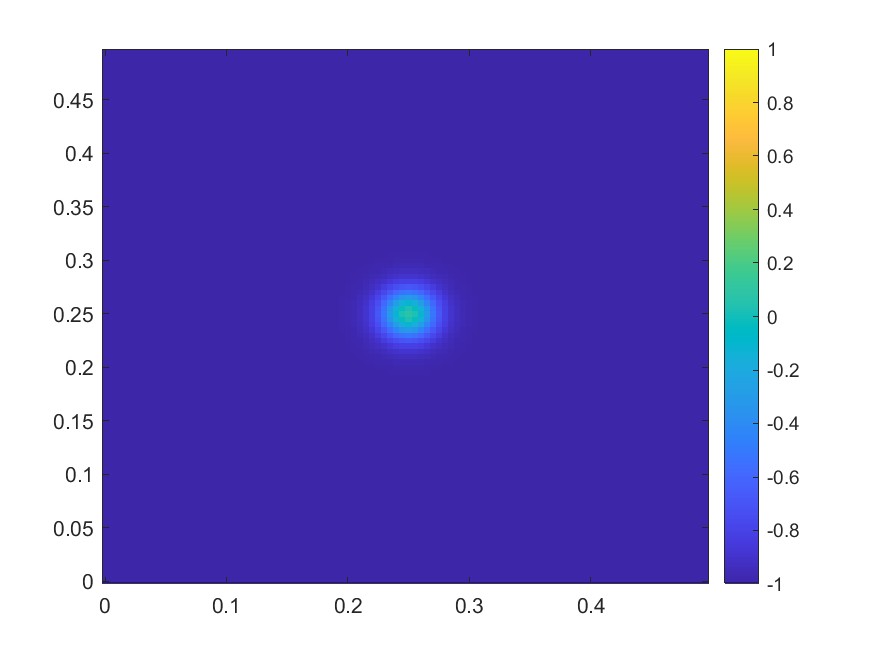}
  \caption{$t=2.0$}
\end{subfigure}
\begin{subfigure}{.161\textwidth}
  \centering
  \includegraphics[width=\linewidth]{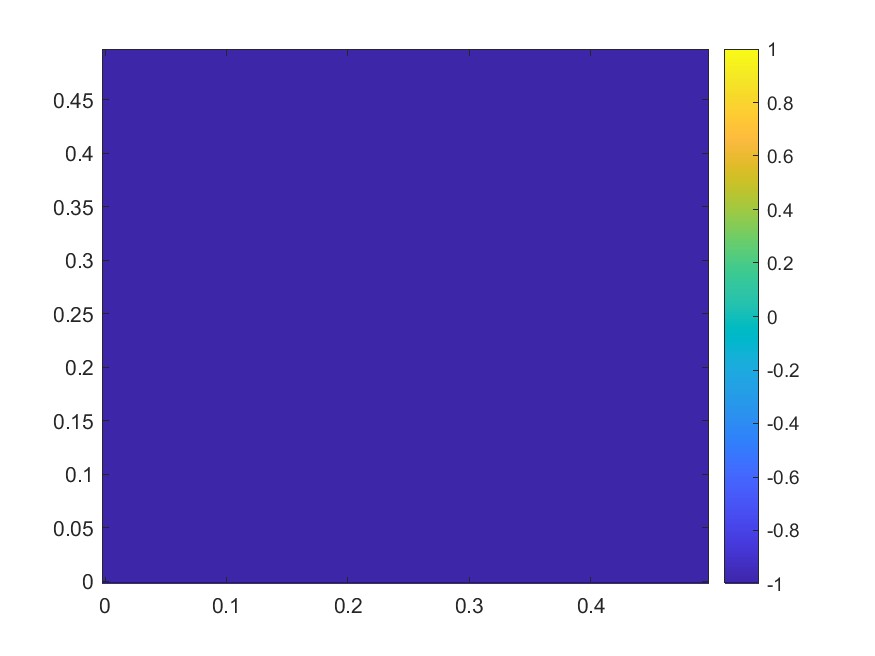}
  \caption{$t=2.5$}
\end{subfigure}
\caption{Numerical solution of \eqref{allen cahn 1} at different times with initial condition $u_0=2\chi_B-1.$}\label{allen cahn 1}
\end{figure}

\begin{figure}[htb!]
\begin{subfigure}[t]{.32\textwidth}
  \centering
  \includegraphics[width=\linewidth]{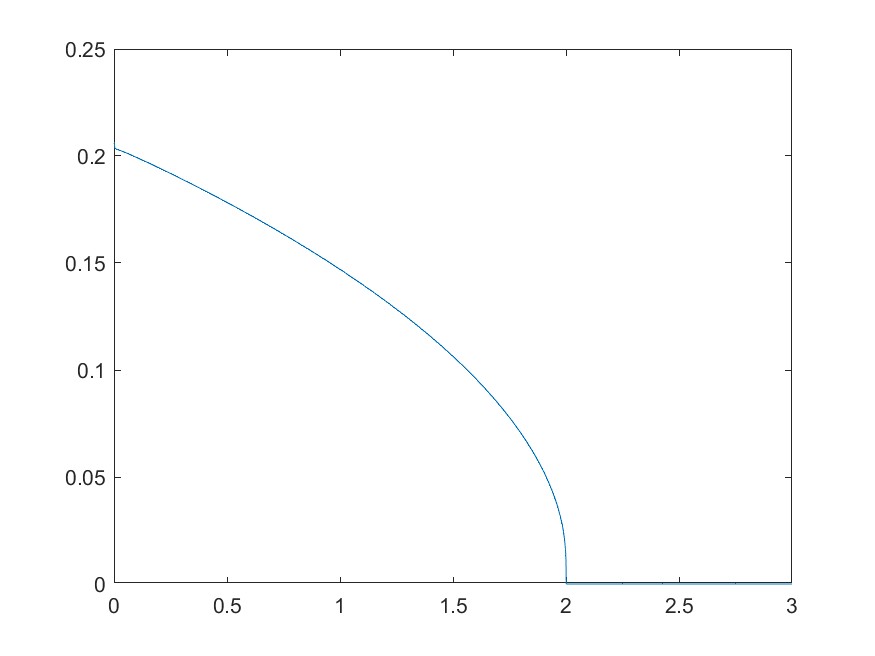}
  \caption{Plot of the front position (calculated from linear interpolation of the numerical solution) in radial direction versus time}
\end{subfigure}
\begin{subfigure}[t]{.32\textwidth}
  \centering
  \includegraphics[width=\linewidth]{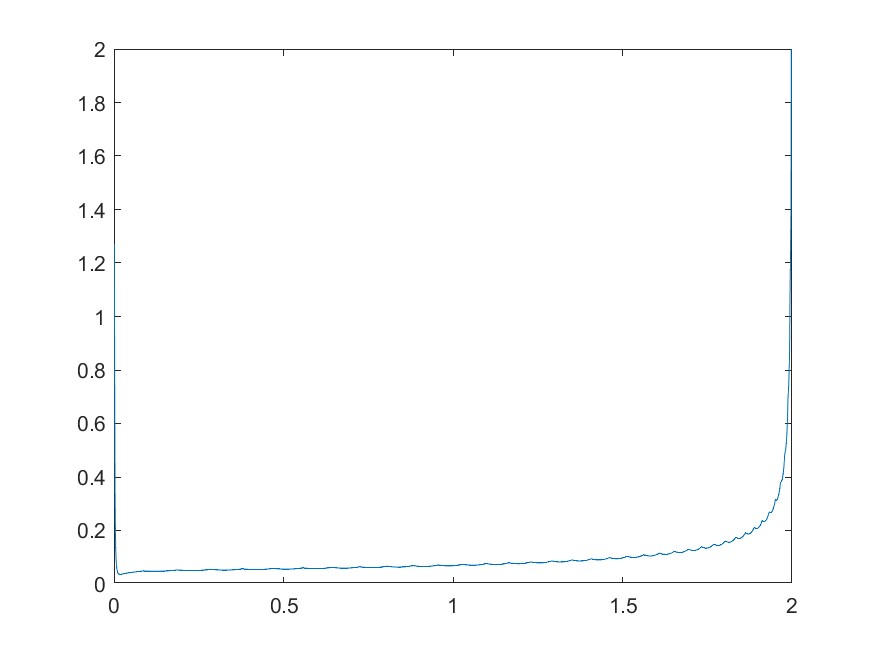}
  \caption{Plot of the front speed (calculated from finite difference) versus time}
\end{subfigure}
\begin{subfigure}[t]{.32\textwidth}
  \centering
  \includegraphics[width=\linewidth]{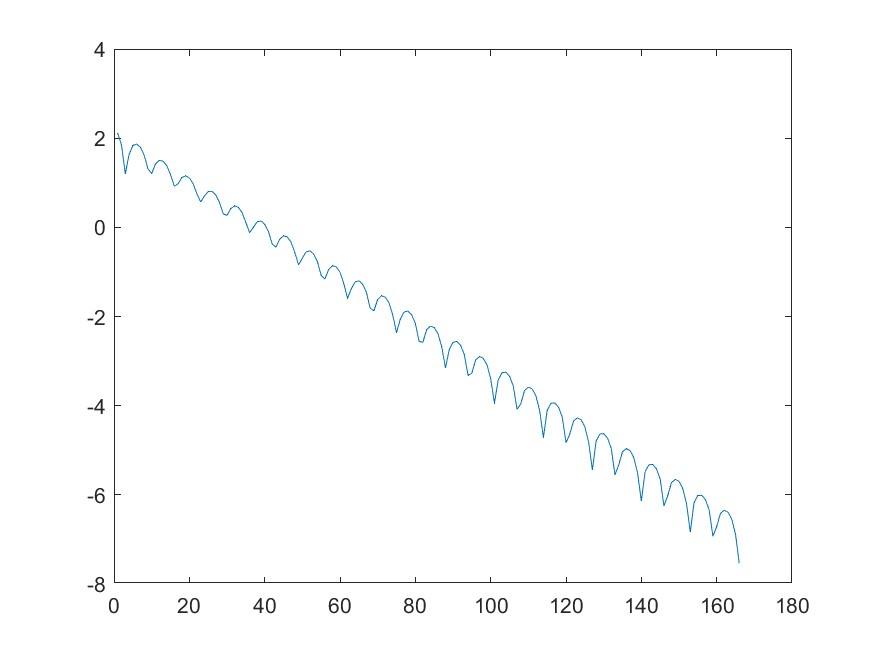}
  \caption{Plot of the front speed (calculated from finite difference) versus time}
\end{subfigure}
\caption{Plots of front position and speed (Left \& Middle); Plots of $\log_{10}\mathrm{Res}(U_n)$ versus PDHG iterations at time stage $t=1.0$ (Right).}\label{allen cahn 1 front and convergence}
\end{figure}

As suggested in section 4.4 of \cite{merriman1992diffusion}, it is important to choose the spatial stepsize $h_x$ small enough so that $h_x$ no larger than $O(\epsilon)$ to capture the shrinkage of level set curve, otherwise, the numerical solution may get stuck at some intermediate stage. In this example, we solve the equation on time interval $[0, 3]$. We choose $N_t = 3000$, thus $h_t = 1/1000$; $N_x=100$ with $h_x = L/N_x=1/200$. Recall $h_x<\epsilon$. We choose $\tau_u = \tau_p = 0.5$ as the stepsize for the PDHG iteration. Some computed results are shown in Figure \ref{allen cahn 1}. Plots of the radial position as well as the moving speed of the front (zero level set circle) of the numerical solution are presented in Figure \ref{allen cahn 1 front and convergence}. 

We apply the PDHG method described in Algorithm \ref{Algorithm_main} to this problem. Although equation \eqref{allen cahn 1} contains a nonlinear term with a significant coefficient $\frac{1}{\epsilon}$, our proposed PDHG method still solves the nonlinear $F(U)=0$ efficiently. To be more specific, we set the PDHG-threshold $\delta=10^{-7}$. Most PDHG iterations will terminate in less than $200$ steps for each discrete time step. The right plot of Figure \ref{allen cahn 1 front and convergence} indicates the linear convergence of the proposed method.

\begin{figure}[htb!]
\begin{subfigure}{.161\textwidth}
  \centering
  \includegraphics[width=\linewidth]{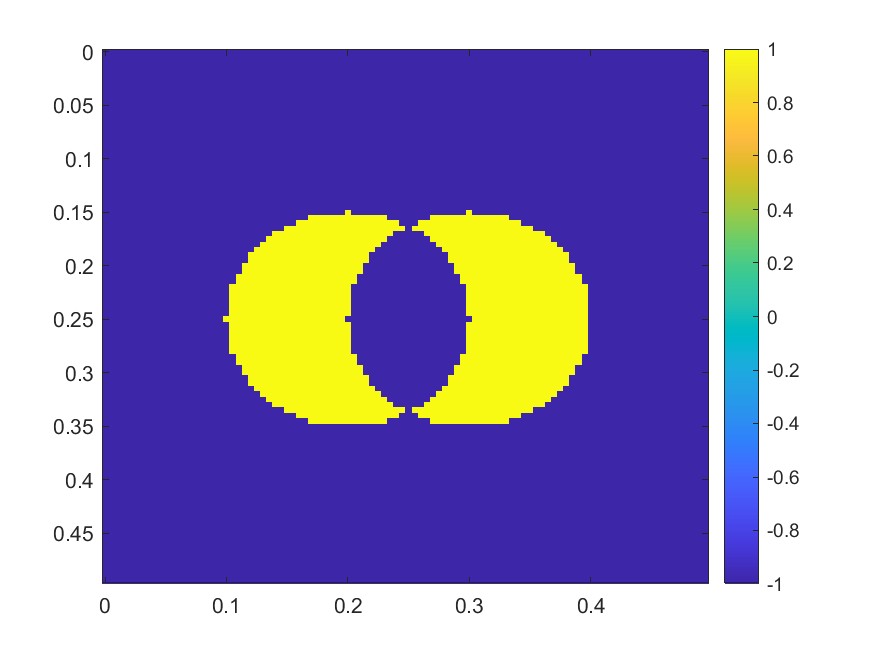}
  \caption{$t=0.0$}
\end{subfigure}
\begin{subfigure}{.161\textwidth}
  \centering
  \includegraphics[width=\linewidth]{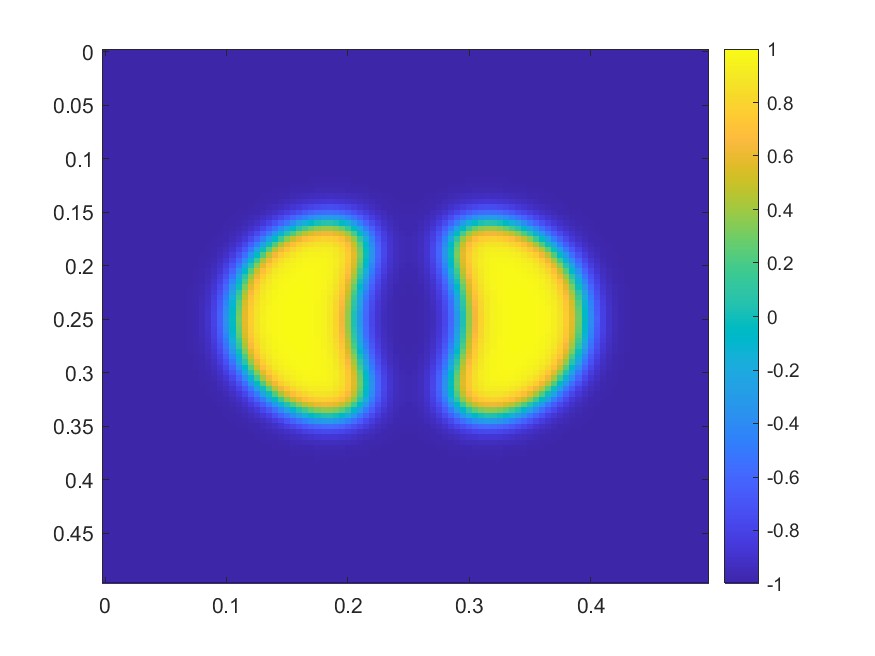}
  \caption{$t=0.05$}
\end{subfigure}
\begin{subfigure}{.161\textwidth}
  \centering
  \includegraphics[width=\linewidth]{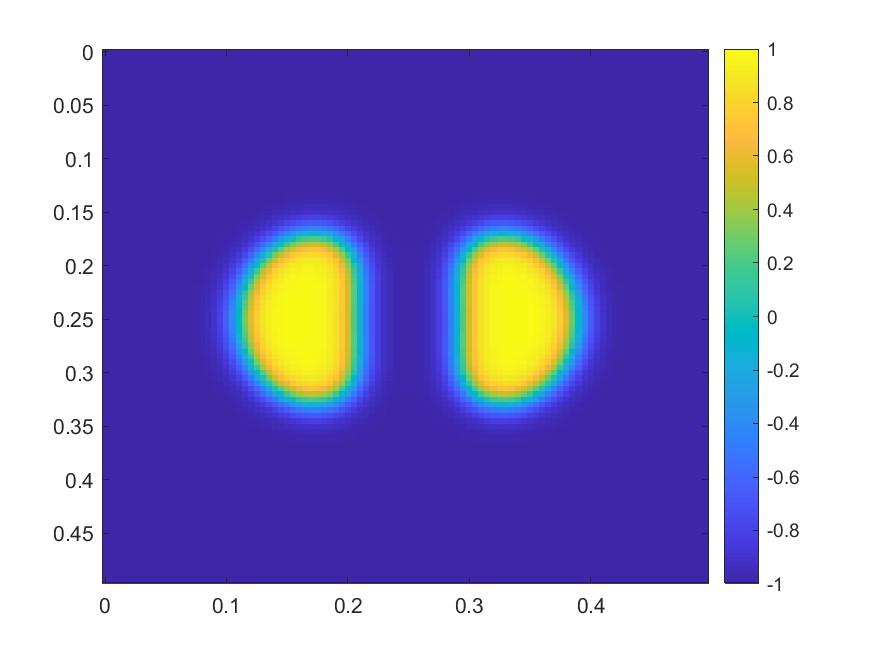}
  \caption{$t=0.1$}
\end{subfigure}
\begin{subfigure}{.161\textwidth}
  \centering
  \includegraphics[width=\linewidth]{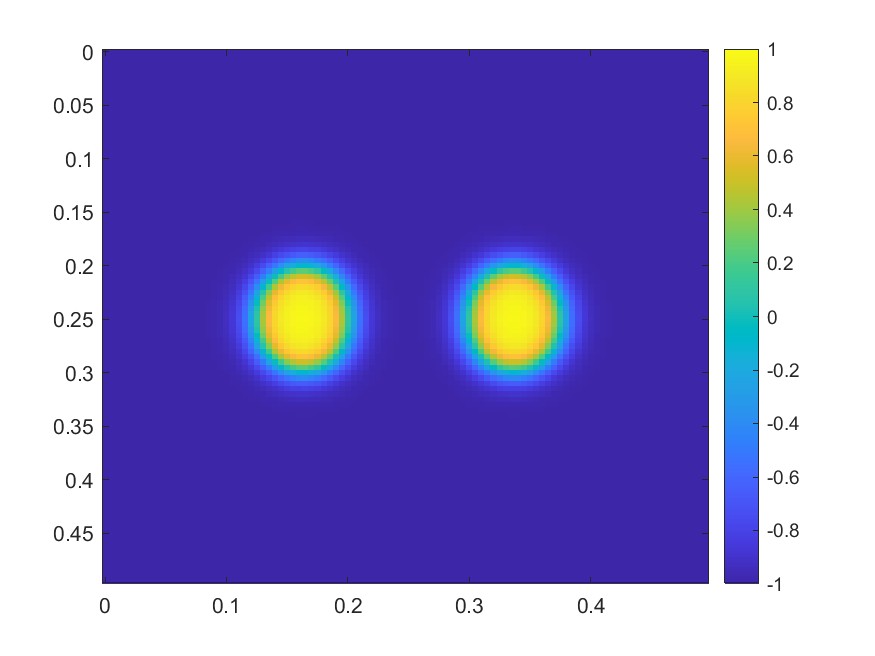}
  \caption{$t=0.2$}
\end{subfigure}
\begin{subfigure}{.161\textwidth}
  \centering
  \includegraphics[width=\linewidth]{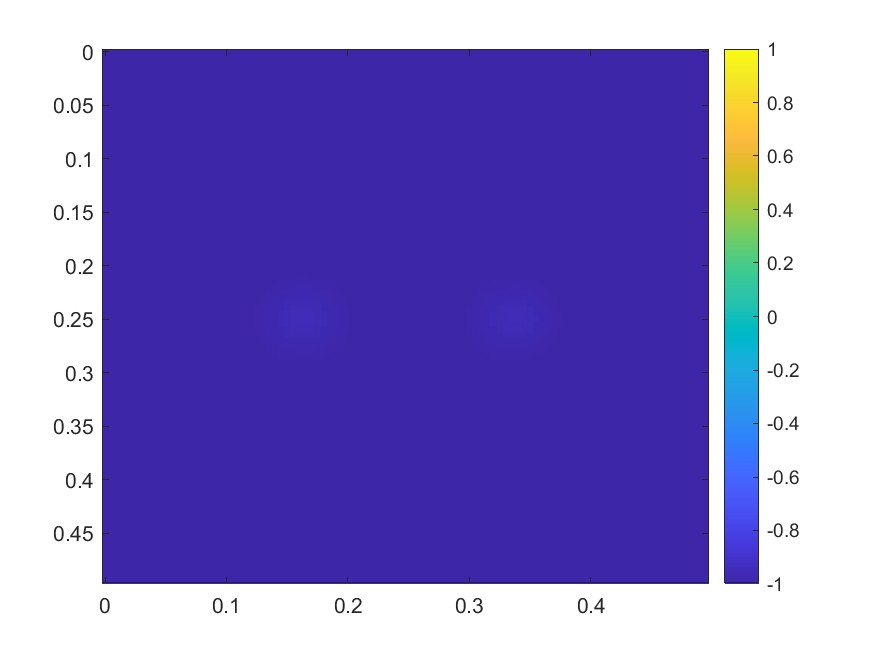}
  \caption{$t=0.3$}
\end{subfigure}
\begin{subfigure}{.161\textwidth}
  \centering
  \includegraphics[width=\linewidth]{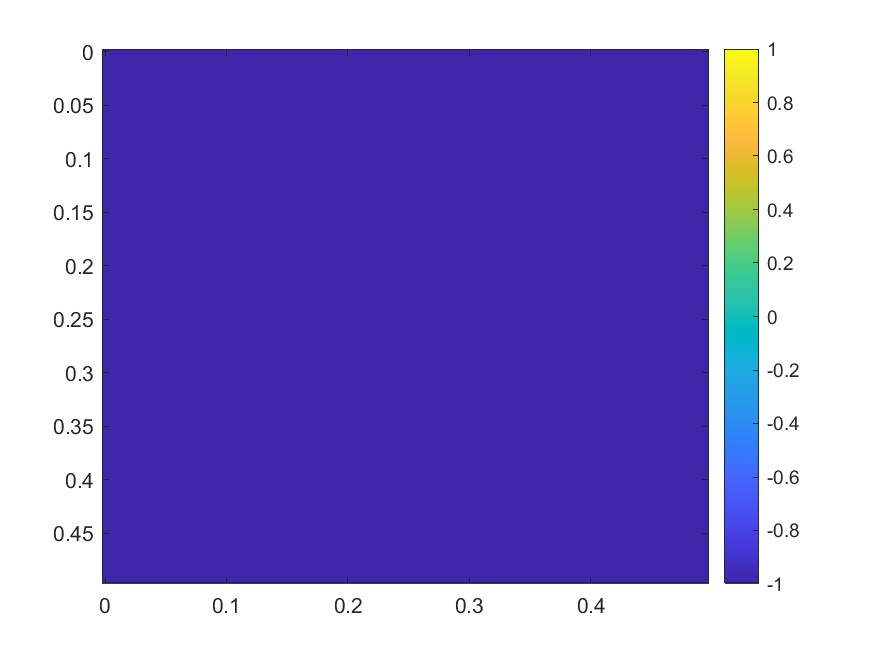}
  \caption{$t=0.4$}
\end{subfigure}
\caption{Numerical solution at different times with initial condition $u_0=2\chi_E-1$. Notice that in the last plot, we have almost converged to the equilibrium solution $u=-1$.}\label{allen cahn 2}
\end{figure}

We also solve equation \eqref{allen cahn 1} on $\Omega = [0, 0.5]^2$ within the time interval $[0, 0.5]$. We still set $a = \epsilon, b = \frac{1}{\epsilon}$ with $\epsilon = 0.01$. We pick $N_x = 100$, $N_t = 500$, thus $h_x = 1/200$, $h_t = 1/1000$. We consider the initial condition $u_0=2\chi_E-1$ with the region $E=(B_1\setminus B_2)\cup(B_2 \setminus B_1)$, where $B_1, B_2$ are disks centered at $(0.2, 0.25)$, $(0.3, 0.25)$ with radius both equal to $0.1$. We apply the PDHG method with $\tau_u=\tau_p = 0.5$ and obtain the numerical results in Figure \ref{allen cahn 2}. In this example, the PDHG method takes no more than $200$ iterations for each time step update.

\subsection{Cahn-Hilliard equations}
We now switch to another well-known reaction-diffusion equation known as the Cahn-Hilliard equation \cite{cahn1961spinodal}, which takes the following form.
\begin{equation}
  \frac{\partial u(x,t)}{\partial t} = - a\Delta\Delta u(x,t) + b\Delta W'(u(x,t)), ~\textrm{on } \Omega\subset\mathbb{R}^2, \quad u(\cdot, 0) = u_0.  \label{CH equ}
\end{equation}
Here we assume $a, b > 0$, and $W(u)$ defined the same as in the Allen-Cahn equation. In this section, we will restrict our discussion to periodic boundary conditions. Similar to the Allen-Cahn equation, the Cahn-Hilliard equation can be treated as the $H^{-1}-$gradient flow of the functional $\mathcal{E}(u)$ defined in \eqref{AC, CH energy}. Due to this reason, compared with the Allen-Cahn equation, the Cahn-Hilliard equation involves one extra operator $-\Delta$ on the right-hand side of the equation. This difference leads to several slight modifications to our original algorithm.

The functional $F(U)$ introduced in \eqref{F(U) rd equation} is now
\begin{equation}
   F(U) =  (I-\lambda h_t\PLap\PLap)U - U^k - h_t \PLap f(U).  \label{F(U) CH}
\end{equation}
Thus $L(U, P) = P^\top((I-\lambda h_t\PLap\PLap)U - U^k) - h_t P^\top\PLap f(U).$ It is then natural to choose precondition matrix $G$ as $(I - \lambda h_t\PLap\PLap)^2$.
The three-step PDHG update for Cahn-Hilliard equation can thus be formulated as
\begin{align}
  P_{n+1} & =  P_n + \tau_p G^{-1}(U_n - \lambda h_t \PLap \PLap U_n -h_t\PLap f(U_n) - U^k); \label{G-prox RD pdhg dynamic CH 1}\\
  \widetilde{P}_{n+1} & = P_{n+1} + \omega (P_{n+1} - P_n); \label{G-prox RD pdhg dynamic CH 2}\\
  U_{n+1} & = U_n - \tau_u(\widetilde{P}_{n+1} - \lambda h_t \PLap\PLap\widetilde{P}_{n+1} - h_t f'(U_n)\odot \PLap\widetilde{P}_{n+1}).\label{G-prox RD pdhg dynamic CH 3}
\end{align}
By carefully investigating the steps among \eqref{G-prox RD pdhg dynamic CH 1} - \eqref{G-prox RD pdhg dynamic CH 3}, one can tell that both the linear equation involving $G$ and the matrix-vector multiplication involving $\PLap$ can be computed via FFT, which indicates the effectiveness of the computational scheme when applied to Cahn-Hilliard equations.

We demonstrate several numerical examples below.
\subsubsection{Example with seven circles}
Inspired by the second example introduced in \cite{church2019high}, we consider Cahn-Hilliard equation \eqref{CH equ} on periodic domain $\Omega = [0, 2\pi]^2$ with $a=0.1^2$ and $b=1$. We set the initial condition $u_0$ as
\[u_0(x, y) = -1 + \sum_{i=1}^7 \varphi(\sqrt{(x-x_i)^2 + (y-y_i)^2} - r_i),\]
where the mollifier function $\varphi$ is defined as
\[ \varphi(s) = \begin{cases}
    2e^{-\frac{\epsilon^2}{s^2}} \quad s<0;\\
    0 \quad s \geq 0
\end{cases}, \quad \textrm{with } \epsilon = 0.1. \]
One can think of $u_0(x,y)$ as an indicator function whose value equals $+1$ if $(x,y)$ falls into any of the seven circles; and equals $-1$ otherwise. Furthermore, we set the centers and radii of the seven circles as in Table \ref{tab: seven circles}.
\begin{table}[]
    \centering
    \begin{tabular}{|c||c|c|c|c|c|c|c|}
    \hline
          $i$ &  $1$     &  $2$        &  $3$        & $4$        & $5$      & $6$   & $7$  \\
          \hline
          \hline
       $x_i$  &  $\pi/2$ &  $\pi/4$    &  $\pi/2$   & $\pi$       & $3\pi/2$ & $\pi$ & $3\pi/2$\\
       $y_i$  &  $\pi/2$ &  $3\pi/4$   &  $5\pi/4$  & $\pi/4$     & $\pi/4$  & $\pi$ & $3\pi/2$\\
       $r_i$  &  $\pi/5$ &  $2\pi/15$  &  $\pi/15$  & $\pi/10$    & $\pi/10$ & $\pi/4$ & $\pi/4$\\
    \hline
    \end{tabular}
    \caption{data 7 circles}
    \label{tab: seven circles}
\end{table}

\begin{figure}[htb!]
\begin{subfigure}{.161\textwidth}
  \centering
  \includegraphics[width=\linewidth]{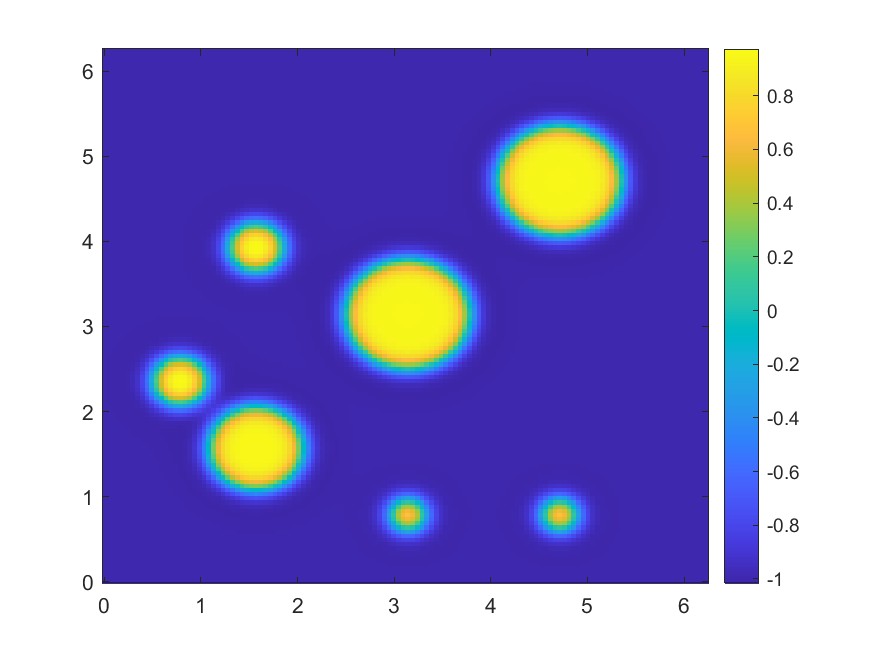}
\end{subfigure}
\begin{subfigure}{.161\textwidth}
  \centering
  \includegraphics[width=\linewidth]{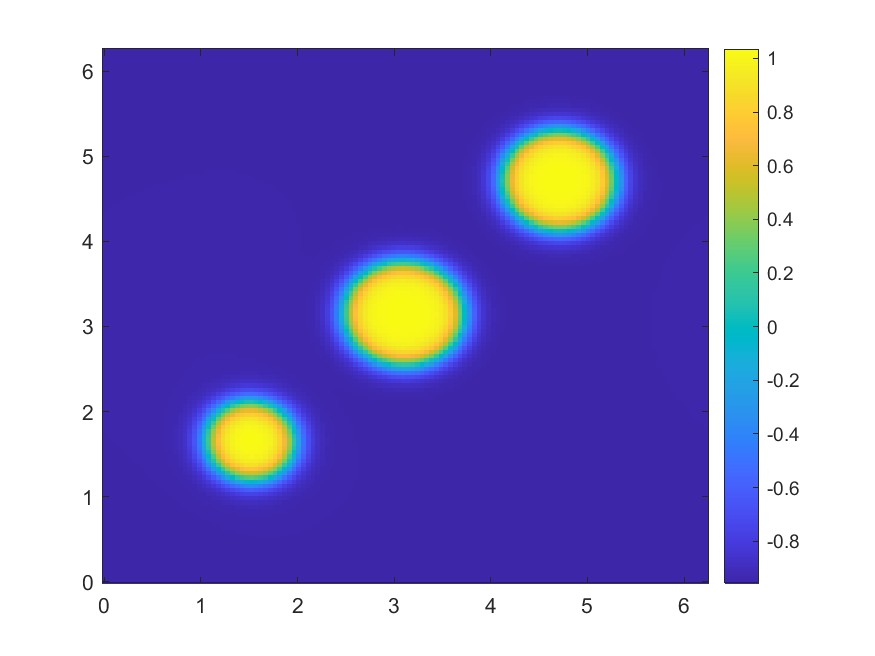}
\end{subfigure}
\begin{subfigure}{.161\textwidth}
  \centering
  \includegraphics[width=\linewidth]{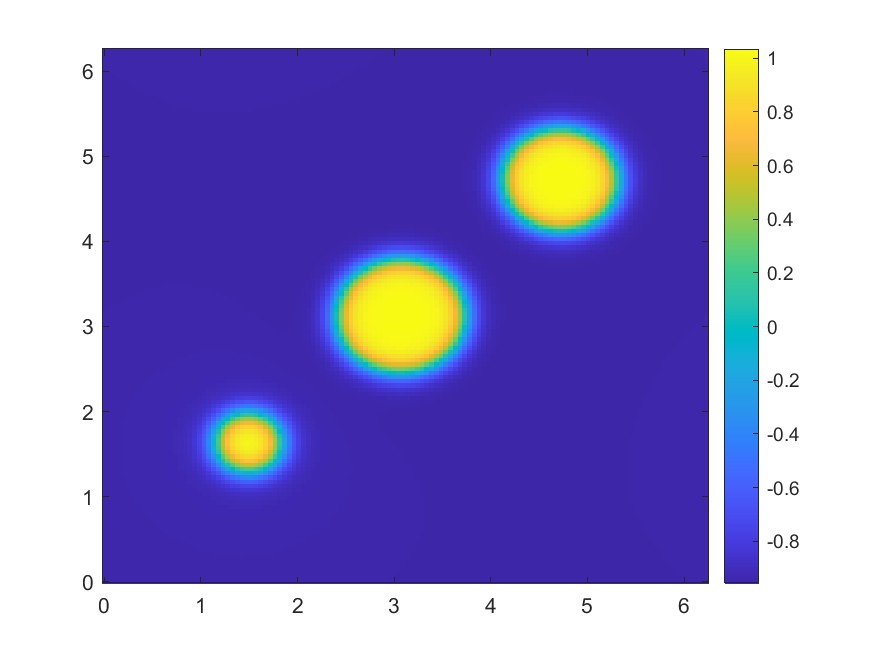}
\end{subfigure}
\begin{subfigure}{.161\textwidth}
  \centering
  \includegraphics[width=\linewidth]{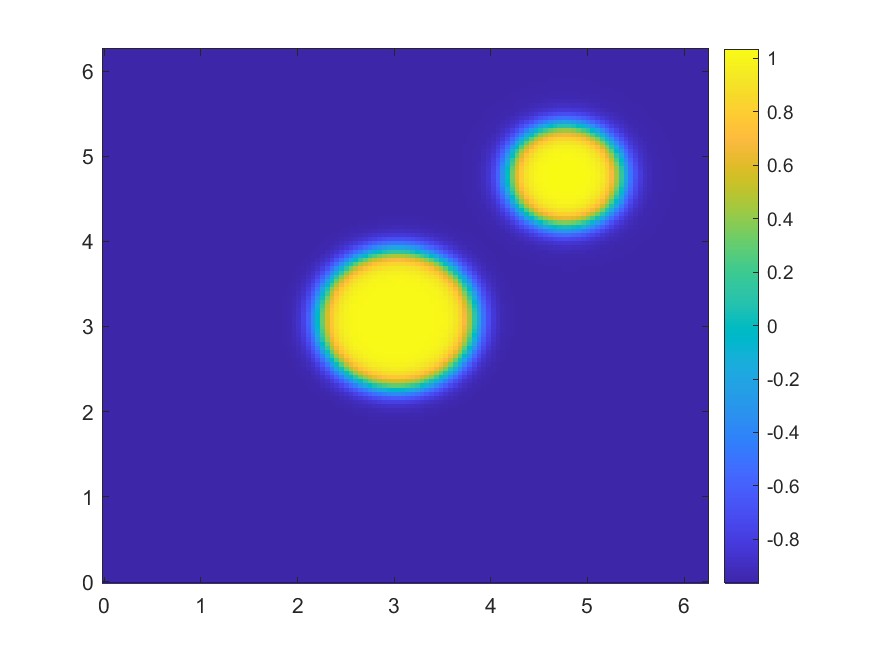}
\end{subfigure}
\begin{subfigure}{.161\textwidth}
  \centering
  \includegraphics[width=\linewidth]{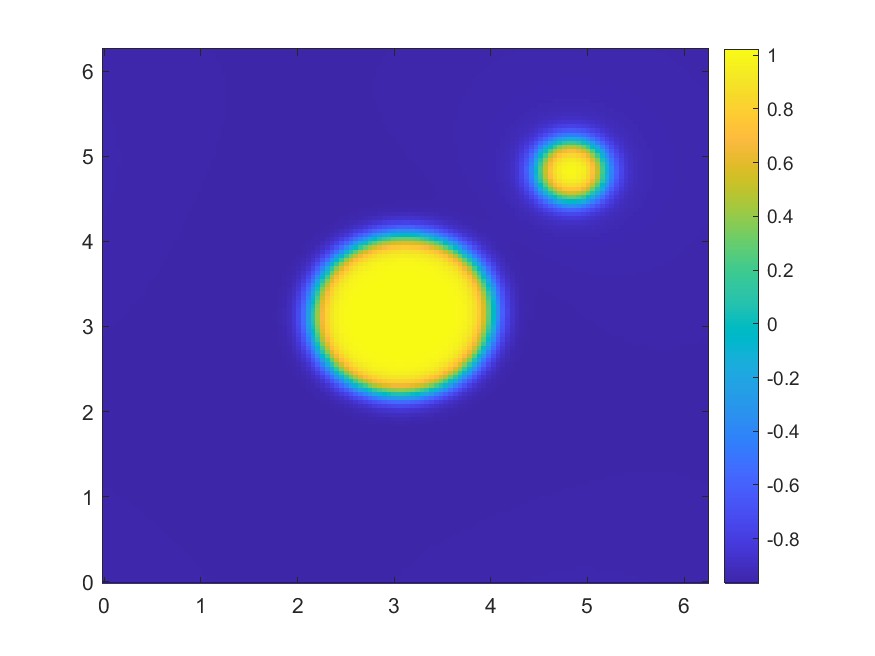}
\end{subfigure}
\begin{subfigure}{.161\textwidth}
  \centering
  \includegraphics[width=\linewidth]{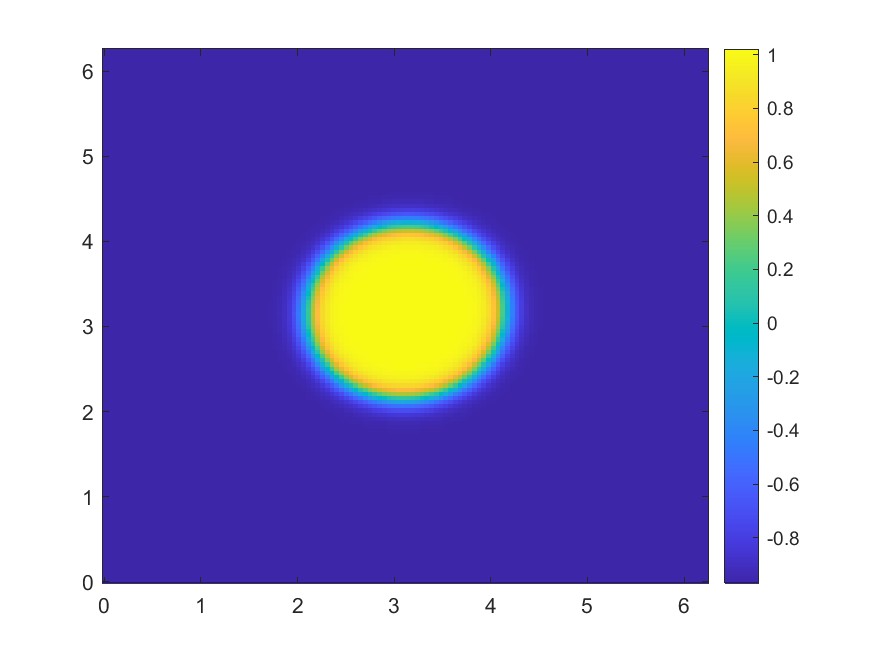}
\end{subfigure}

\begin{subfigure}{.161\textwidth}
  \centering
  \includegraphics[width=\linewidth]{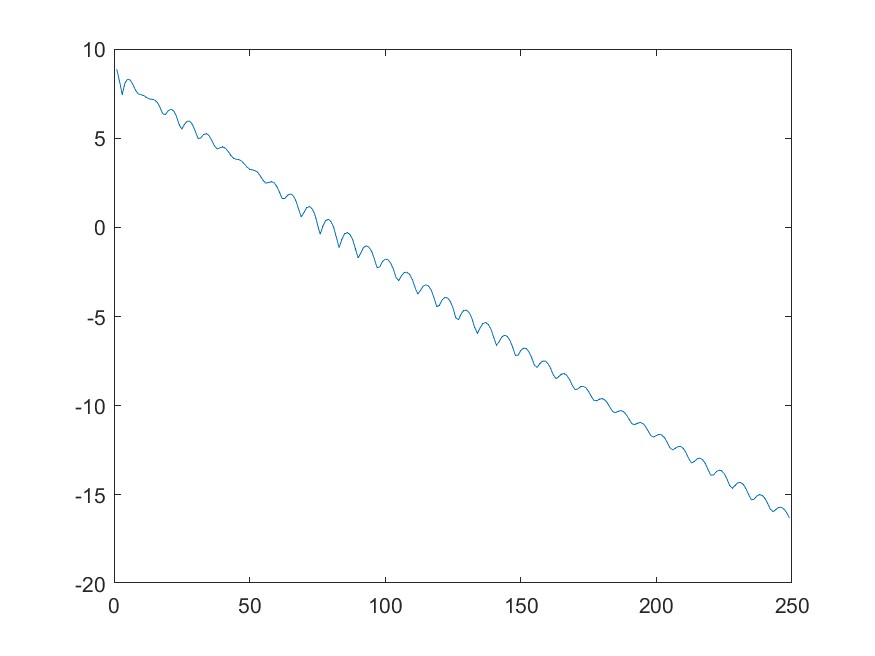}
  \caption{$t=0.0$}
\end{subfigure}
\begin{subfigure}{.161\textwidth}
  \centering
  \includegraphics[width=\linewidth]{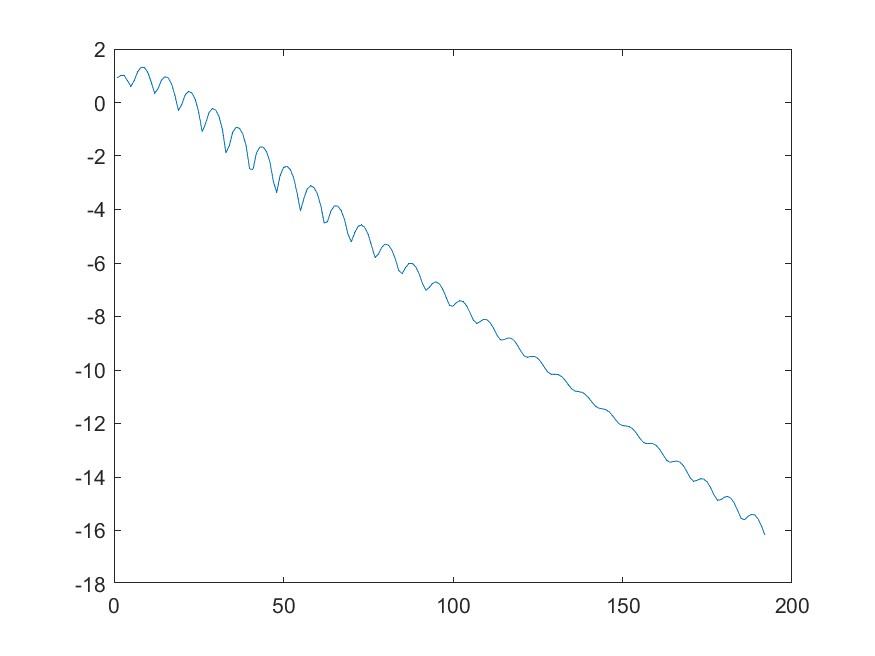}
  \caption{$t=1.0$}
\end{subfigure}
\begin{subfigure}{.161\textwidth}
  \centering
  \includegraphics[width=\linewidth]{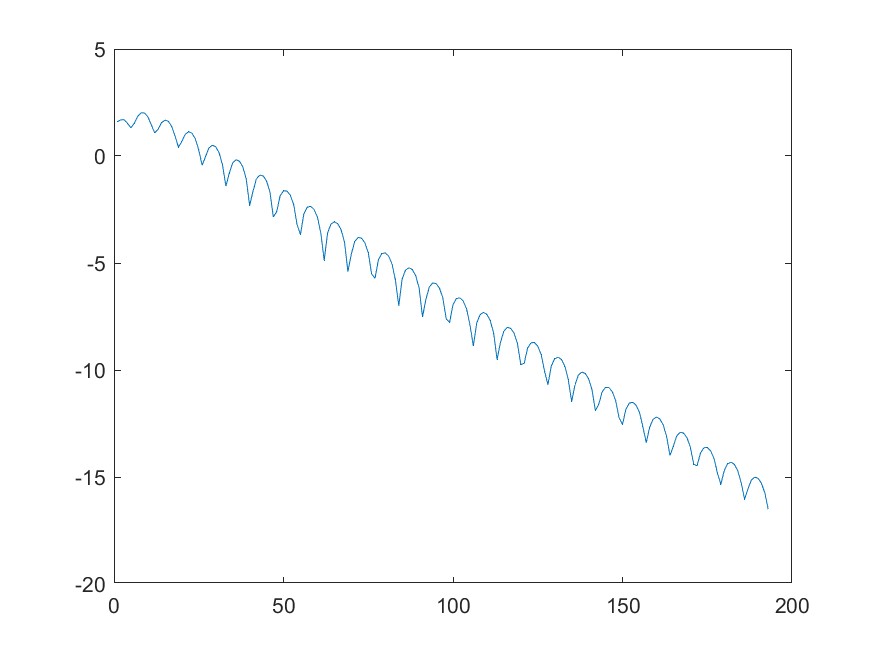}
  \caption{$t=5.0$}
\end{subfigure}
\begin{subfigure}{.161\textwidth}
  \centering
  \includegraphics[width=\linewidth]{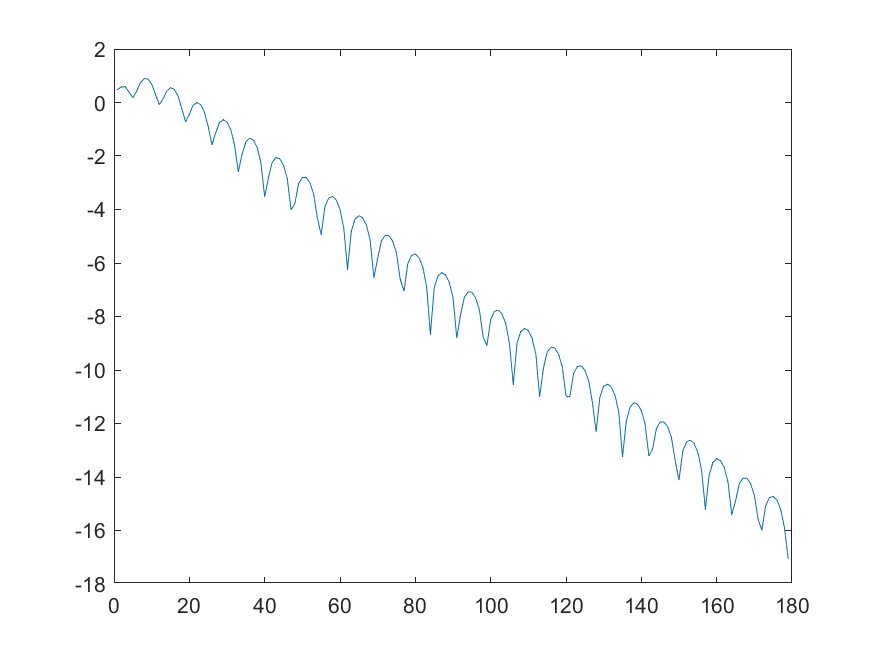}
  \caption{$t=15.0$}
\end{subfigure}
\begin{subfigure}{.161\textwidth}
  \centering
  \includegraphics[width=\linewidth]{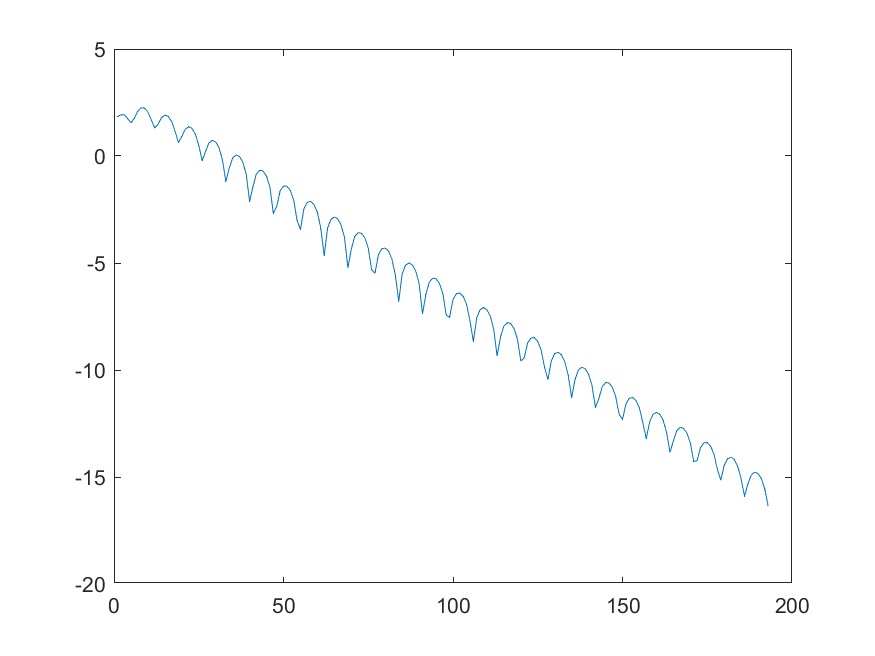}
  \caption{$t=25.0$}
\end{subfigure}
\begin{subfigure}{.161\textwidth}
  \centering
  \includegraphics[width=\linewidth]{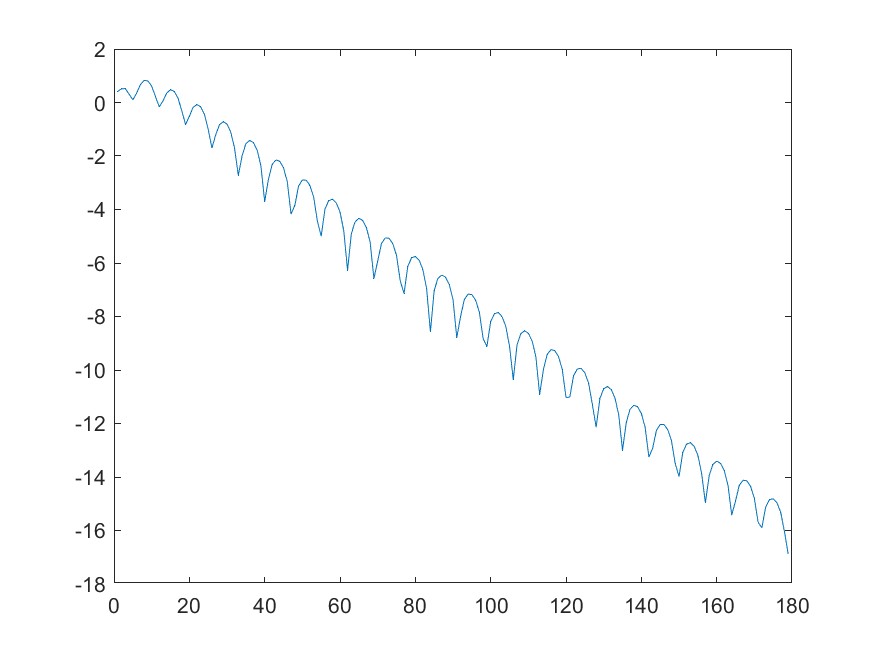}
  \caption{$t=30.0$}
\end{subfigure}

\caption{Numerical solution and $\log_{10}\mathrm{Res}(U_n)$ plot at different time stages for the seven circle example. The residual plots indicate the linear convergence of PDHG method for the nonlinear objective functions used in this example.}\label{cahn hilliard 1}
\end{figure}

We will solve equation \eqref{CH equ} on the time interval $[0, 30]$. In our numerical implementation, we set $N_x = 128$, $h_x = \pi/64$; $N_t = 6000$, $h_t = 1/200$. For the PDHG iteration, we set $\tau_u = \tau_p = 0.5$. The numerical solution to this equation is demonstrated in Figure \ref{cahn hilliard 1}. The plots of $\log$ residuals at different time stages are also presented in Figure \ref{cahn hilliard 1}, which exhibit the linear convergence of the PDHG algorithm. The small circles will gradually fade out, leaving the largest circle in the center of the domain till the end. By analyzing our numerical solution, the time $T_1$ at which the value of our numerical solution is evaluated at $(\pi/2, \pi/2)$ passes $0$ is located in the interval $[6.340, 6.345]$; while the time $T_2$ at which our numerical value evaluated at $(3\pi/2, 3\pi/2)$ passes $0$ is located in the interval $[26.015, 26.020]$. Both times meet the accuracy proposed in \cite{church2019high}.

\subsubsection{Example with sinusoidal initial condition}
\begin{figure}[htb!]
\begin{subfigure}{.23\textwidth}
  \centering
  \includegraphics[width=\linewidth]{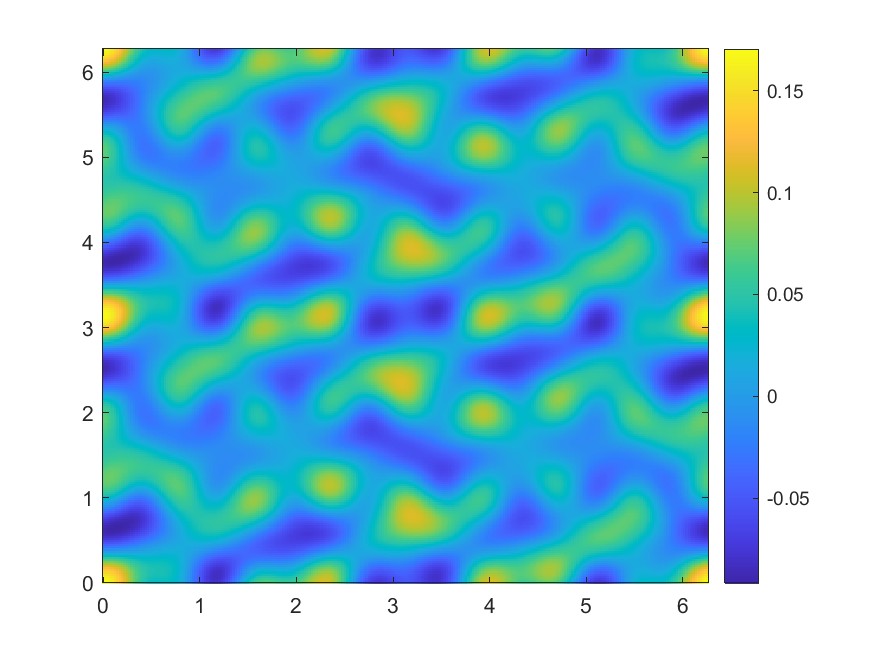}
  \caption{$t=0.0$}
\end{subfigure}
\begin{subfigure}{.23\textwidth}
  \centering
  \includegraphics[width=\linewidth]{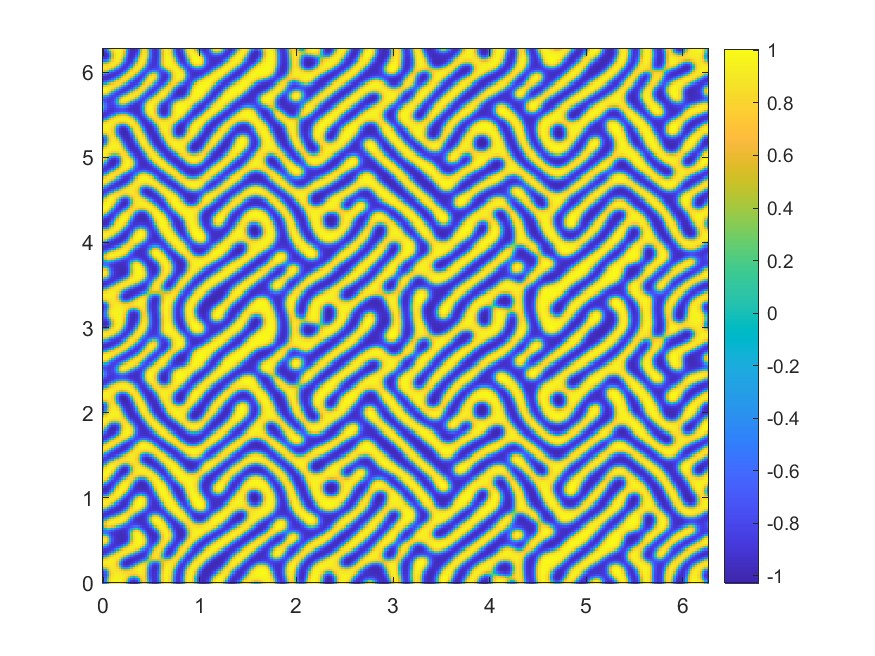}
  \caption{$t=0.05$}
\end{subfigure}
\begin{subfigure}{.23\textwidth}
  \centering
  \includegraphics[width=\linewidth]{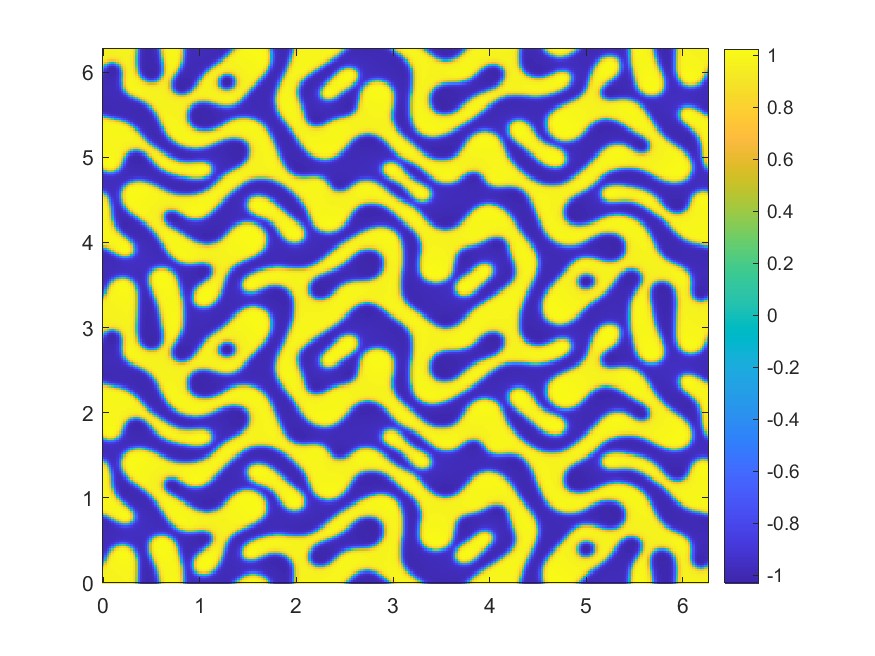}
  \caption{$t=0.30$}
\end{subfigure}
\begin{subfigure}{.23\textwidth}
  \centering
  \includegraphics[width=\linewidth]{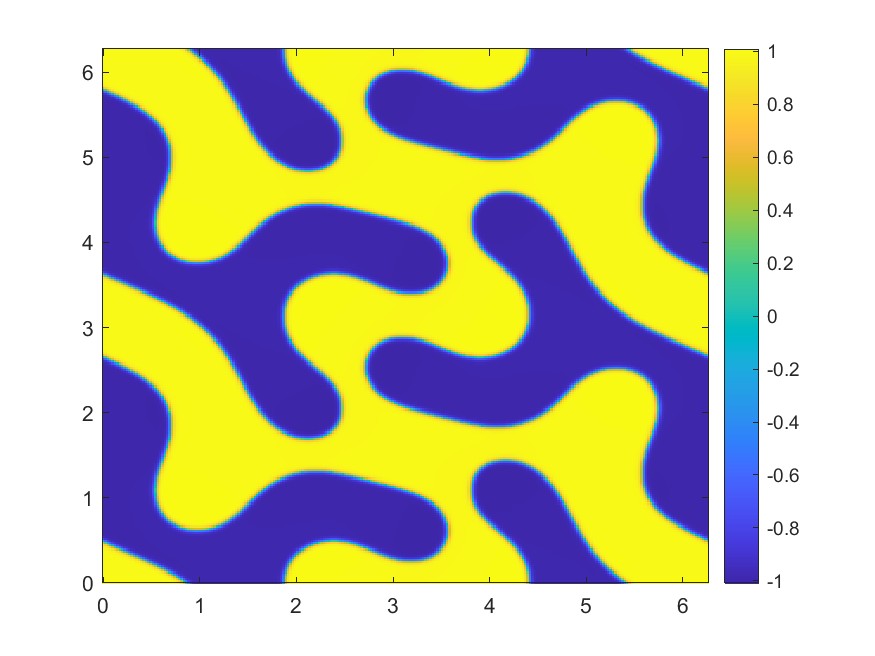}
  \caption{$t=8.0$}
\end{subfigure}
\caption{Numerical solution at different time stages with sinusoidal initial condition.}\label{cahn hilliard 2}
\end{figure}

In this section, we follow example 4 proposed in \cite{church2019high} to compute \eqref{CH equ} on $\Omega = [0, 2\pi]^2$. We set $a = \frac{\pi^2}{25000}$, $b=1$. The initial condition is set as
\[u_0(x,y) = 0.05(\cos(3x)\cos(4y)+(\cos(4x)\cos(3y))^2+\cos(x-5y)\cos(2x-y)).\]
We solve the equation on the time interval $[0, 8]$. We set $N_x = 256$, $h_x = \pi/128$; $N_t=24000$, $h_t = 1/3000$. For the PDHG part, we choose $\tau_u = \tau_p = 0.5$. The PDHG iteration is working efficiently at every time stepsize. Some numerical plots are shown in Figure \ref{cahn hilliard 2}.

\subsubsection{Example with random initial condition}
\begin{figure}[htb!]
\begin{subfigure}{.19\textwidth}
  \centering
  \includegraphics[width=\linewidth]{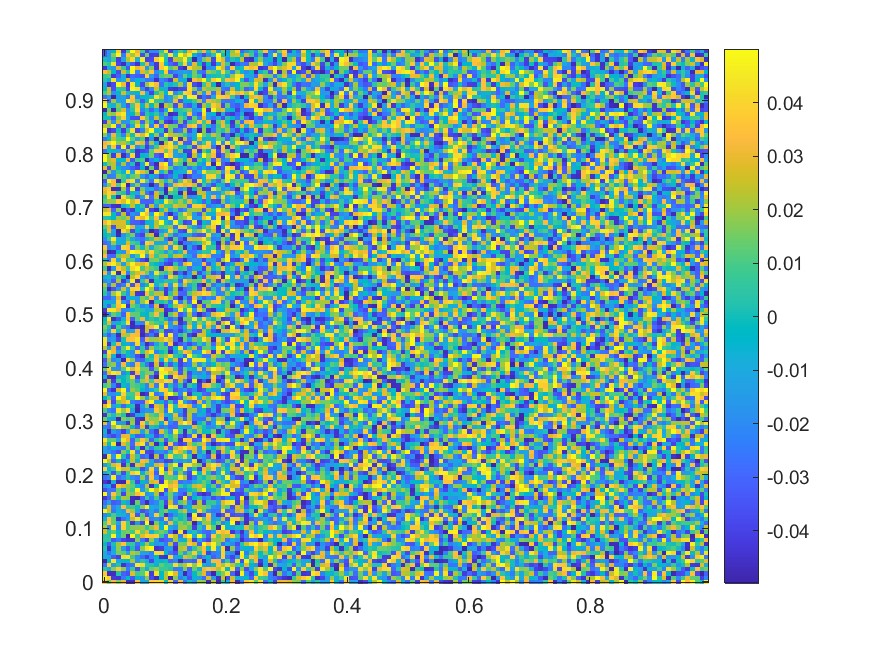}
  \caption{$t=0.0$}
\end{subfigure}
\begin{subfigure}{.19\textwidth}
  \centering
  \includegraphics[width=\linewidth]{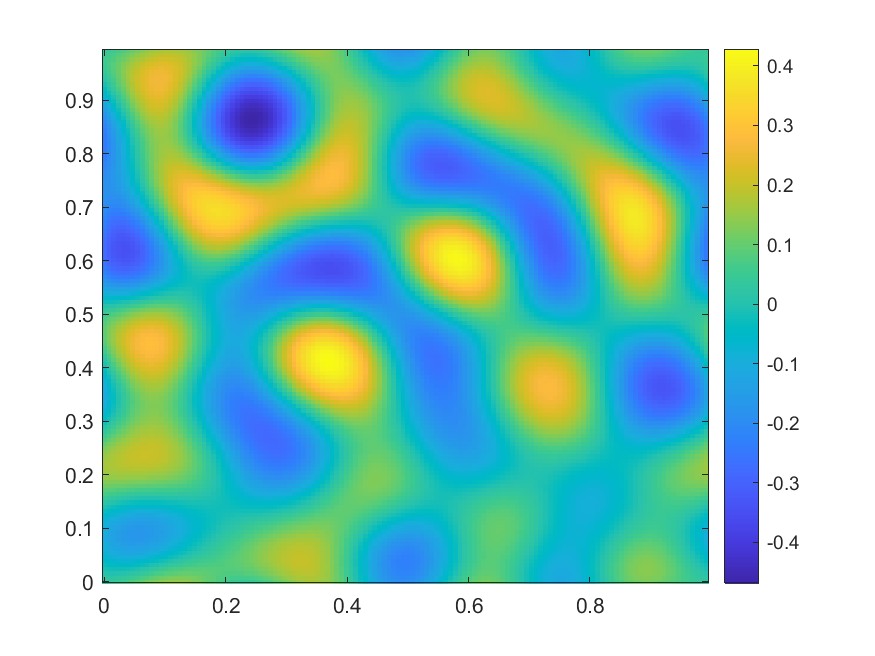}
  \caption{$t=0.001$}
\end{subfigure}
\begin{subfigure}{.19\textwidth}
  \centering
  \includegraphics[width=\linewidth]{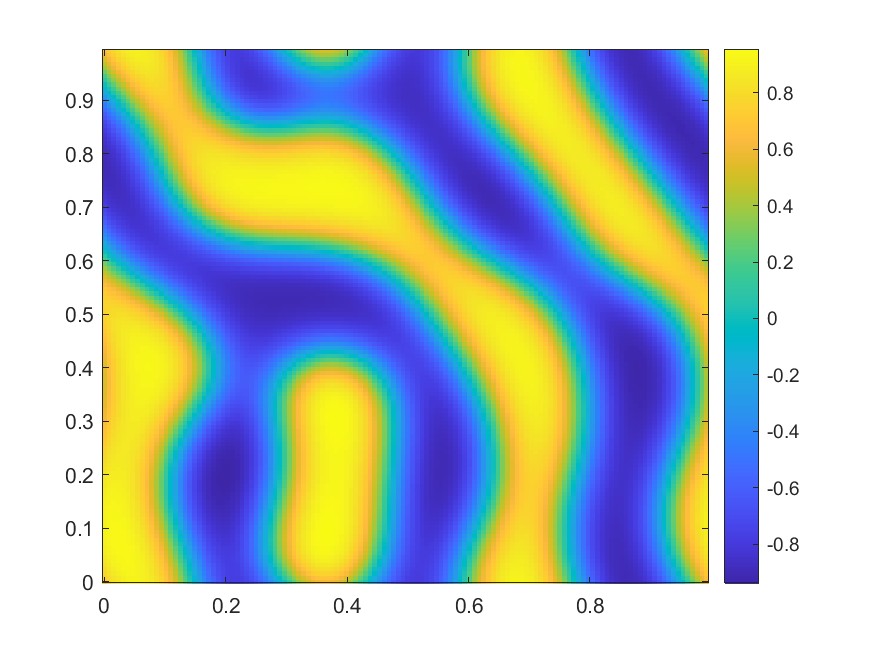}
  \caption{$t=0.003$}
\end{subfigure}
\begin{subfigure}{.19\textwidth}
  \centering
  \includegraphics[width=\linewidth]{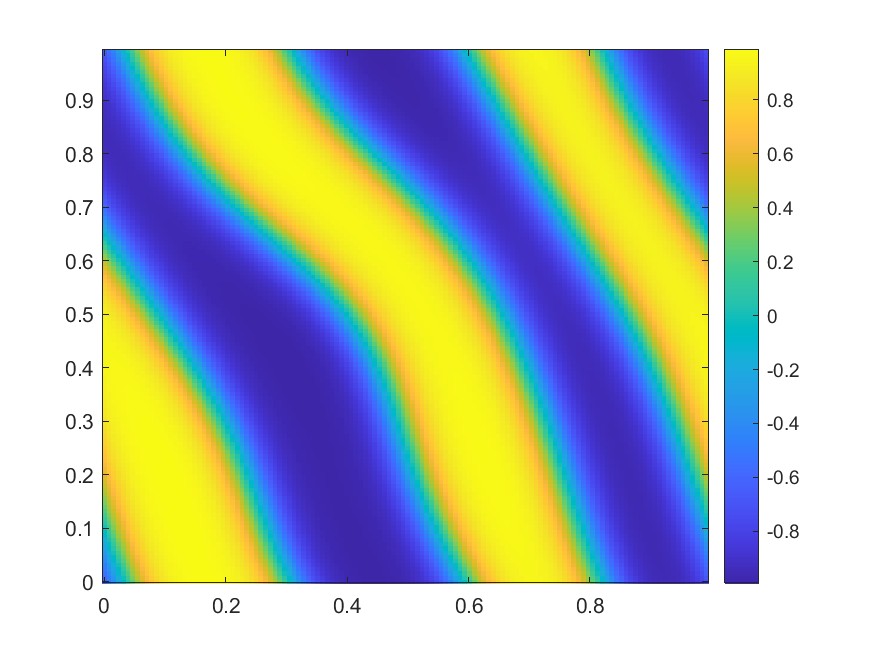}
  \caption{$t=0.01$}
\end{subfigure}
\begin{subfigure}{.19\textwidth}
  \centering
  \includegraphics[width=\linewidth]{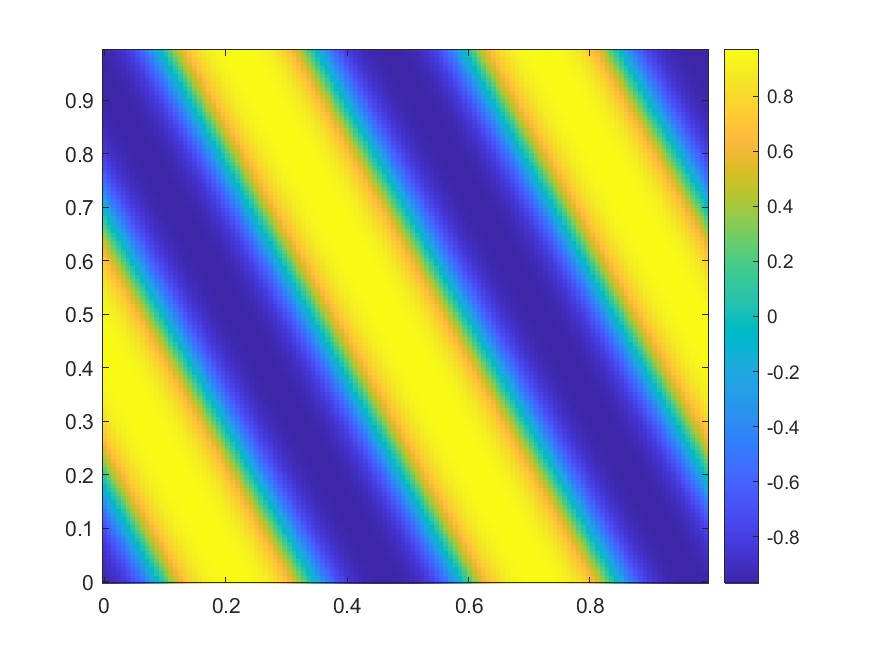}
  \caption{$t=1.0$}
\end{subfigure}
\caption{Numerical solution at different time stages with random initial condition.}\label{cahn hilliard 3}
\end{figure}

\begin{figure}[htb!]
\begin{subfigure}[t]{.45\textwidth}
  \centering
  \includegraphics[width=\linewidth]{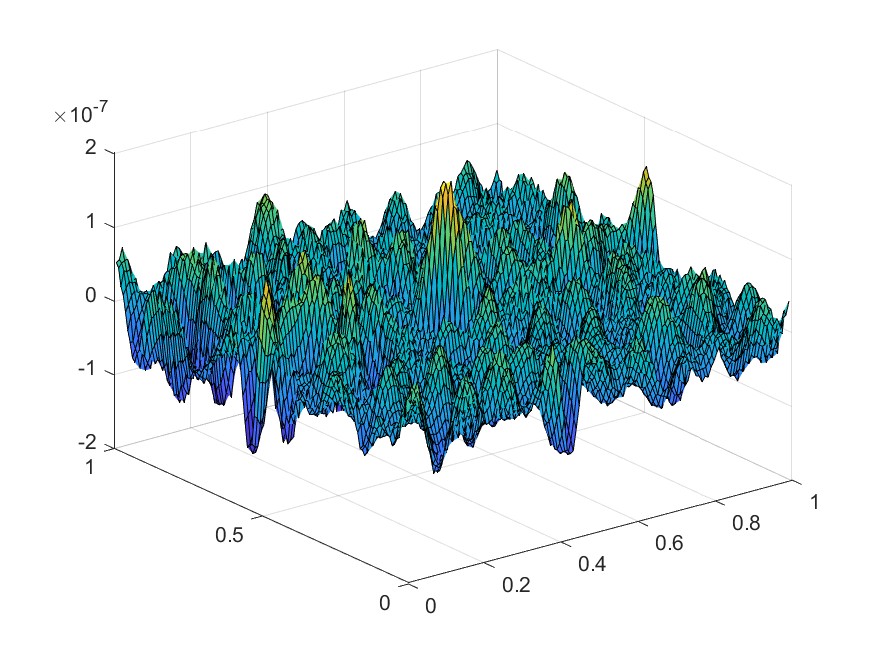}
  \caption{Plot of the residual $\mathrm{Res}(U)$ at $t=0.01$}
\end{subfigure}
\begin{subfigure}[t]{.45\textwidth}
  \centering
  \includegraphics[width=\linewidth]{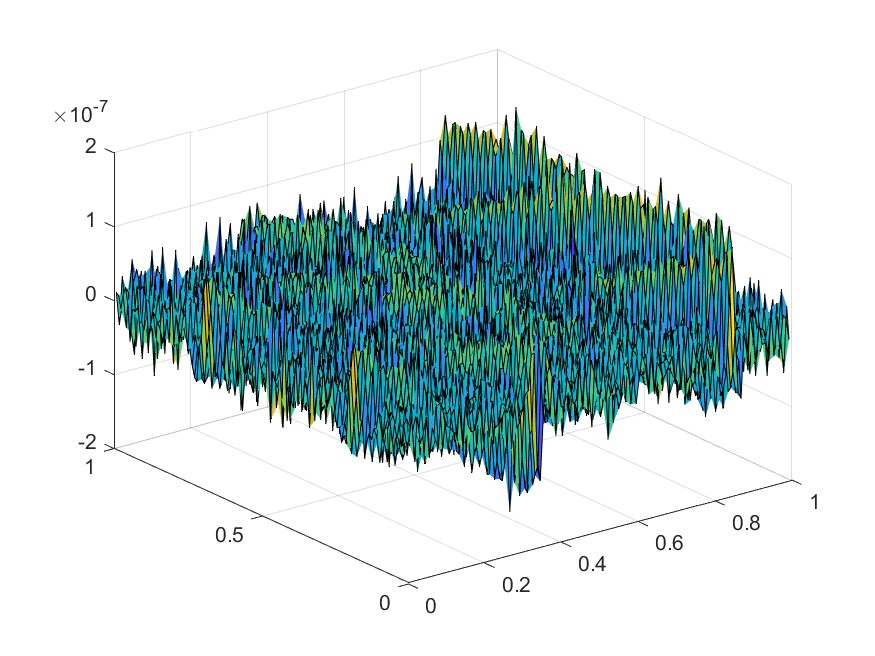}
  \caption{Plot of the residual $\mathrm{Res}(U)$ at $t=1.0$}
\end{subfigure}
\caption{Plots of residual at $t=0.01$ and $t=1.0$.}\label{cahn hilliard3 residue}
\end{figure}

One can also consider the Cahn-Hilliard equation \eqref{CH equ} with random initial condition. This may impose more challenges to our computation since the randomness will remove the regularity of $u_0$ and make the numerical computation unstable. We will solve the equation \eqref{CH equ} on $[0, 1]$ in this example. We let the periodic domain $\Omega = [0, 1]^2$. Then we choose $N_x=128$, $h_x = 1/128$; $N_t = 100000$ with $h_t = 1/100000$. We choose a rather small time step size in this example in order to guarantee the accuracy of our numerical solution. For the initial condition, we choose $u_0$ as a random scalar field that takes i.i.d. values uniformly distributed on $[-0.05, 0.05].$ We evolve the PDHG dynamic with stepsize $\tau_u = \tau_p = 0.75$. We plot the numerical solutions at certain time stages in Figure \ref{cahn hilliard 3}. The reaction-diffusion system reaches the equilibrium state at $t=1$.  The residual plots of $\mathrm{Res}(U)$ at time $t=0.01$ and $t=1$ are provided in Figure \ref{cahn hilliard3 residue}.

\subsection{Higher-order Reaction-Diffusion Equations}
In addition to the Allen-Cahn and Cahn-Hilliard equations, we test the method on the following 6th-order Cahn-Hilliard-type equation. 
\begin{equation}
  \frac{\partial u(x,t)}{\partial t} = \Delta (\epsilon^2\Delta - W''(u) + \epsilon^2) (\epsilon^2\Delta u - W'(u)) \textrm{ On } \Omega, \quad u(\cdot,0) = u_0.  \label{6thorder}
\end{equation}
The above equation was first proposed in \cite{polym4010630} which depicts the pore formation in functionalized polymers. This equation was later studied in the numerical examples of \cite{christlieb2014high}. In this example, we set $\Omega = [0, 2\pi]^2$. We choose parameter $\epsilon = 0.18$. The potential function $W(u)$ is the same as defined in \eqref{AC CH potential}. Thus $W'(u)=u^3-u, W''(u)=3u^2$. Similar to Allen-Cahn or Cahn-Hilliard equations, equation \eqref{6thorder} can also be treated as a flow that dissipates the energy $\mathcal{E}(u)$ with $a=\epsilon^2, b=1$. 

In our numerical implementation of the PDHG method, the functional $F(U)$ is now
\begin{equation*}
  F(U) = U - h_t \PLap(\epsilon^2\PLap - \mathrm{diag}(W''(U)) + \epsilon^2 I)(\epsilon^2 \PLap U - W'(U)) - U^k.
\end{equation*}
Similar to Allen-Cahn and Cahn-Hilliard equations, we pick the preconditioner $G$ as the square of the matrix in the dominating linear part of $F(U)$. However, if we directly keep the diagonal matrix $\mathrm{diag}(W''(U))$, we will not be able to invert $G$ efficiently by using FFT. Since the value of $W''(u)$ close to the equilibrium states $\pm 1$ is approximately $2$, we follow a similar idea in \cite{christlieb2014high} to replace such matrix with $2I$. Thus, in this problem, we set $G=(I - h_t \epsilon^2 \PLap(\epsilon^2\PLap - (2-\epsilon^2)I)\PLap)^2$ which can be inverted via FFT algorithm. Now the 3-line PDHG update is formulated as
\begin{align*}
  P_{n+1} & =  P_n + \tau_p G^{-1}( U_n - h_t \PLap\widetilde{\mathrm{Lap}}_{h_x}(\epsilon, U_n))(\epsilon^2 \PLap U_n - W'(U_n)) - U^k); \\ 
  \widetilde{P}_{n+1} & = P_{n+1} + \omega (P_{n+1} - P_n); \\   
  U_{n+1} & = U_n - \tau_u(\widetilde{P}_{n+1} - h_t  (\epsilon^2 \PLap - \mathrm{diag}(W''(U))) \widetilde{\mathrm{Lap}}_{h_x}(\epsilon, U_n)\PLap \tilde{P}_{n+1} \\ 
  & \qquad \qquad \qquad \quad + h_t  (\epsilon^2 \PLap U_n - W'(U_n)) \odot W'''(U_n) \odot \PLap \tilde{P}_{n+1})).
\end{align*}
Here we denote $\widetilde{\mathrm{Lap}}_{h_x}(\epsilon, U) = \epsilon^2\PLap - \mathrm{diag}(W''(U)) + \epsilon^2 I$. If the size of $U$ is $N_x^2$, one can verify that all calculations among the PDHG iteration can be computed with complexity $O(N_x^2\log(N_x))$ via the FFT method. 

\begin{figure}[htb!]
\begin{subfigure}{.24\textwidth}
  \centering
  \includegraphics[width=\linewidth]{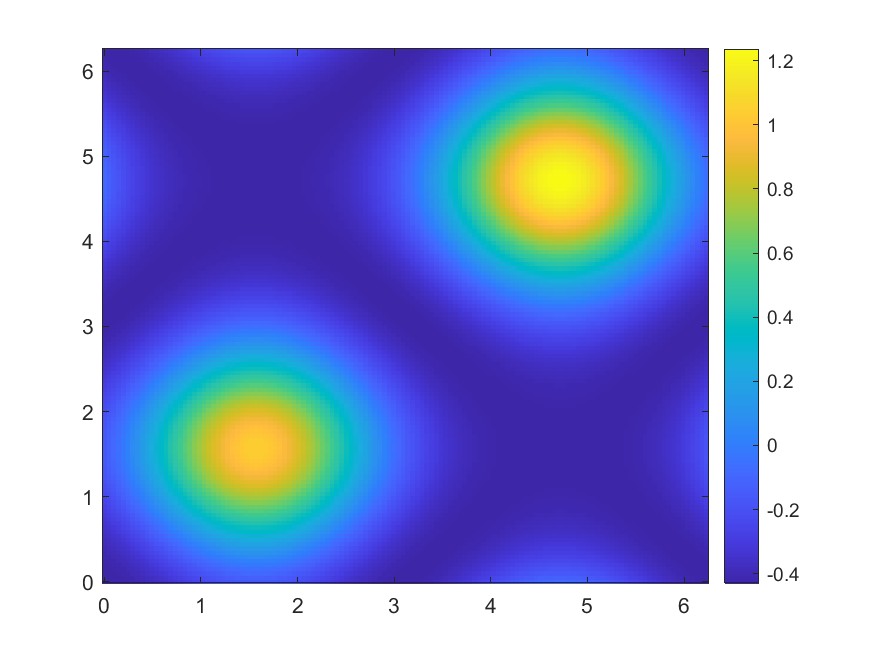}
  \caption{$t=0.0$}
\end{subfigure}
\begin{subfigure}{.24\textwidth}
  \centering
  \includegraphics[width=\linewidth]{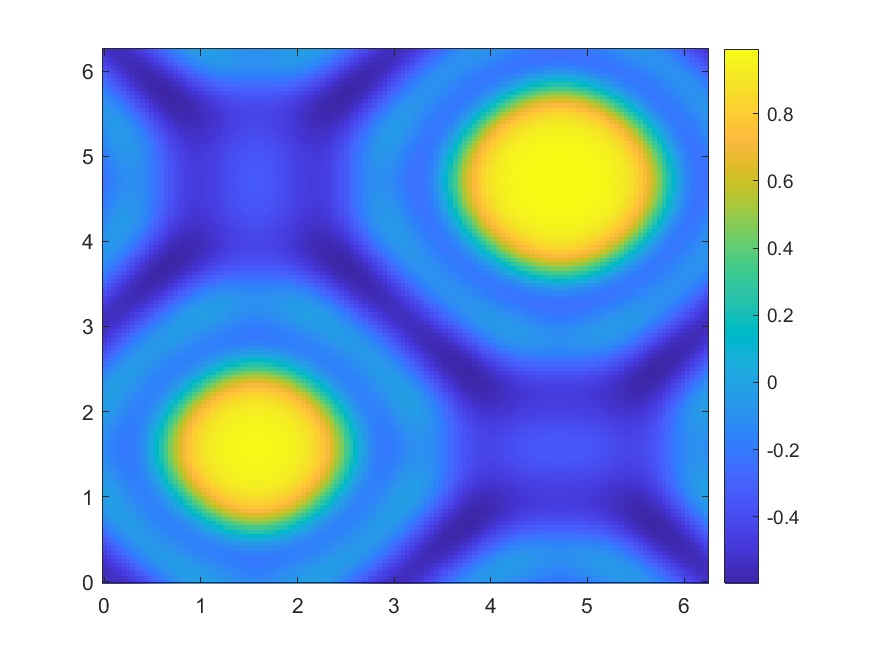}
  \caption{$t=0.1$}
\end{subfigure}
\begin{subfigure}{.24\textwidth}
  \centering
  \includegraphics[width=\linewidth]{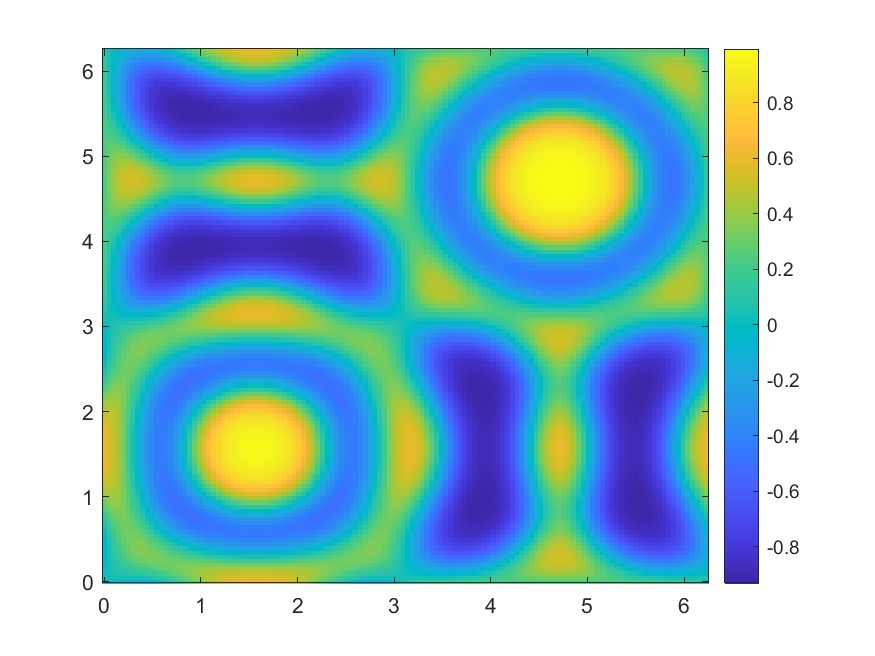}
  \caption{$t=2.0$}
\end{subfigure}
\begin{subfigure}{.24\textwidth}
  \centering
  \includegraphics[width=\linewidth]{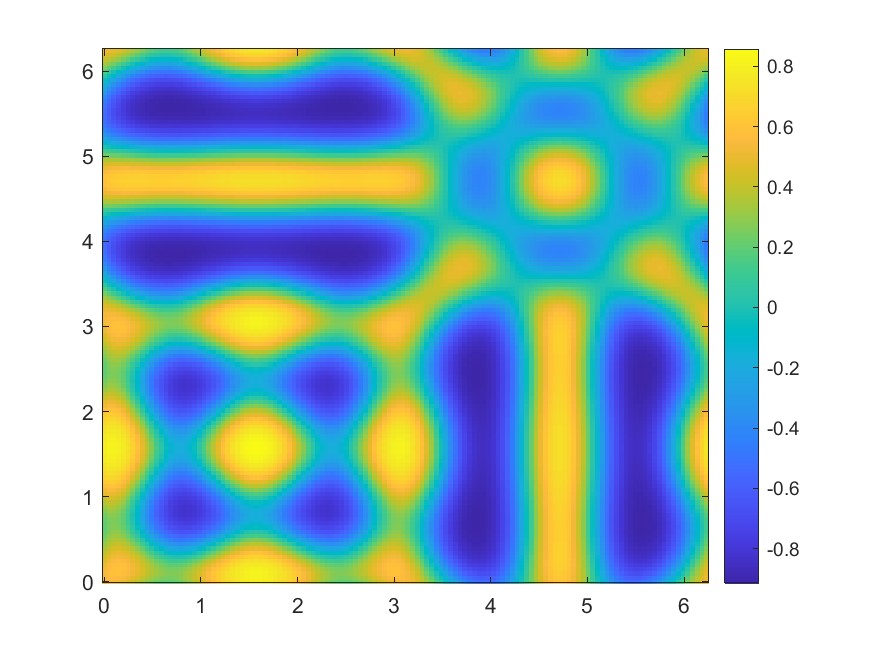}
  \caption{$t=20.0$}
\end{subfigure}
\caption{Numerical solution at different time stages with initial condition \eqref{6th order initial}.}\label{6th order 1}
\end{figure}

\begin{figure}[htb!]
\begin{subfigure}{.32\textwidth}
  \centering
  \includegraphics[width=\linewidth]{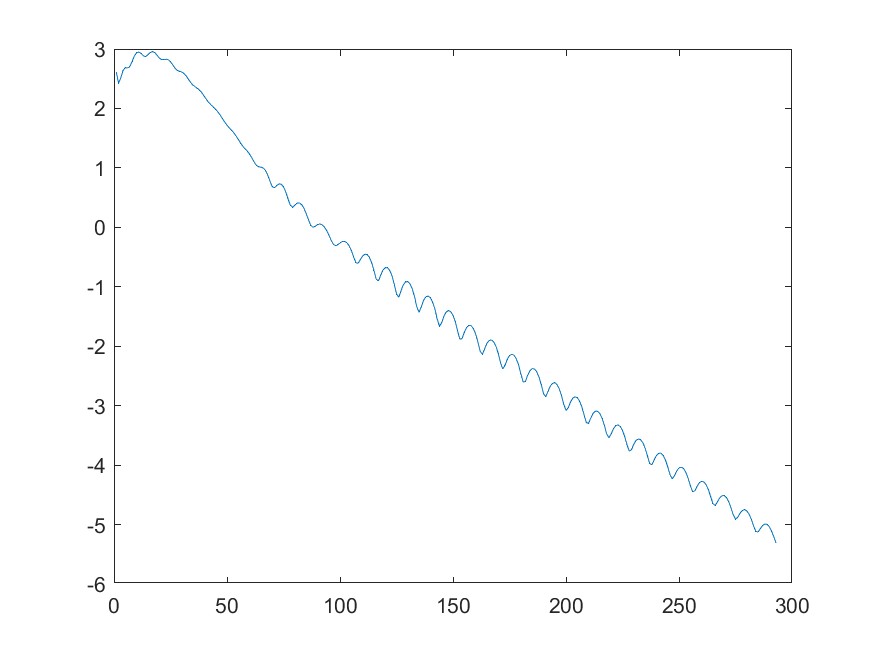}
  \caption{Residual decay at $t=0.1$}
\end{subfigure}
\begin{subfigure}{.32\textwidth}
  \centering
  \includegraphics[width=\linewidth]{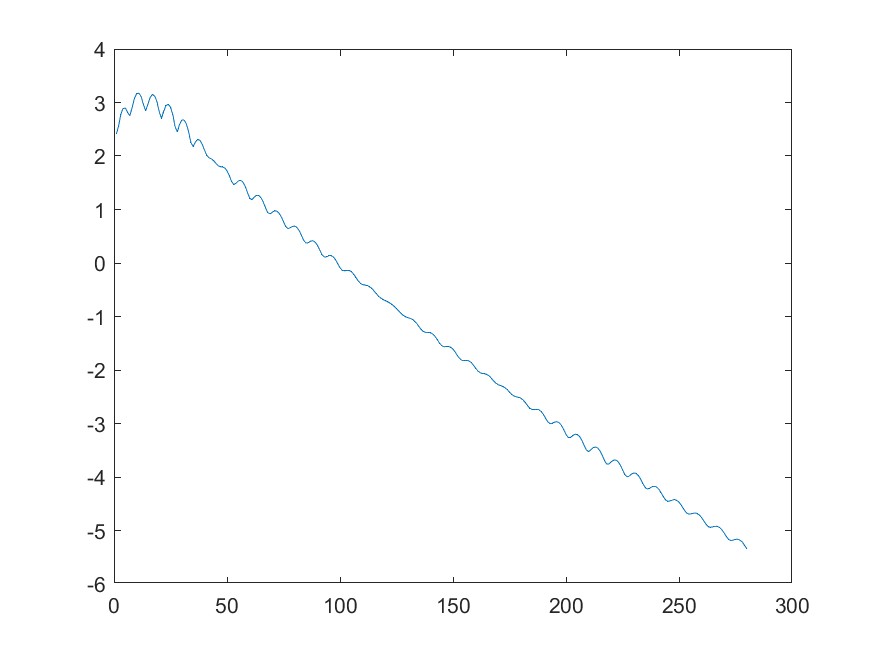}
  \caption{Residual decay at $t=2.0$}
\end{subfigure}
\begin{subfigure}{.32\textwidth}
  \centering
  \includegraphics[width=\linewidth]{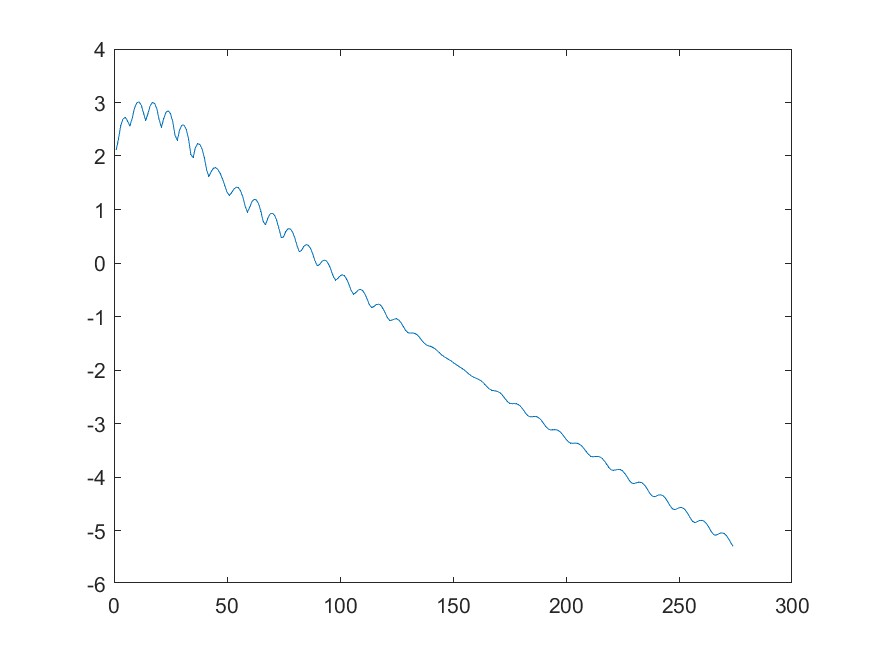}
  \caption{Residual decay at $t=20.0$}
\end{subfigure}
\begin{subfigure}{.32\textwidth}
  \centering
  \includegraphics[width=\linewidth]{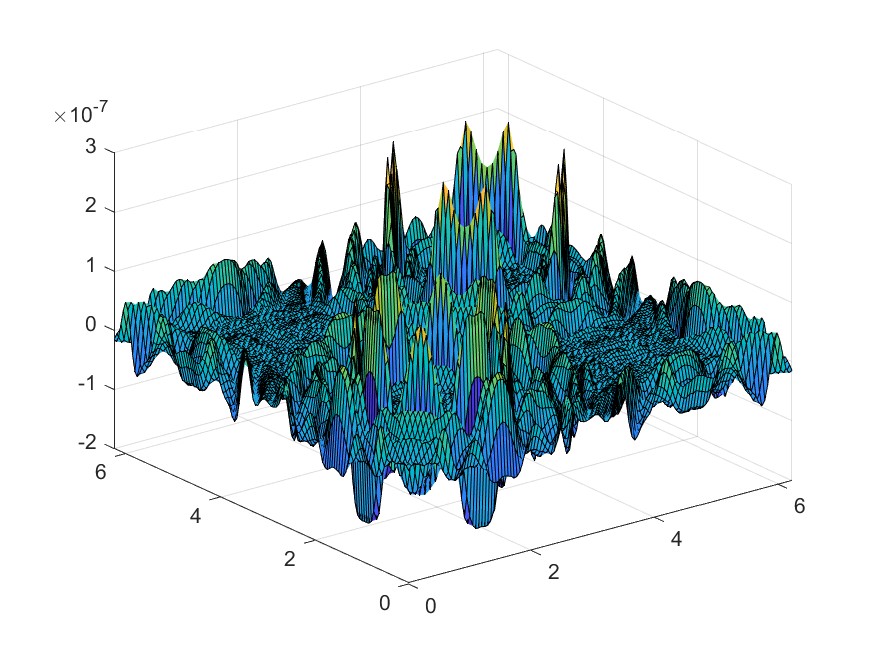}
  \caption{Residual plot at $t=0.1$}
\end{subfigure}
\begin{subfigure}{.32\textwidth}
  \centering
  \includegraphics[width=\linewidth]{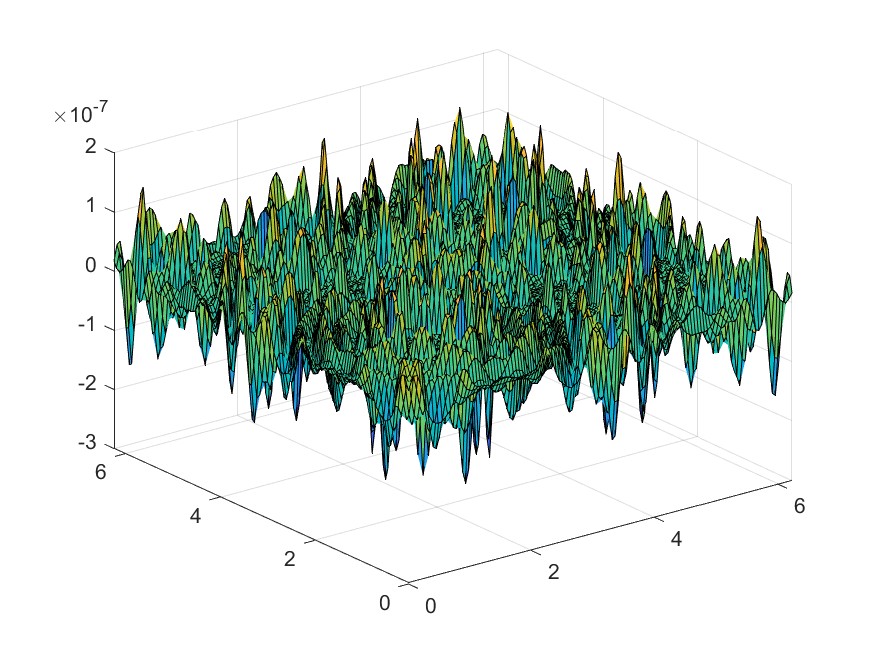}
  \caption{Residual plot at $t=2.0$}
\end{subfigure}
\begin{subfigure}{.32\textwidth}
  \centering
  \includegraphics[width=\linewidth]{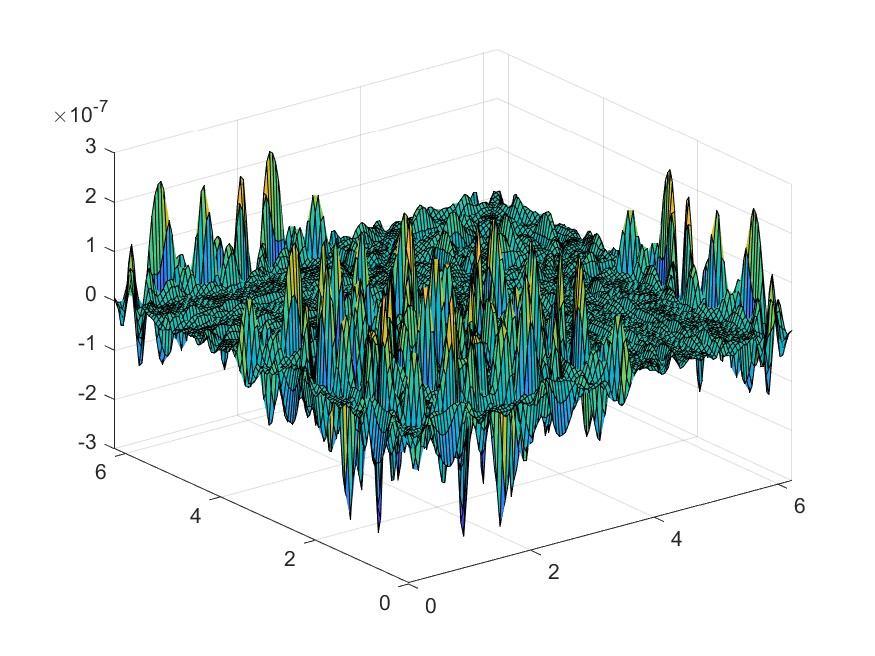}
  \caption{Residual plot at $t=20.0$}
\end{subfigure}
\caption{$\log-$residual decay \& plots of residual functional $\mathrm{Res}(U^n)$ at different time stages $t=0.1,2.0,20.0$.}\label{6th order 2}
\end{figure}

Similar to \cite{christlieb2014high}, we choose initial condition 
\begin{equation}
  u_0(x,y) = 2 e^{\sin x+\sin y-2} + 2.2 e^{-\sin x-\sin y-2}-1.  \label{6th order initial}
\end{equation}
In the numerical implementation, we solve the equation \eqref{6thorder} from $t=0$ to $t=20$. We choose $N_x = 128$, $h_x = \pi / 64$; $N_t = 20000$, $h_t = 1/1000$. We choose the PDHG stepsizes $\tau_u = \tau_p = 0.58$. We choose the threshold for terminating the iteration as $\delta = 0.5 \times 10^{-5}$. We present some of the results in Figure \ref{6th order 1} and Figure \ref{6th order 2}.

\subsection{Reaction-diffusion systems}
We have already shown some reaction-diffusion equation examples in the previous sections. We now apply the method to compute reaction-diffusion systems. 

\subsubsection{Schnakenberg Model}
The Schnakenberg model is first considered in \cite{schnakenberg1979simple} to model the limit-cycle behavior in a two-component chemical reaction system. In the discussion, we consider the following reaction-diffusion PDE system defined on unit square $\Omega = [0,1]^2$ where $u$, $v$ represent the density concentration of two chemicals. Such a PDE system is also investigated in references \cite{hundsdorfer2003numerical, zhu2009application}.
\begin{align}
  \frac{\partial u(x,y,t)}{\partial t} & = D_1\Delta u(x,y,t) + \kappa(a - u + u^2v), \label{Schnakenberg system 1}\\
  \frac{\partial v(x,y,t)}{\partial t} & = D_2\Delta v(x,y,t) + \kappa(b - u^2v). \label{Schnakenberg system 2}
\end{align}
The initial condition of the system is
\begin{equation}
  u(x,y,0) = a + b + 10^{-3}*e^{-100((x-\frac{1}{3})^2+(y-\frac{1}{2})^2)}, \quad v(x,y,0) = \frac{b}{(a+b)^2}.\label{initial condition Schnakenberg system}
\end{equation}
Here we set $\kappa = 100, a = 0.1305, b = 0.7695, D_1 = 0.05, D_2 = 1$.
One can understand the initial data as exerting a tiny perturbation to the equilibrium solution $(a+b, \frac{b}{(a+b)^2})$ of the Schnakenberg system \eqref{Schnakenberg system 1}, \eqref{Schnakenberg system 2}. Such equilibrium state is unstable, the small perturbation will lead to the formation of certain dotted patterns in both components $u$ and $v$.

We assume the Neumann boundary condition $\frac{\partial u}{\partial \boldsymbol{n}} = 0$ on $\partial \Omega$ where $\frac{\partial}{\partial \boldsymbol{n}}$ denotes the directional derivative w.r.t. the outer pointing normal direction $\boldsymbol{n}$.

Suppose we apply the one-step implicit scheme to solve this problem, recall the discrete Laplacian with Neumann BC introduced in \eqref{def Lap N}, at the $k-$th time step, we consider 
\begin{align*}
  & F_u(U, V) = U - U^k - h_t(D_1\NLap U + \kappa (a\boldsymbol{1} - U + U^2\odot V)); \\
  & F_v(U, V) = V - V^k - h_t(D_2\NLap V + \kappa (b\boldsymbol{1} - U^2\odot V)).
\end{align*}
At each time step, our purpose is to solve $F_u(U,V)=0, F_v(U,V)=0$ for updating $U^{k}, V^{k}$. By treating $\widetilde{U} = (U, V)\in \mathbb{R}^{2N_x^2}$ as an entity; and by denoting $\widetilde{F}: \mathbb{R}^{2N_x^2}\rightarrow \mathbb{R}^{2N_x^2}, \widetilde{U} \mapsto (F_u(U, V)^\top, F_v(U,V)^\top)^\top$, the problem of solving $\widetilde{F}(\widetilde{U})=0$ reduces to the scenario of solving single $F(U)=0$ discussed before. Hence, it is natural to introduce the dual variable $\widetilde{P} = (P, Q) \in  \mathbb{R}^{2N_x^2}$; The stiff Laplacian terms can be treated as dominating linear terms of both functions ${F}, {G}$, thus we set our preconditioner matrix $\widetilde{G} = \left(\begin{array}{cc}
    G_u &  \\
     & G_v
\end{array}\right)$ with $G_u = (I-h_t D_1 \NLap)^2$ and $G_v = (I - h_t D_2 \NLap)^2$. The corresponding PDHG iteration for solving $\widetilde{F}(\widetilde{U})=0$ is formulated as follows.
\begin{align*}
  P_{n+1} & =  P_n + \tau_p G_u^{-1}(U_n - U^k - h_t(D_1\NLap U_n + \kappa (a\boldsymbol{1} - U_n + U_n^2\odot V_n))); \\
  Q_{n+1} & =  Q_n + \tau_p G_v^{-1}(V_n - V^k - h_t(D_2\NLap V_n + \kappa (b\boldsymbol{1} - U_n^2\odot V_n))); \\
  \widetilde{P}_{n+1} & = P_{n+1} + \omega (P_{n+1} - P_n); \quad \widetilde{Q}_{n+1} = Q_{n+1} + \omega (Q_{n+1} - Q_n);\\
  U_{n+1} & = U_n - \tau_u(\widetilde{P}_{n+1} - h_t( D_1\NLap\widetilde{P}_{n+1} + \kappa (-\widetilde{P}_{n+1} + 2 U_n \odot V_n \odot (\widetilde{P}_{n+1}-\widetilde{Q}_{n+1}))).\\
  V_{n+1} & = V_n - \tau_u(\widetilde{Q}_{n+1} - h_t( D_2\NLap\widetilde{Q}_{n+1} + \kappa (U_n^2 \odot (\widetilde{P}_{n+1} - \widetilde{Q}_{n+1}))). \\  
\end{align*}
We recall that the Discrete Cosine Transform mentioned in Remark \ref{remark: NLap and DCT} can be used to compute matrix-vector multiplication involving $\NLap$ as well as inverting the preconditioners $G_u, G_v$ within $O(N_x^2\log N_x)$ complexity. Thus every step of the above PDHG iterations can be computed efficiently.

In our numerical implementation, we solve this PDE system, on time interval $[0, 2]$. We choose $N_x = 128$, $h_x = 1/128$; and $N_t = 10000$, $h_t = 1/5000$. We then choose $\tau_u=\tau_p = 0.9$. We terminate the PDHG iteration when $\|\mathrm{Res}(U_n)\|_2+\|\mathrm{Res}(V_n)\|_2<\delta$, where we pick threshold $\delta = 10^{-7}$. Our numerical solutions are presented in the following Figure \ref{Schnakenberg 1}.

\begin{figure}[htb!]
\begin{subfigure}{.162\textwidth}
  \centering
  \includegraphics[width=\linewidth]{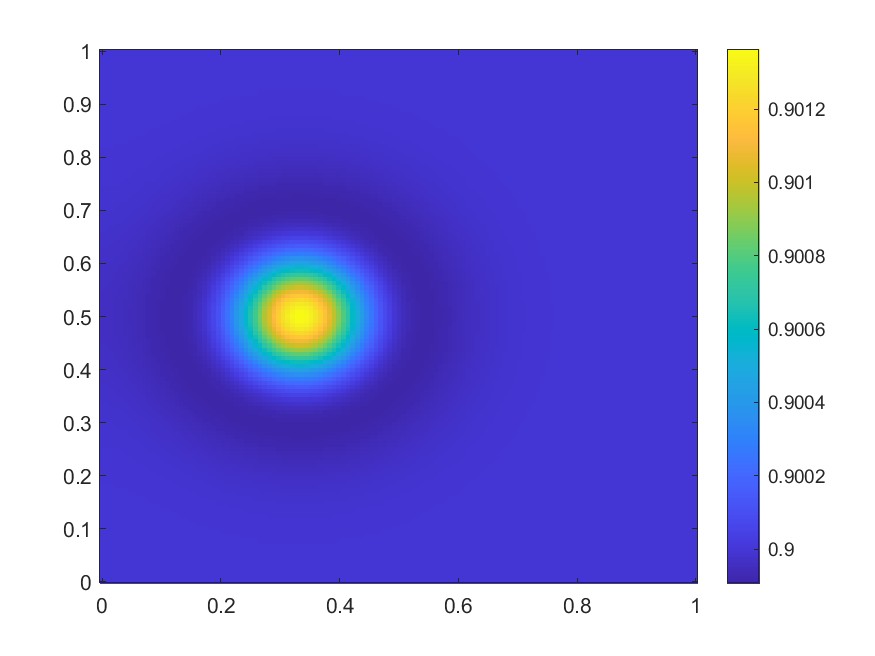}
\end{subfigure}
\begin{subfigure}{.162\textwidth}
  \centering
  \includegraphics[width=\linewidth]{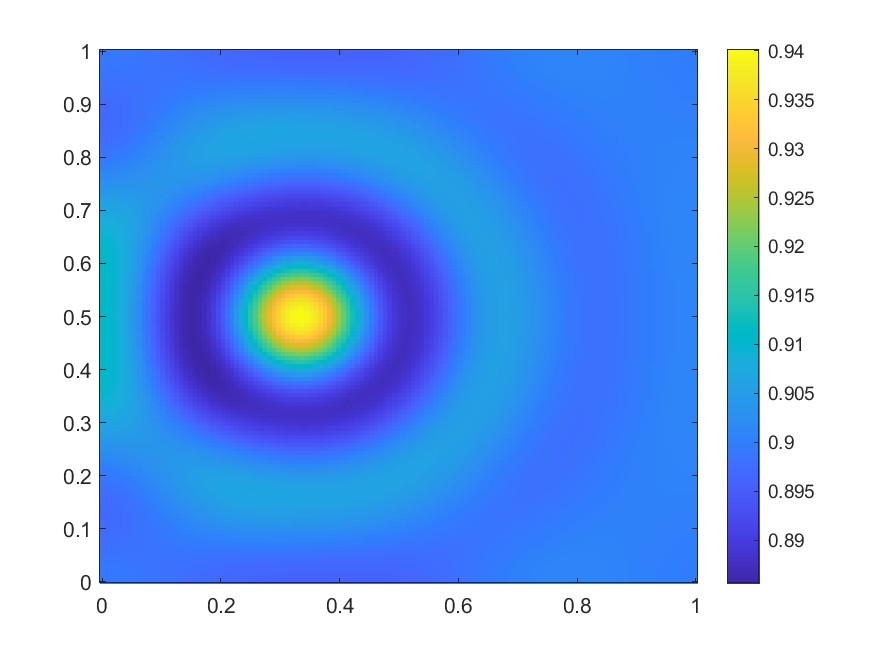}
\end{subfigure}
\begin{subfigure}{.162\textwidth}
  \centering
  \includegraphics[width=\linewidth]{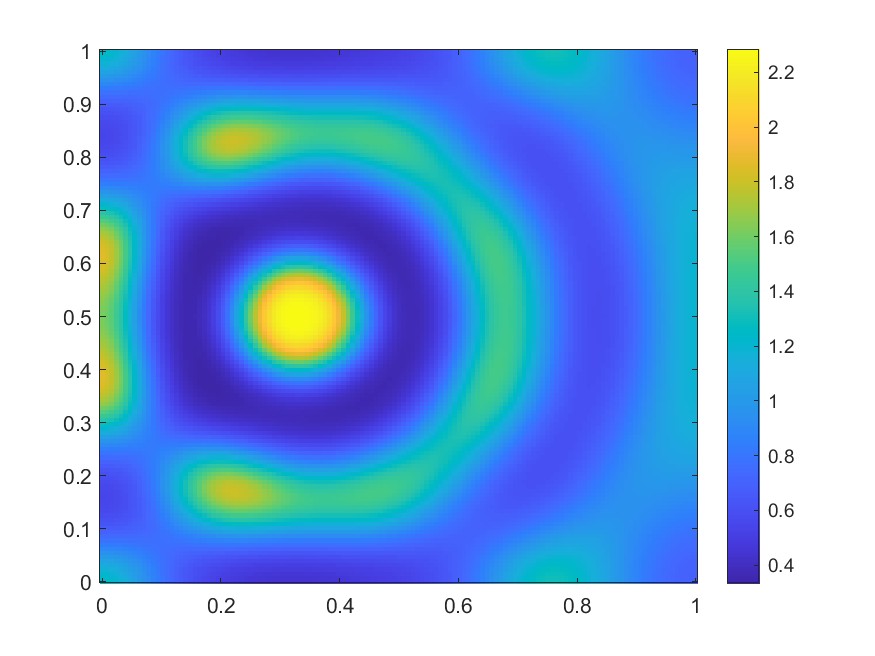}
\end{subfigure}
\begin{subfigure}{.162\textwidth}
  \centering
  \includegraphics[width=\linewidth]{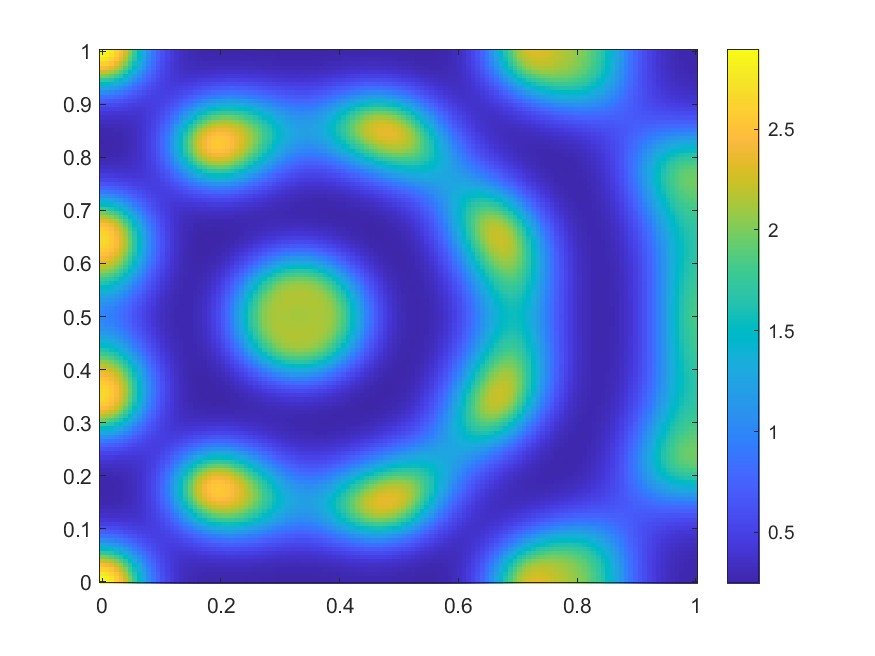}
\end{subfigure}
\begin{subfigure}{.162\textwidth}
  \centering
  \includegraphics[width=\linewidth]{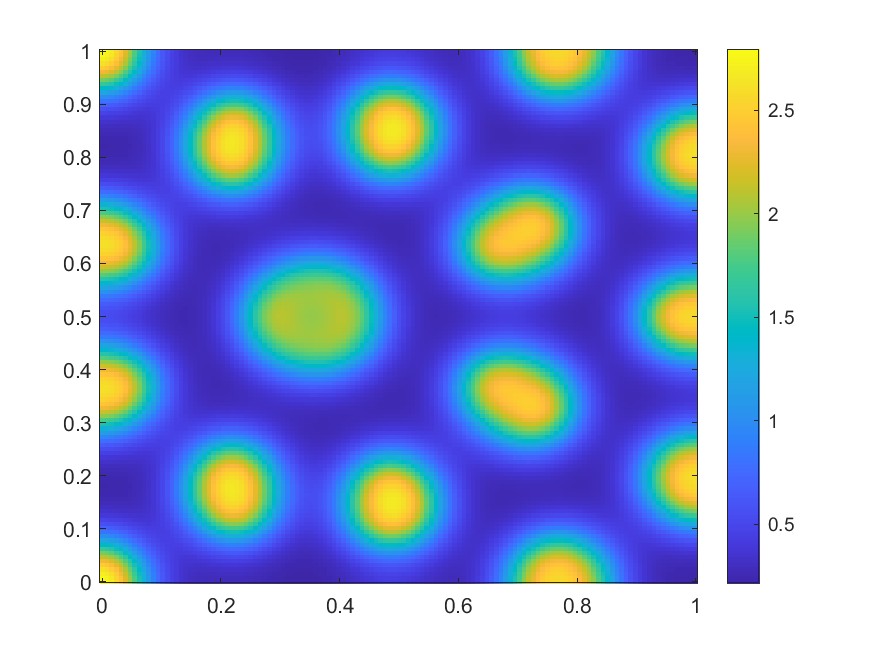}
\end{subfigure}
\begin{subfigure}{.162\textwidth}
  \centering
  \includegraphics[width=\linewidth]{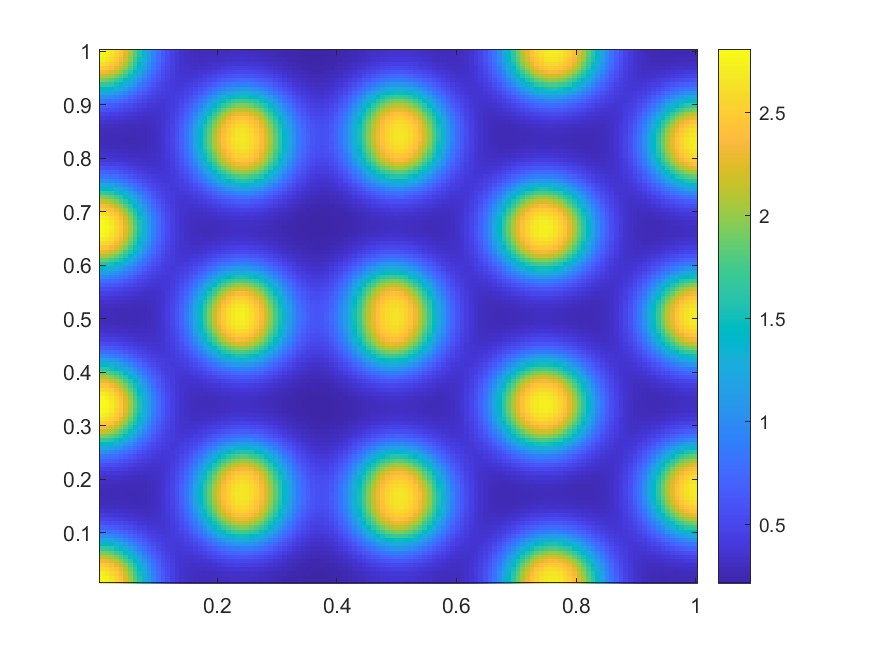}
\end{subfigure}
\begin{subfigure}{.162\textwidth}
  \centering
  \includegraphics[width=\linewidth]{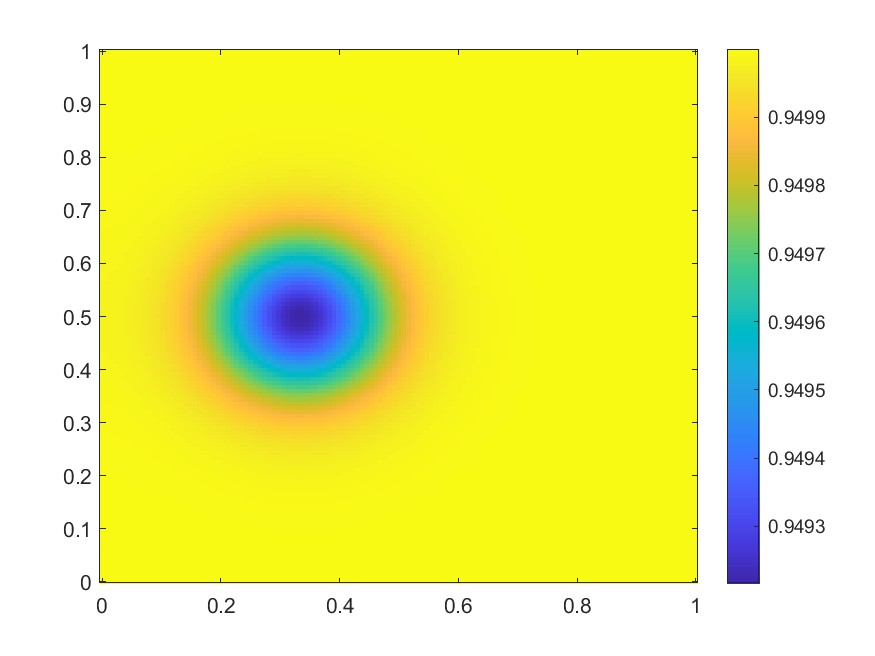}
  \caption{$t=0.0$}
\end{subfigure}
\begin{subfigure}{.162\textwidth}
  \centering
  \includegraphics[width=\linewidth]{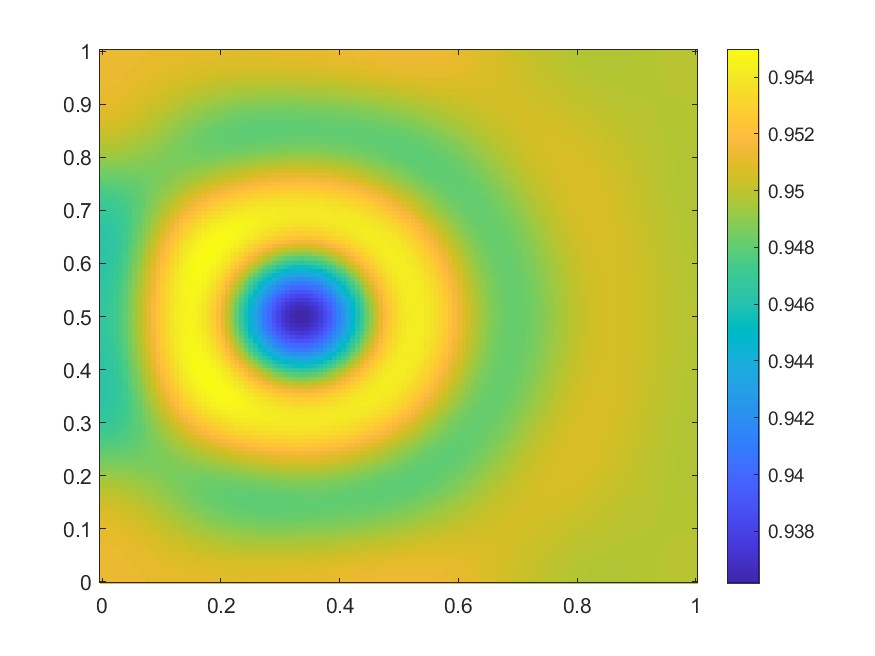}
  \caption{$t=0.2$}
\end{subfigure}
\begin{subfigure}{.162\textwidth}
  \centering
  \includegraphics[width=\linewidth]{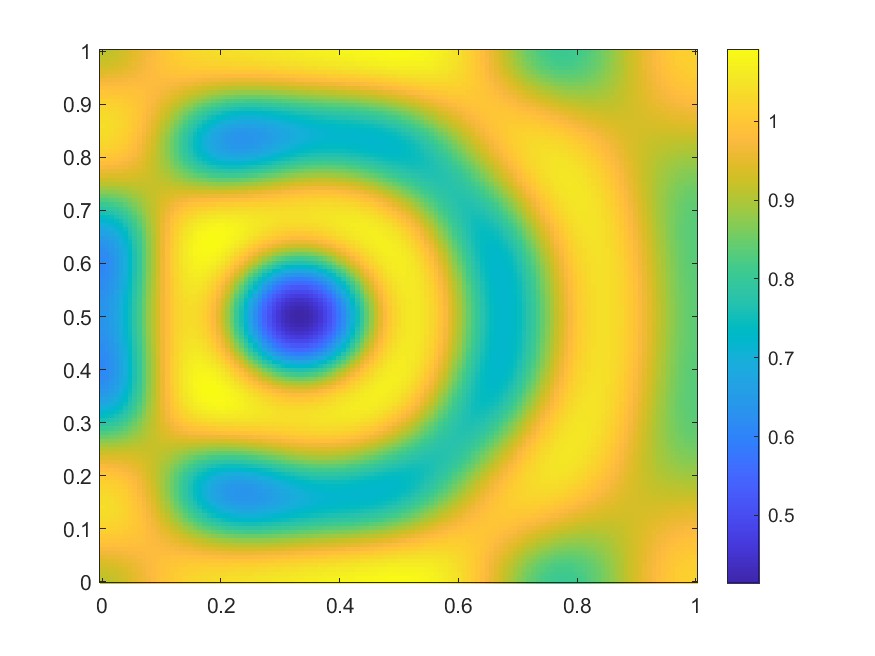}
  \caption{$t=0.4$}
\end{subfigure}
\begin{subfigure}{.162\textwidth}
  \centering
  \includegraphics[width=\linewidth]{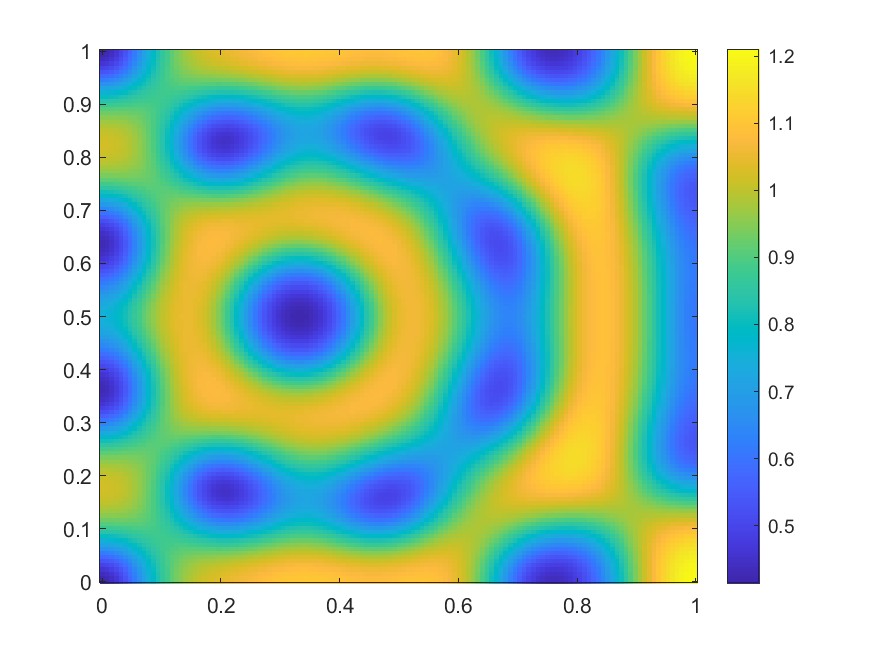}
  \caption{$t=0.5$}
\end{subfigure}
\begin{subfigure}{.162\textwidth}
  \centering
  \includegraphics[width=\linewidth]{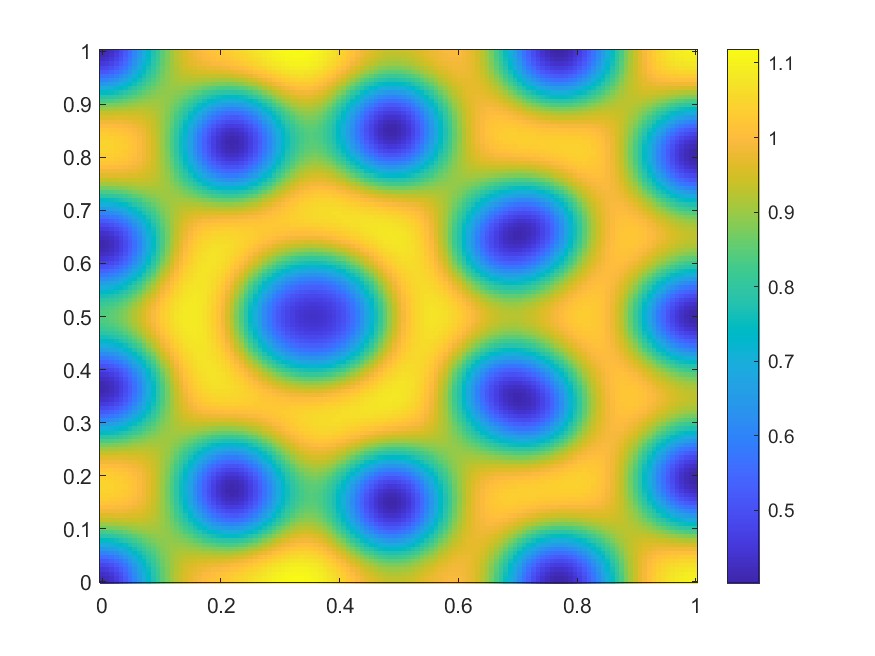}
  \caption{$t=1.0$}
\end{subfigure}
\begin{subfigure}{.162\textwidth}
  \centering
  \includegraphics[width=\linewidth]{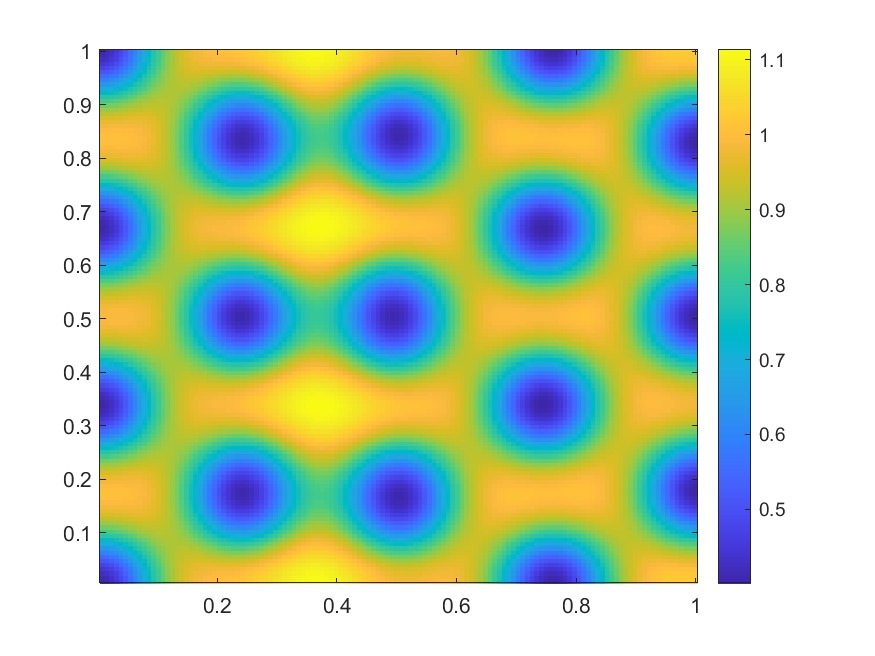}
  \caption{$t=2.0$}
\end{subfigure}
\caption{Numerical solution of $u$ (upper row), and $v$ (lower row) at different time stages with initial condition \eqref{initial condition Schnakenberg system}.}\label{Schnakenberg 1}
\end{figure}

\begin{figure}[htb!]
\begin{subfigure}{.24\textwidth}
  \centering
  \includegraphics[width=\linewidth]{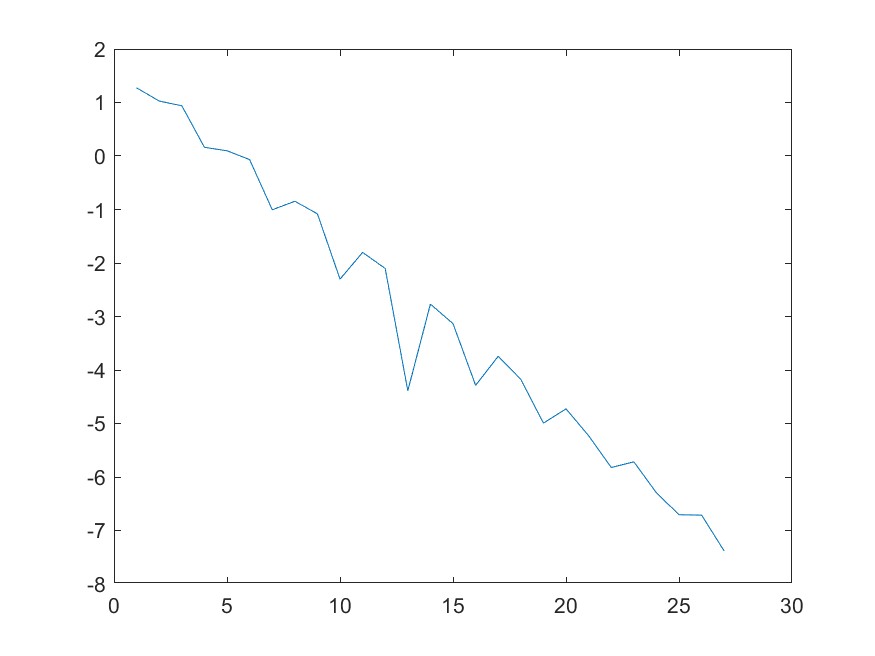}
\end{subfigure}
\begin{subfigure}{.24\textwidth}
  \centering
  \includegraphics[width=\linewidth]{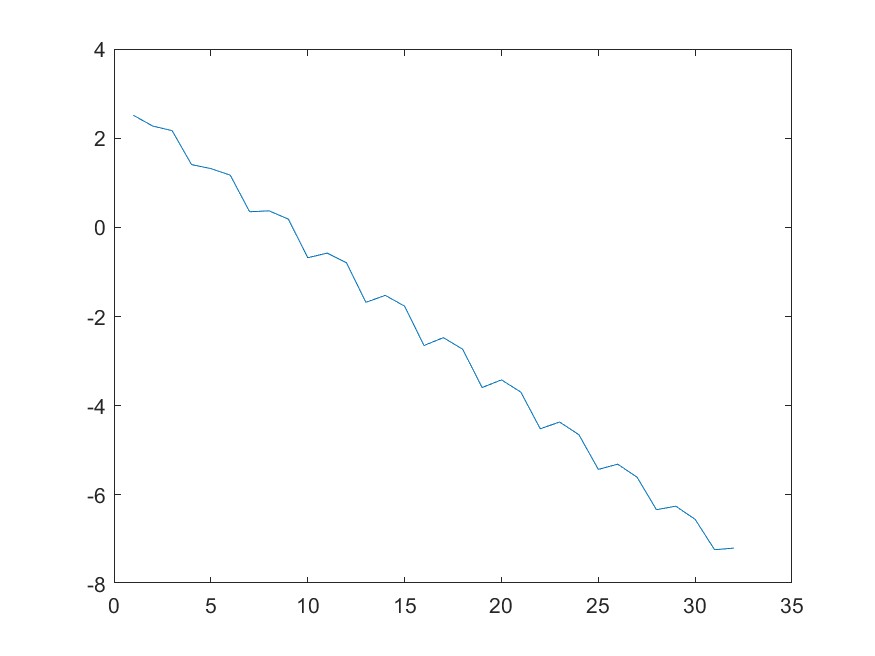}
\end{subfigure}
\begin{subfigure}{.24\textwidth}
  \centering
  \includegraphics[width=\linewidth]{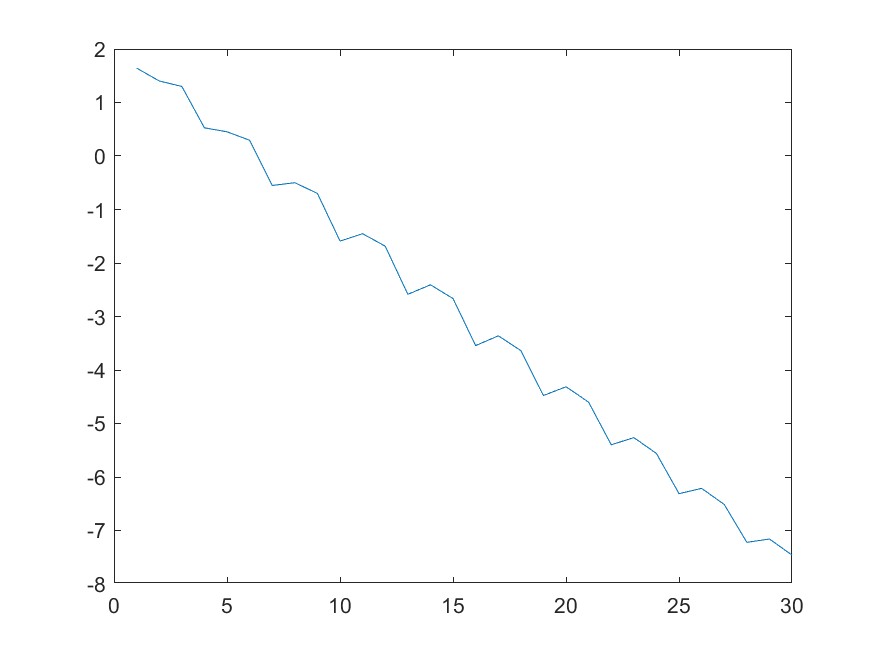}
\end{subfigure}
\begin{subfigure}{.24\textwidth}
  \centering
  \includegraphics[width=\linewidth]{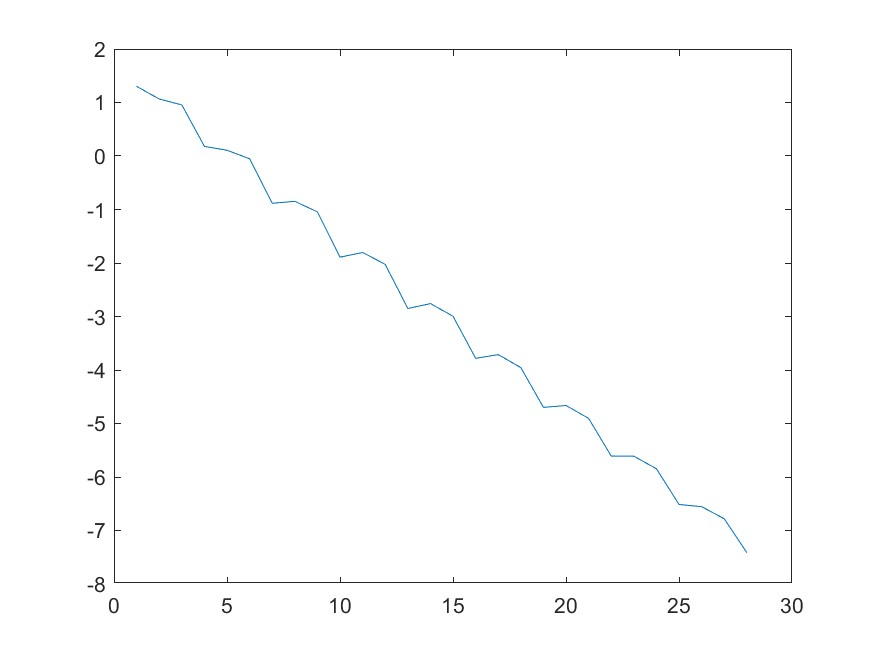}
\end{subfigure}
\begin{subfigure}{.24\textwidth}
  \centering
  \includegraphics[width=\linewidth]{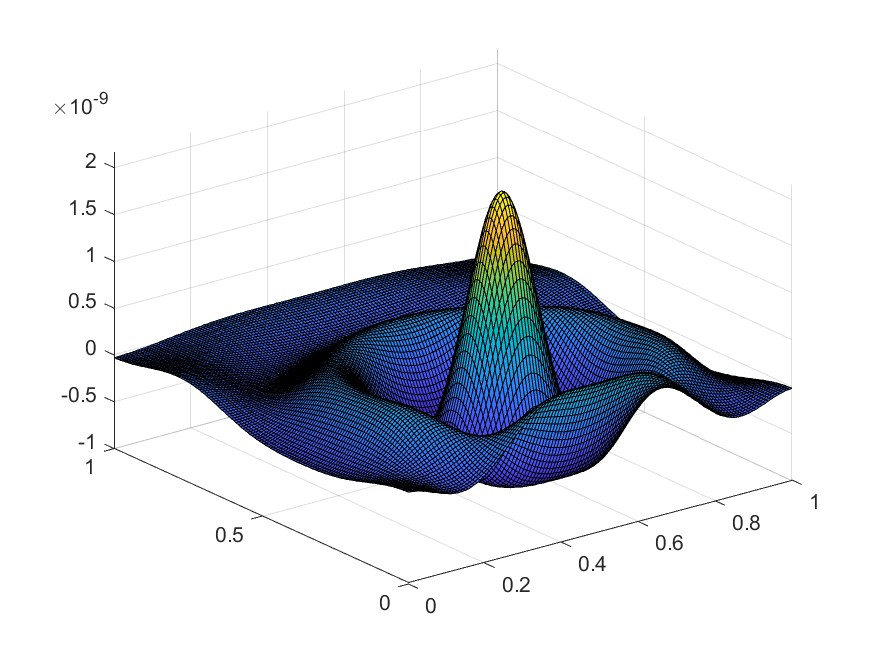}
  \caption{$t=0.2$}
\end{subfigure}
\begin{subfigure}{.24\textwidth}
  \centering
  \includegraphics[width=\linewidth]{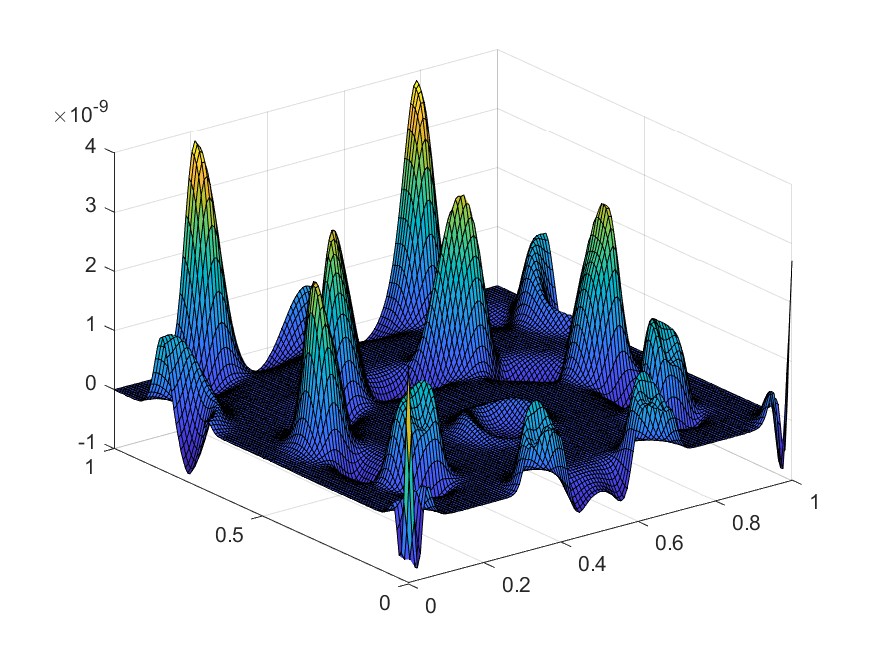}
  \caption{$t=0.5$}
\end{subfigure}
\begin{subfigure}{.24\textwidth}
  \centering
  \includegraphics[width=\linewidth]{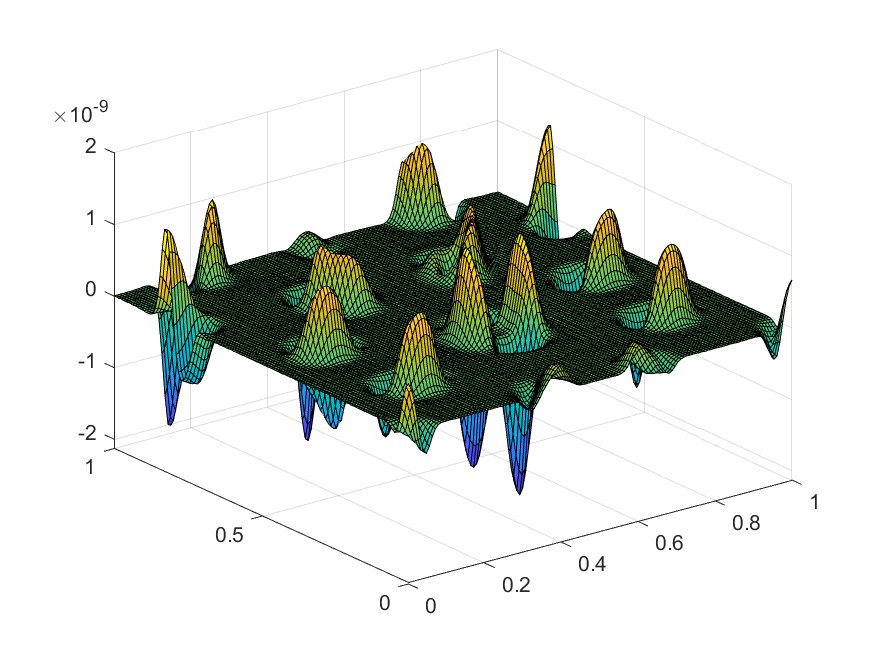}
  \caption{$t=1.0$}
\end{subfigure}
\begin{subfigure}{.24\textwidth}
  \centering
  \includegraphics[width=\linewidth]{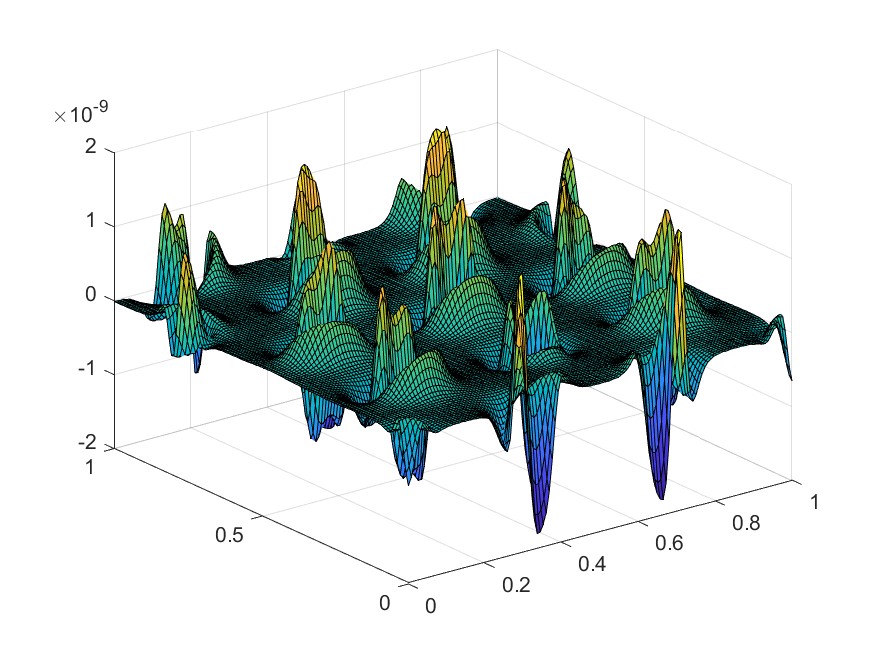}
  \caption{$t=2.0$}
\end{subfigure}
\caption{$\log-$residual decay of $U$ \& plots of residual $\mathrm{Res}(U)$ at different time stages $t=0.2,0.5,1.0, 2.0$.}\label{Schnakenberg 2}
\end{figure}

In this example, the performance of the PDHG method is stable and the method terminates in around $30$ iterations for all $10000$ time steps. We plot the loss as well as the residual term $\mathrm{Res}(U)$ at different time stages in Figure \ref{Schnakenberg 2}.

We compare the computational speed of our PHDG method with the commonly used Newton-SOR method \cite{sherman1978newton, merriman1994motion}. We fix all the parameters the same for both methods, typically, we set the termination threshold for both methods to be $\delta = 10^{-7}$. We solve the PDE system on $[0, 1]$ with $N_t=5000$ and $N_x = 128$. The time cost for the Newton-SOR method is 8122.81s, while the time cost for the PDHG method is 1121.81s.

\subsubsection{Wolf-deer model}
At last, let us consider an equation system describing the evolution of predator (wolves) and prey (deer) distributions in an ecology system \cite{murray2001mathematical, gallouet2019unbalanced}. The PDE system is defined on the region $\Omega=[-L, L]^2$ and takes the following form,
\begin{align}
  \frac{\partial \rho_1}{\partial t} & = D\Delta \rho_1 + \nabla\cdot(\rho_1\nabla \mathcal{V}_1(\rho_1, \rho_2)) + A\rho_1(1-\rho_1) - B\frac{\rho_1\rho_2}{1+\rho_1}; \label{wd 1}\\ 
  \frac{\partial \rho_2}{\partial t} & = D\Delta \rho_2 + \nabla\cdot(\rho_2\nabla \mathcal{V}_2(\rho_1, \rho_2)) + B\frac{\rho_1\rho_2}{1+\rho_1} - C\rho_2. \label{wd 2}
\end{align}

Here we set $D = \frac{1}{2}$, $A = 5$, $B = 35$, $C=\frac{5}{2}$. We also define the interacting potentials $\mathcal{V}_1, \mathcal{V}_2$ as
\[ \mathcal{V}_1(\rho_a, \rho_b)(\cdot) = 
V*\rho_1 - V*\rho_2; \quad \mathcal{V}_2(\rho_a, \rho_b)(\cdot) = V*\rho_1 + V*\rho_2.   \]
Here the convolution is defined as $V*\rho(x,y)=\iint_{\Omega\times\Omega} V((x,y)-(x',y'))\rho(x',y')~dx'dy'$ with potential $V(x,y)=\frac{x^2+y^2}{2}.$

We choose Neumann boundary condition for both $\rho_1$ and $\rho_2$, and set the initial condition
\begin{equation}
  \rho_i(x, 0) = \frac{1}{\pi}\left(\frac{\pi}{2} + \mathrm{arctan}\left(\frac{R^2 - |X-\vec{\mu}_i|^2}{\epsilon}\right)\right),~i=1,2,  \label{init wd}
\end{equation}
where $\vec{\mu}_1 = (\frac{3}{2}, \frac{3}{2}),~\vec{\mu}_2 = (-\frac{3}{2},  - \frac{3}{2}),$ $R = 1$ and $\epsilon = 0.1$.

In system \eqref{wd 1}, \eqref{wd 2}, $\rho_1$ represents the distribution of deer, and $\rho_2$ stands for the distribution of wolf.
In addition to the diffusion and reaction terms affecting $\rho_a, \rho_2$, the PDE system \eqref{wd 1}, \eqref{wd 2} contain non-local drift terms $\nabla \mathcal{V}_1(\rho_a, \rho_2), \nabla \mathcal{V}_2(\rho_1, \rho_2)$ that depict the interactions among the individuals of wolves and deer: The deer are attracting each other to dodge wolves' predation, while the wolves are gathering together to chase the flock of deer.

Suppose we discretize each side of $\Omega$ into $N_x-1$ equal subintervals, we denote the numerical solutions at the $n-$th time step as $\rho_1^{n}, \rho_2^{n}\in\mathbb{R}^{N_x^2}$. We consider the following implicit, central-difference scheme,
{\small
\begin{align*}
 \frac{\rho_{1}^{n+1} - \rho_{1}^n}{h_t} - 
 D\NLap\rho_1^{n+1} & - \Dx^\top (\barrhoax \odot \Dx (K\rho_1^{n+1} - K\rho_2^{n+1}) )\\
 & - \Dy^\top (\barrhoay \odot \Dy (K\rho_1^{n+1} - K\rho_2^{n+1}) ) + R_1(\rho_1^{n+1}, \rho_2^{n+1}) = 0;\\
 \frac{\rho_2^{n+1} - \rho_2^n}{h_t} - D\NLap\rho_2^{n+1} & - \Dx^\top (\barrhobx \odot \Dx (K\rho_1^{n+1} + K\rho_2^{n+1}))\\ 
 & - \Dy^\top (\barrhoby \odot \Dy (K\rho_1^{n+1} + K\rho_2^{n+1}))  + R_2(\rho_1^{n+1}, \rho_2^{n+1}) = 0.
\end{align*}}
Here $D_x$ is an $(N_x+1)N_x\times N_x^2$ matrix which can be treated as the discrete gradient with respect to $x$, i.e., for any $u\in\mathbb{R}^{N_x^2}$, $(\Dx u)_{(i, j+\frac{1}{2})}$ equals $\frac{u_{i,j+1} - u_{i,j}}{h_x}$ for $1 \leq i \leq N_x$, $1\leq j\leq N_x-1$; for $j=0, N_x,$ we define $(\Dx u)_{(i, \frac{1}{2})} = \frac{u_{i, 2} - u_{i, 1}}{h_x}$ and $(\Dx u)_{(i, N_x+\frac{1}{2})} = \frac{u_{i, N_x} - u_{i, N_x - 1}}{h_x}.$ $\Dy$ can also be defined in a similar way. 

The notation $\barrhoax\in\mathbb{R}^{(N_x+1)N_x}$ denotes the average value of $\rho^{n+1}_1$ at midpoints, i.e., $(\barrhoax)_{(i,j+\frac{1}{2})} = \frac{\rho_{1, (i,j)}^{n+1} + \rho_{1, (i,j+1)}^{n+1}}{2}$ for $1\leq i\leq N_x, ~1\leq j \leq N_x-1$; for $j=0, N_x$, we define $\barrhoax_{(i, \frac{1}{2})} = \rho^{n+1}_{1, i, 1}$ and $\barrhoax_{(i, N_x+\frac{1}{2})} = \rho^{n+1}_{i, N_x}.$ $\barrhoay, \barrhobx, \barrhoby$ can be defined in the similar way.

 $K$ is an $N_x^2\times N_x^2$ matrix used for approximating the convolution $V*\rho_1, V*\rho_2$, to precisely, for any $u\in \mathbb{R}^{N_x^2}$ defined on the mesh grid of $\Omega$, $Ku$ is defined as 
 $$(Ku)_{(i,j)} = \sum_{1\leq k,l \leq N_x} h_x^2 V(v_{i,j} - v_{k,l}) u_{k,l} = \sum_{1\leq k,l \leq N_x} \frac{h_x^4}{2} ((i-k)^2 + (j-l)^2) u_{k, l}.$$
The above discrete convolution can be reduced to Toeplitz matrix-vector multiplication computation, which can be efficiently computed by FFT algorithm \cite{strang1986proposal}.

Furthermore, for $\rho_1, \rho_2\in\mathbb{R}^{N_x^2}$ the reaction terms are defined as $R_1(\rho_1,\rho_2)=A\rho_1\odot(1-\rho_1) - B \frac{\rho_1}{1+\rho_1}\odot\rho_2$; $R_2(\rho_1, \rho_2) = B\frac{\rho_1}{1+\rho_1}\odot \rho_2 - C\rho_2$.

Given the above discrete scheme of the PDE system, similar to the discussion made in the Schnakenberg model, our purpose is to solve two nonlinear equations $F_{\rho_1}(\rho_1, \rho_2)=0,~ F_{\rho_2}(\rho_1, \rho_2)=0$ at each time step $n$. We apply two dual variables $P,Q$, and compose the corresponding PDHG dynamic for solving the two equations. According to the previous discussion, one can verify that each PDHG step can be computed within $O(N_x^2\log N_x)$ complexity. This guarantees the efficiency of the computation. To keep the discussion concise, we omit the exact formulas for the PDHG dynamic here.

\begin{figure}[htb!]
\begin{subfigure}{.162\textwidth}
  \centering
  \includegraphics[width=\linewidth]{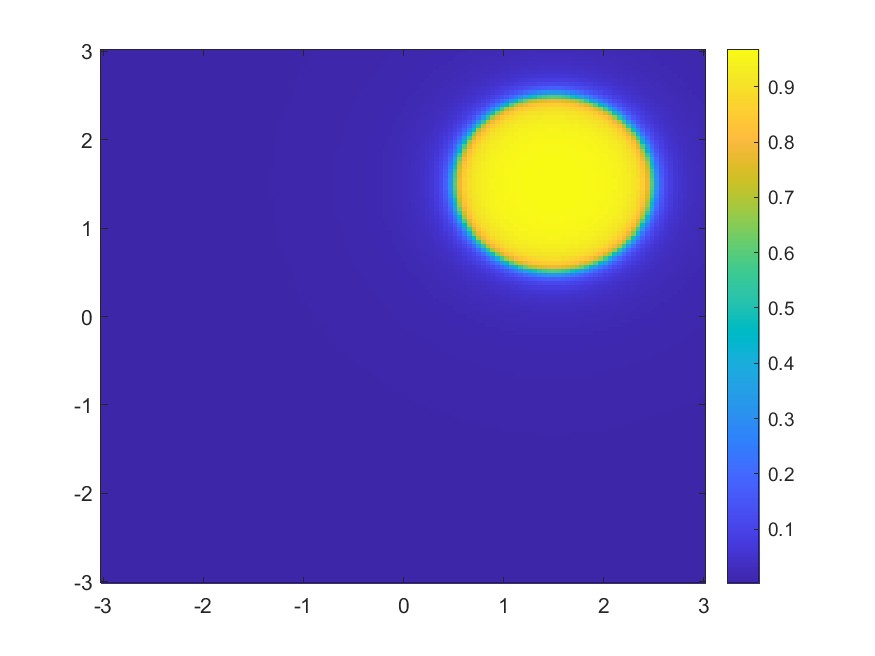}
\end{subfigure}
\begin{subfigure}{.162\textwidth}
  \centering
  \includegraphics[width=\linewidth]{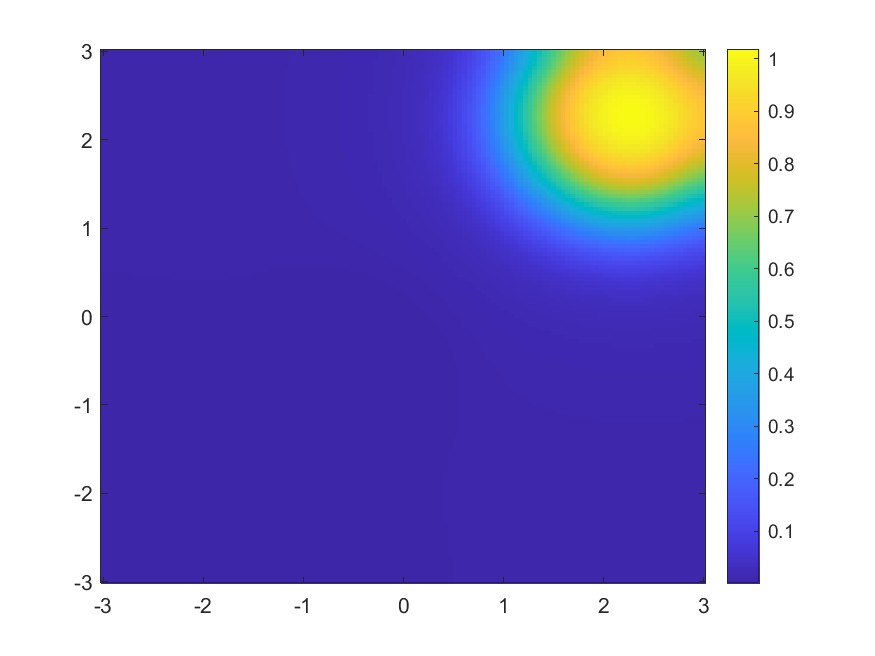}
\end{subfigure}
\begin{subfigure}{.162\textwidth}
  \centering
  \includegraphics[width=\linewidth]{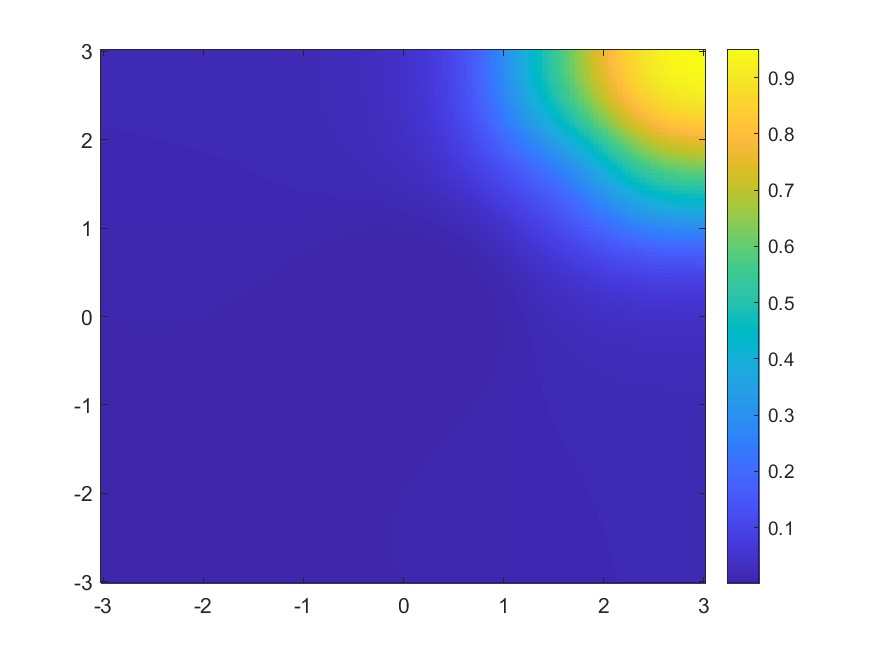}
\end{subfigure}
\begin{subfigure}{.162\textwidth}
  \centering
  \includegraphics[width=\linewidth]{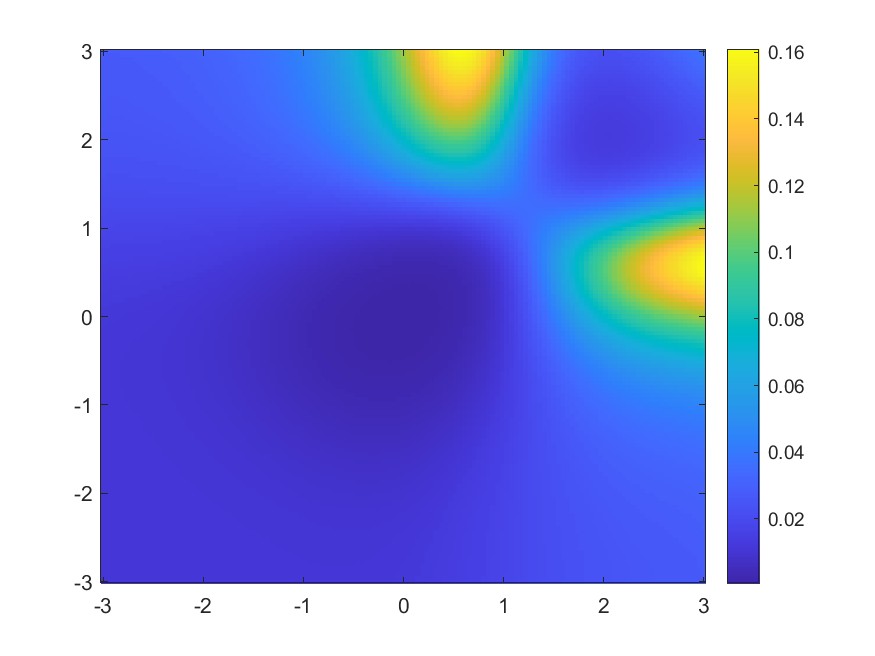}
\end{subfigure}
\begin{subfigure}{.162\textwidth}
  \centering
  \includegraphics[width=\linewidth]{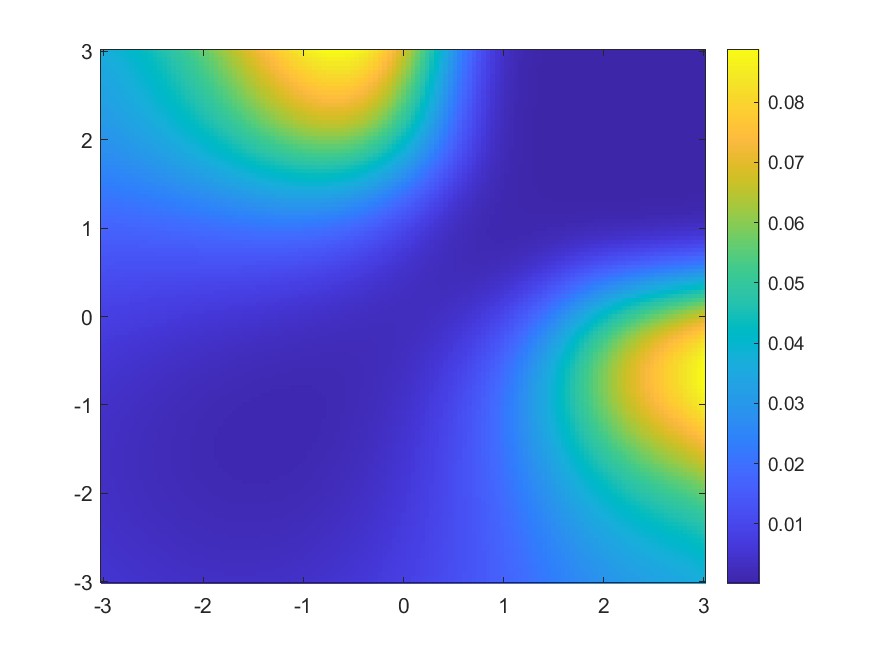}
\end{subfigure}
\begin{subfigure}{.162\textwidth}
  \centering
  \includegraphics[width=\linewidth]{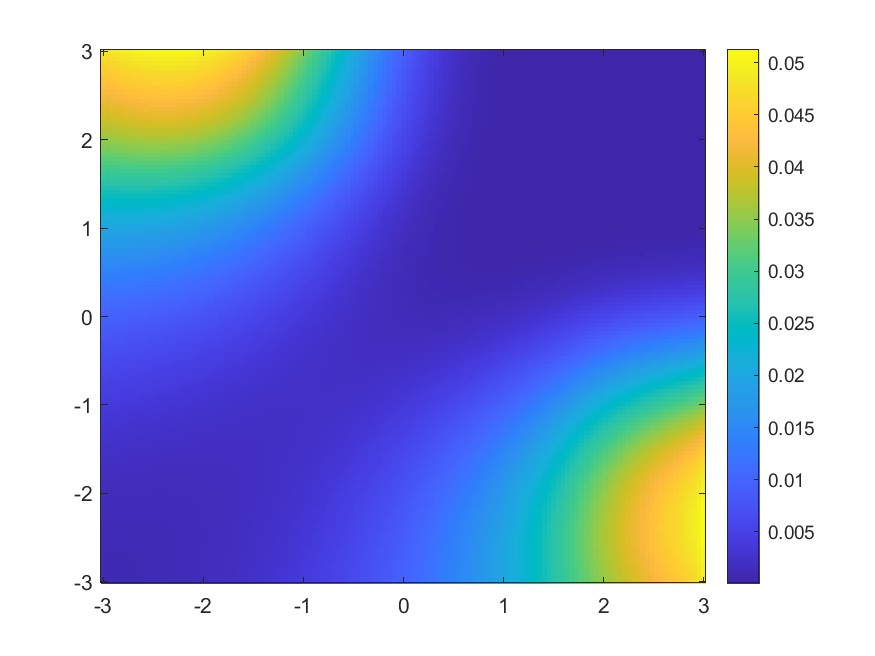}
\end{subfigure}
\begin{subfigure}{.162\textwidth}
  \centering
  \includegraphics[width=\linewidth]{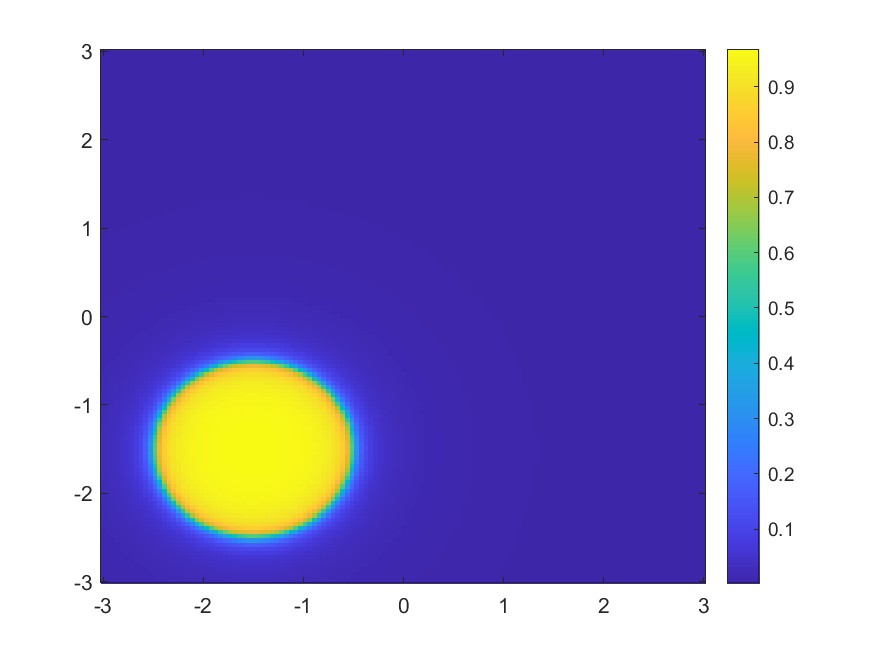}
  \caption{$t=0.0$}
\end{subfigure}
\begin{subfigure}{.162\textwidth}
  \centering
  \includegraphics[width=\linewidth]{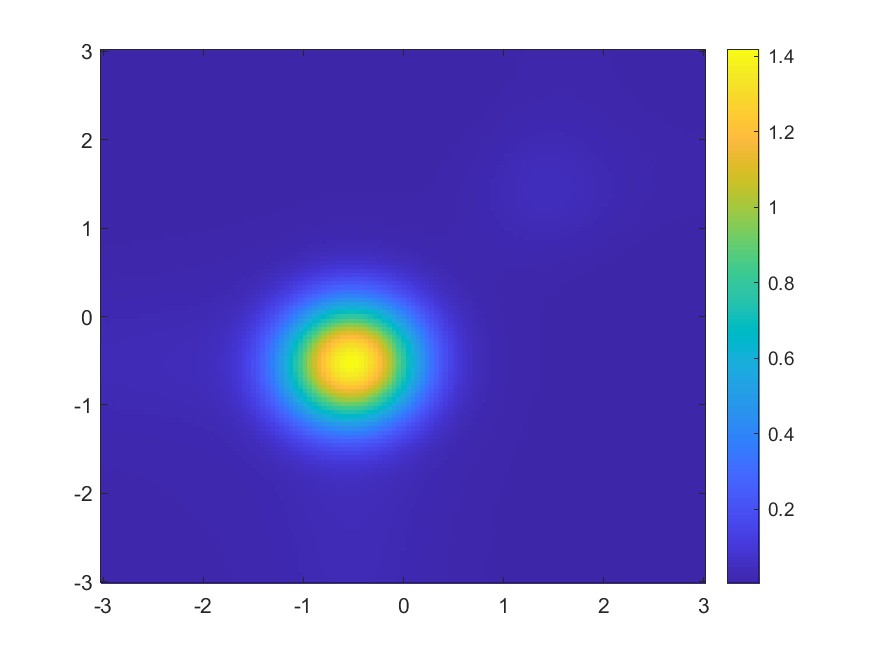}
  \caption{$t=0.2$}
\end{subfigure}
\begin{subfigure}{.162\textwidth}
  \centering
  \includegraphics[width=\linewidth]{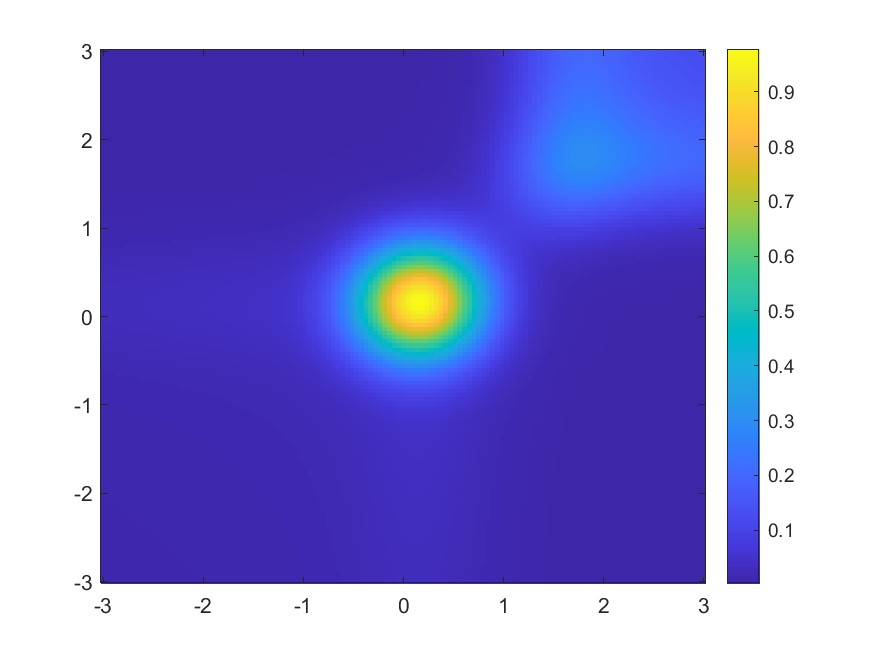}
  \caption{$t=0.4$}
\end{subfigure}
\begin{subfigure}{.162\textwidth}
  \centering
  \includegraphics[width=\linewidth]{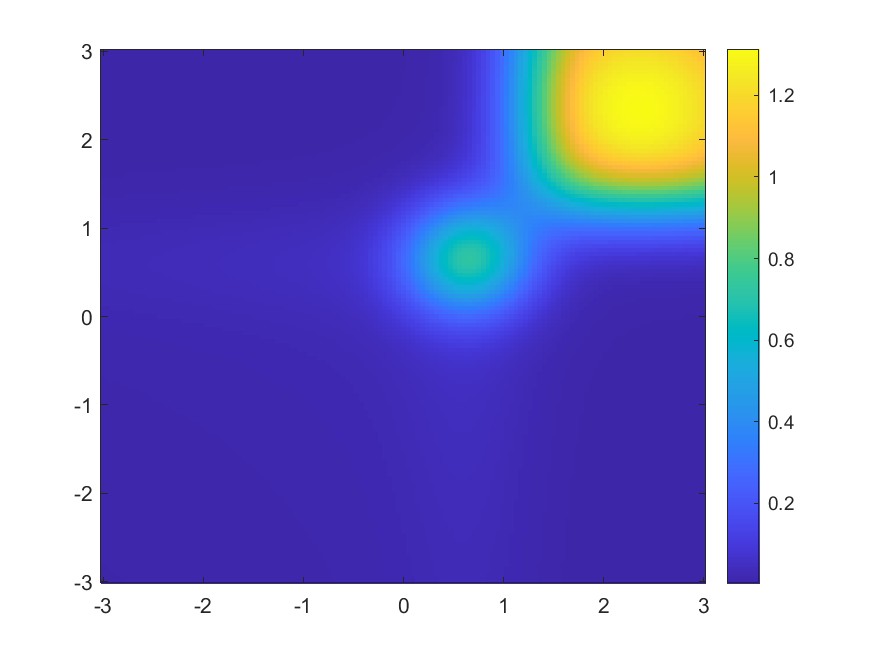}
  \caption{$t=0.585$}
\end{subfigure}
\begin{subfigure}{.162\textwidth}
  \centering
  \includegraphics[width=\linewidth]{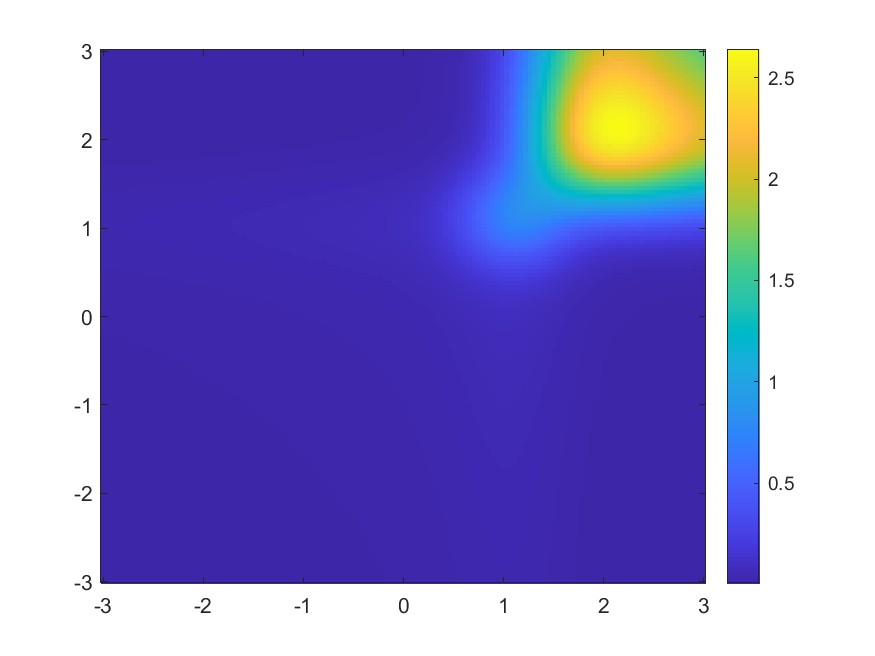}
  \caption{$t=0.759$}
\end{subfigure}
\begin{subfigure}{.162\textwidth}
  \centering
  \includegraphics[width=\linewidth]{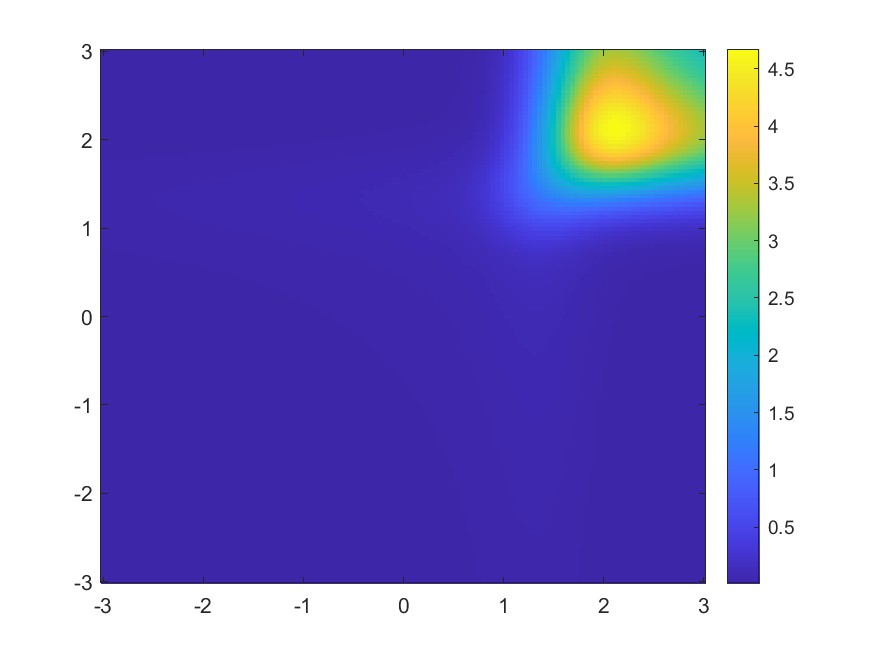}
  \caption{$t=0.866$}
\end{subfigure}
\caption{Numerical solution of $\rho_1$ (upper row), and $\rho_2$ (lower row) at different time stages with initial condition \eqref{init wd}.}\label{wd fig 1}
\end{figure}

In our implementation, we set $L = 3$. We solve the equation system \eqref{wd 1}, \eqref{wd 2} on time interval $[0, 1]$. We set $N_x = 128$, $h_x = 3/64$. We practice the method of adaptive time stepsize $h_t$ in this example. We set both our initial time stepsize $h_t$ and maximum stepsize $h_t^0$ equal to $1/500$ with shrinkage/enlarge
coefficient $\eta$ as $0.75$. The thresholding iteration numbers of shrinking and enlarging $h_t$ are set to be $N^*=100$, and $N_*=20$. For the PDHG iteration, we set stepsize $\tau_u=\tau_p = 0.95$, and pick the threshold $\delta = 5\times 10^{-6}$. We present the numerical results in Figure \ref{wd fig 1}.

\begin{figure}[htb!]
\begin{subfigure}{.19\textwidth}
  \centering
  \includegraphics[width=\linewidth]{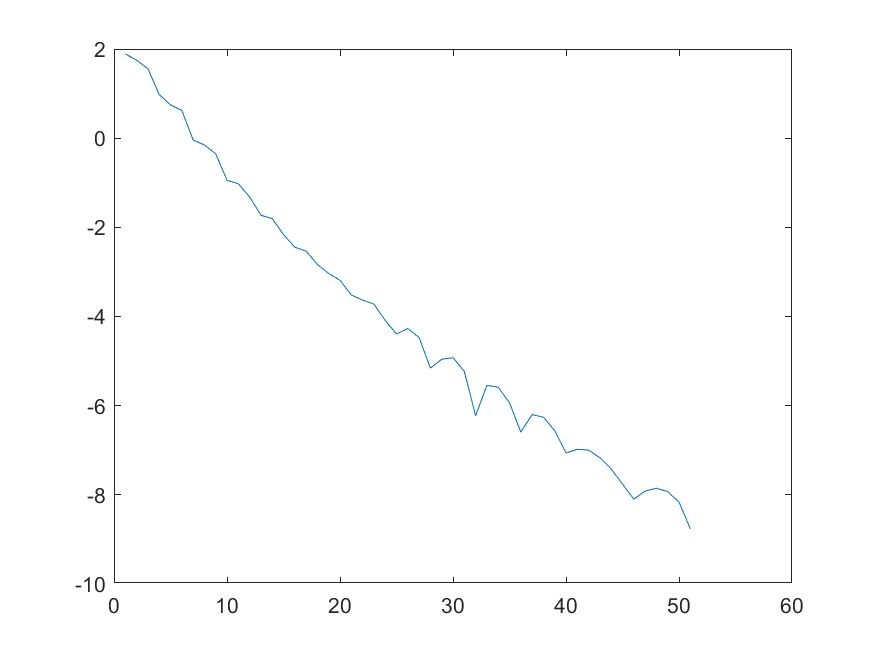}
  \caption{$t=0.2$}
\end{subfigure}
\begin{subfigure}{.19\textwidth}
  \centering
  \includegraphics[width=\linewidth]{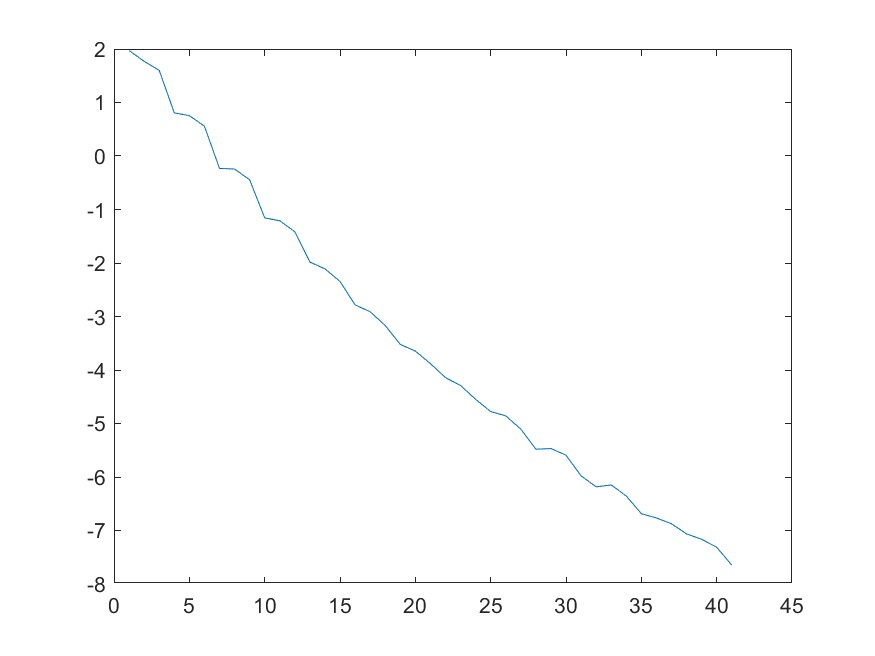}
  \caption{$t=0.4$}
\end{subfigure}
\begin{subfigure}{.19\textwidth}
  \centering
  \includegraphics[width=\linewidth]{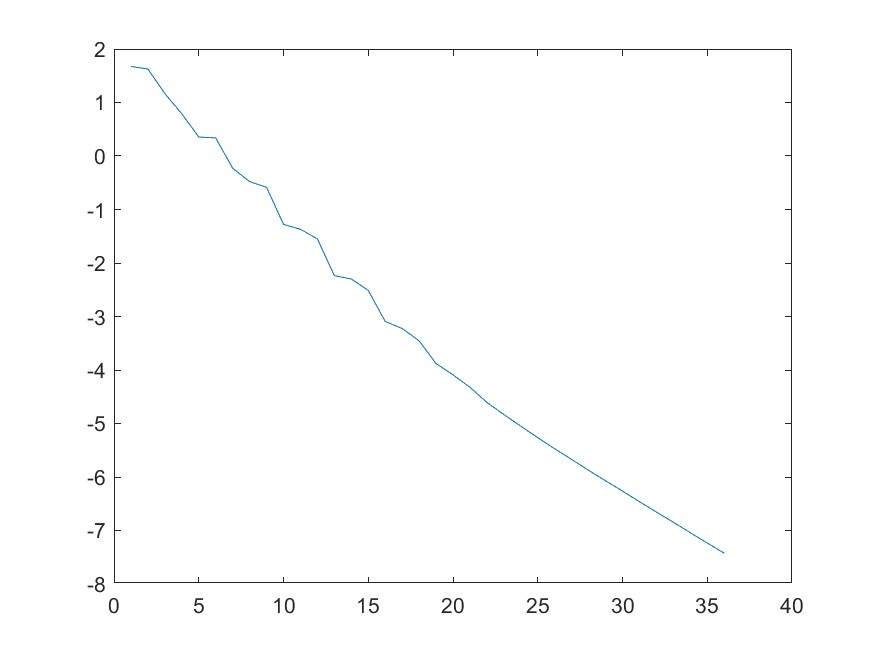}
  \caption{$t=0.585$}
\end{subfigure}
\begin{subfigure}{.19\textwidth}
  \centering
  \includegraphics[width=\linewidth]{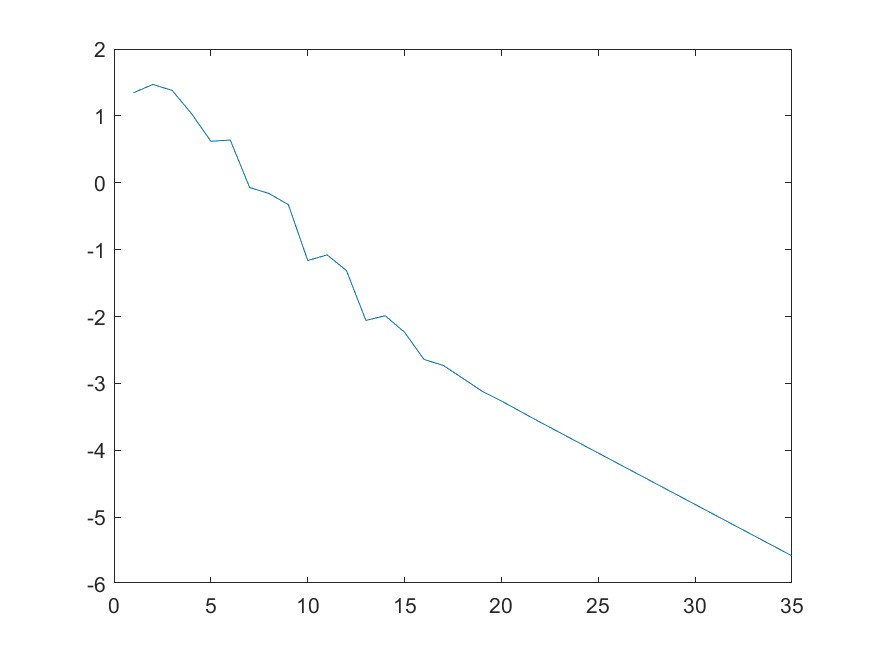}
  \caption{$t=0.759$}
\end{subfigure}
\begin{subfigure}{.19\textwidth}
  \centering
  \includegraphics[width=\linewidth]{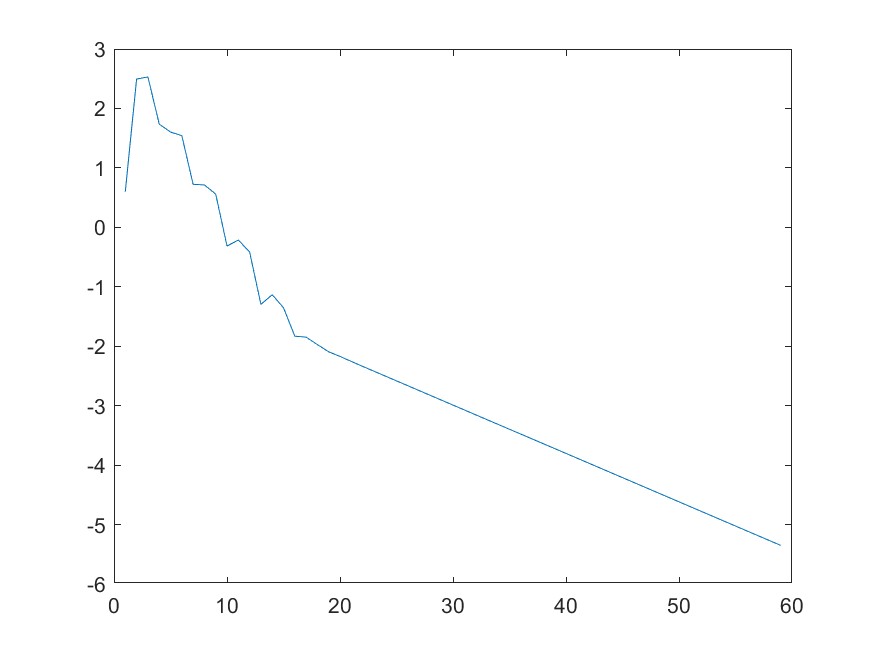}
  \caption{$t=0.997$}
\end{subfigure}
\caption{$\log_{10}\mathrm{Res}(U_n)$ vs PDHG iteration number $n$ at different time stages.}\label{wd fig 2}
\end{figure}

The linear convergence of the residual term $\|\mathrm{Res}(U)\|_2$ is reflected from the residual decay plots in Figure \ref{wd fig 2}.

\begin{figure}[htb!]
    \begin{subfigure}{.45\linewidth}
      \centering
      \includegraphics[width=\linewidth]{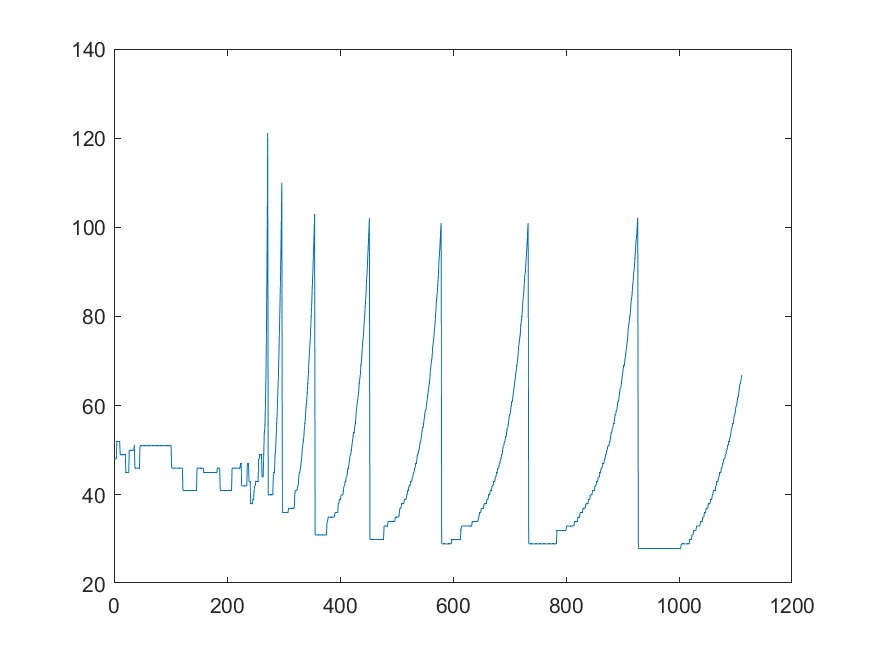}
      \caption{Number of PDHG iterations at each time step.}
    \end{subfigure}
    \begin{subfigure}{.45\linewidth}
      \centering
      \includegraphics[width=\linewidth]{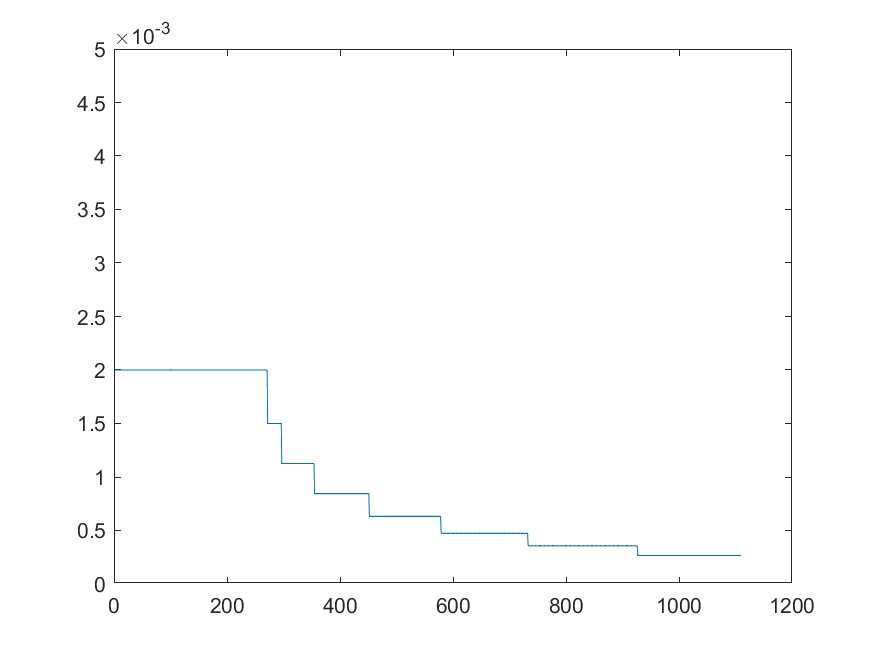}
      \caption{Time stepsize $h_t$ used at each time step, initial $h_t=0.002$, terminal $h_t\approx 0.00027$.}
    \end{subfigure}
    \caption{Changes in number of PDHG iterations \& time stepsize $h_t$.}\label{wd fig 3}
\end{figure}

The changes in PDHG iterations at each time stepsize as well as the changes in time step size $h_t$ are demonstrated via Figure \ref{wd fig 3}. As reflected from the plots, in this example, we are gradually shrinking the time stepsize $h_t$ as the accumulated time increases to guarantee the computational efficiency of the PDHG method. Our method takes a total of $1106$ time to complete the computation. the algorithm experiences $7$ stepsize shrinkage among our computation. The initial $h_t$ is set as $0.002$ while when we finish the computation, $h_t = 0.00027.$ 

\section{Conclusion and Future Study}\label{conclude }
This research proposes an iterative method as a convenient but efficient gadget for solving the implicit (or semi-implicit) numerical schemes arising in time-evolution PDEs, especially the reaction-diffusion type equations.  Our method recasts the nonlinear equation from the discrete numerical scheme at each time step as a min-max saddle point problem and applies the Primal-Dual Hybrid Gradient method. The algorithm can flexibly fit into various numerical schemes, {such as semi-implicit and fully implicit schemes, etc.} Furthermore, the method is easy to implement since it gets rid of the computation of large-scale linear systems involving Jacobian matrices, which are usually required by Newton's methods. The performance of our method on accuracy and efficiency is satisfying and is comparable to the commonly used Newton-type methods. This has been verified by the numerical examples presented in this paper.

There are three main future research directions of our work. We summarize them below:
\begin{itemize}
    \item Conduct theoretical analysis on the convergence of the PDHG method for nonlinear RD equations. We are interested in the necessary condition on $h_t, h_x, \tau_u, \tau_p$ that can guarantee the convergence of our method for certain types of RD equations. 
    \item Generalize the method to nonlinear time-evolution equations, especially the advection-reaction-diffusion dynamics from GENERIC (General Equation for Non-Equilibrium Reversible-Irreversible Coupling)  \cite{duong2023non, grmela1997dynamics, ottinger1997dynamics}.
    \item  Apply the method to high-dimensional time-evolution PDEs by leveraging deep learning techniques and PDHG algorithms. 
\end{itemize}
\noindent\textbf{Acknowledgement:} S. Liu and S. Osher's work was partly supported by AFOSR MURI FP 9550-18-1-502 and ONR
grants:  N00014-20-1-2093 and N00014-20-1-2787. W. Li's work was supported by AFOSR MURI FP 9550-18-1-502, AFOSR YIP award No. FA9550-23-1-0087, and NSF RTG: 2038080.

\appendix
\section{Proof of Theorem \ref{thm linear convergence}}\label{thm pf}
\begin{proof}
It is not hard to verify that the dynamic \eqref{illustrate pdhg dynamic 1}, \eqref{illustrate pdhg dynamic 2} and \eqref{illustate pdhg dynamic 3 simplified} can be formulated as
\begin{equation}
  \left[\begin{array}{c}
     U_{n+1} \\
     P_{n+1}
  \end{array}\right] = \left[\begin{array}{cc}
     I - 2\tau_u\tau_p A^\top A   &  -\tau_u A^\top \\
     \tau_p A  &   I
  \end{array}\right]  \left[\begin{array}{c}
    U_n \\
    P_n
  \end{array}\right] + \left[\begin{array}{c}
       2\tau_u\tau_p A^\top b \\
       - \tau_p b
  \end{array}\right], \quad n\geq 0\label{iterequ}
\end{equation}
This equation admits a unique fixed point $X_* = (U_*, P_*)=(A^{-1}b, 0)$. We denote $X_n = [U_n^\top, P_n^\top]^\top$ and the above recurrence equation as $X_{n+1} = M X_n + y$ (or equivalently, $X_{n+1}-X_*=M(X_n-X_*)$) for shorthand. Suppose $A$ has spectral decomposition $A = Q \Lambda Q^\top$, then $M$ is decomposed as 
\begin{equation*}
 M = \left[\begin{array}{cc}
     Q  &  \\
      &   Q
 \end{array}\right] \left[\begin{array}{cc}
    I - 2\tau_u\tau_p\Lambda^2  &  -\tau_u\Lambda \\
     \tau_p \Lambda  &  I
 \end{array}\right]\left[\begin{array}{cc}
     Q  &  \\
      &   Q
 \end{array}\right]^\top
\end{equation*}
We denote $N_A$ as the size of $A$. The middle matrix is composed of four $N_A\times N_A$ diagonal matrices, by rearranging the rows and columns of it, one can show that it is orthogonally equivalent to the block diagonal matrix $\Sigma = \mathrm{diag}(D_1,D_2,...,D_{N_A})$ where each $D_k = \left[\begin{array}{cc}
    1-2\tau_u\tau_p\lambda_k^2 & -\tau_u\lambda_k \\
    \tau_p\lambda_k   &    1 
\end{array}\right].$ We now analyze the spectral radius 
$\rho(M)$ of $M$, which equals $\rho(\Sigma) = \underset{1\leq k\leq N_A}{\max}\{\rho(D_k)\}$. We can calculate $$\rho(D_k) = \max\{|\tau_u\tau_p\lambda_k^2 - 1 \pm \sqrt{(\tau_u\tau_p)^2\lambda_k^4-\tau_u\tau_p\lambda_k^2}|\}.$$ 
We denote $f(t) = \max\{|t-1+\sqrt{t^2-t}|, |t-1 - \sqrt{t^2-t}|\},$ with $t>0$. One can directly compute $f(t) = \begin{cases}
    \sqrt{1-t} \quad 0 < t \leq 1\\
    t-1+\sqrt{t^2-t} \quad t > 1
\end{cases}.$ 
Thus we have $\rho(M) = \underset{1\leq k\leq N_A}{\max} \{f(\tau_u\tau_p\lambda_k^2)\}.$ Then $f$ is decreasing on $[0, 1]$ and increasing on $[1,\infty)$ with $f(0)=f(\frac{4}{3})=1$. We know that the convergence of \eqref{iterequ} is guaranteed if and only if $\rho(M)<1$. This is equivalent to $\tau_u\tau_p\lambda_{\max}^2 \leq \frac{4}{3}$, which yields $\tau_u\tau_p \leq \frac{4}{3\lambda_{\max}^2}$. 

Furthermore, $\rho(M)$ is the convergence rate of the dynamic, i.e., $$\|X_n-X_*\|_2 \leq \rho(M)^n \|X_0-X_*\|_2 $$ 
To evaluate the optimal convergence rate, we compute the minimum value of $\rho(M)$ w.r.t. stepsizes $\tau_u, \tau_p$. Suppose we require $\tau_u\tau_p\leq \frac{4}{3\lambda_{\max}}$ to guarantee convergence, if we denote $\eta = \tau_u\tau_p\lambda^2_{\max}$, then $\rho(M) = \underset{1 \leq k \leq N_A}{\max}\left\{f(\frac{\lambda_k^2}{\lambda_{\max}^2}\eta)\right\}$ for $\eta \in (0, \frac{4}{3}).$  The minimum value of $\rho(M)$ will be attained at a unique $\eta = \eta_*\in[1, \frac{4}{3})$ such that $f(\frac{\eta}{\kappa^2})=f(\eta)$. I.e., $\eta^*$ is the solution of 
\begin{equation}
  \sqrt{1-\frac{ \eta}{\kappa^2}} = \eta - 1 + \sqrt{\eta^2 - \eta},\quad \textrm{on } [1, \frac{4}{3}).   \label{equation eta}
\end{equation}
Thus, the optimal convergence rate $\gamma_* = \sqrt{1-\frac{\eta_*}{\kappa^2}}$ and it is achieved when $\tau_u, \tau_p$ satisfy $\tau_u\tau_p = \frac{\eta^*}{\lambda_{\mathrm{max}}^2}$.

\end{proof}

\begin{remark}\label{formula of gamma}
  The equation \eqref{equation eta} can be reduced to a quadratic equation. And it admits a unique solution $\eta_*$ on $[1, \frac{4}{3})$, $\eta_*$ takes the following form 
  \begin{equation*}
     \eta_* = \frac{2\kappa^2}{\left(\frac{3}{4}\kappa^2+\frac{3}{2}-\frac{1}{4\kappa^2}\right) + \frac{\kappa-1}{2\kappa}\sqrt{(\kappa-1)(3\kappa+1)}\sqrt{\frac{3}{4}\kappa^2+\frac{3}{2}+2\kappa-\frac{1}{4\kappa^2}}}.
  \end{equation*}
  The optimal convergence rate $\gamma_* = \sqrt{1-\frac{\eta_*}{\kappa^2}}$ takes the explicit form
  \begin{equation*}
    \gamma_* = \left( 1- \frac{2}{\left(\frac{3}{4}\kappa^2+\frac{3}{2}-\frac{1}{4\kappa^2}\right) + \frac{\kappa-1}{2\kappa}\sqrt{(\kappa-1)(3\kappa+1)}\sqrt{\frac{3}{4}\kappa^2+\frac{3}{2}+2\kappa-\frac{1}{4\kappa^2}}} \right)^{\frac{1}{2}}.
  \end{equation*}
  Notice that $\gamma_*\approx (1-\frac{4}{3\kappa^2})^{\frac{1}{2}}$ when the conditional number $\kappa$ is very large; and $\gamma_*$ will approach $0$ as condition number $\kappa$ approaches $1$. This motivates the preconditioning technique of our method.
\end{remark}

\bibliographystyle{plain}
\bibliography{references}

\end{document}